\documentclass[12pt]{article}
\usepackage{amscd}
\usepackage{amsthm}
\usepackage{extarrows}
\usepackage{latexsym}
\usepackage{mathrsfs}
\usepackage{pst-all,multido,ifthen}
\usepackage{amsmath}
\usepackage{amssymb}
\usepackage[all]{xy}
\usepackage[english]{babel}
\usepackage{mathtools}
\usepackage{upgreek}

\newtheorem{theorem}{Theorem}
\newtheorem{definitions}{Definitions}
\newtheorem{proposition}{Proposition}
\newtheorem{corollary}{Corollary}
\newtheorem{lemma}{Lemma}

\newtheorem{fact}{Fact}

\theoremstyle{remark}
\newtheorem{remark}{Remark}
\def\Spec{\text{Spec}}
\def\Spf{\text{Spf}}

\def\Ker{\text{Ker}}

\def\Lie{\text{Lie}}
\def\endproof{$\hfill \square$}
\def\proof{\par\noindent {{\it Proof:}}\enspace}

\usepackage{scalerel}
\usepackage{stackengine,wasysym}
\stackMath
\usepackage{tikz}

\usepackage{pgffor}
\makeatletter
\newcommand{\medoplus}{\mathbin{\mathpalette\make@small\oplus}}
\newcommand{\medotimes}{\mathbin{\mathpalette\make@small\otimes}}

\newcommand{\make@small}[2]{%
 \vcenter{\hbox{%
 \scalebox{1.6}{$\m@th#1#2$}%
 }}%
}
\makeatother 

\begin{document}
\title{Good reductions of Shimura varieties of Hodge type in arbitrary unramified mixed characteristic. Part I}
\author{Adrian Vasiu}
\maketitle

\centerline{final version, to appear in {\it Mathematische Nachrichten}} 
\centerline{(most alignment issues kept loose to match with the layout of the journal)}

\bigskip\noindent
{\bf ABSTRACT.} We prove the existence of good smooth integral models of Shimura varieties of Hodge type in arbitrary unramified mixed characteristic $(0,p)$. As a first application we provide a smooth solution (answer) to a conjecture (question) of Langlands for Shimura varieties of Hodge type. As a second application we prove the existence in arbitrary unramified mixed characteristic $(0,p)$ of integral canonical models of projective Shimura varieties of Hodge type with respect to h--hyperspecial subgroups as pro-\'etale covers of N\'eron models; this forms progress towards the proof of conjectures of Milne and Reimann. Though the second application was known before in some cases, its proof is new and more of a principle.  

\bigskip\noindent
{\bf KEY WORDS:} abelian scheme, affine group scheme, Barsotti--Tate group, deformation theory, $F$-crystal, Hodge cycle, integral model, and Shimura variety

\bigskip\noindent
{\bf MSC 2010:} 11G10, 11G18, 14F30, 14G35, 14G40, 14K10, 14J10

\section{Introduction}
Let $p\in\mathbb N$ be a prime. Let $\mathbb Z_{(p)}$ be the localization of $\mathbb Z$ at its prime ideal $(p)$. Let $r\in\mathbb N^*$. Let $N\ge 3$ be a natural number relatively prime to $p$. Let $\mathcal A_{r,1,N}$ be the {\it Mumford moduli scheme} over $\mathbb Z_{(p)}$ that parameterizes isomorphism classes of {\it principally polarized abelian schemes} over $\mathbb Z_{(p)}$-schemes of relative dimension $r$ and endowed with a symplectic similitude level-$N$ structure (cf. [43, Thms. 7.9 and 7.10] applied to symplectic similitude level structures instead of simply level structures). 

\subsection{Basic properties} 
The $\mathbb Z_{(p)}$-schemes $\mathcal A_{r,1,N}$ have the following four properties:

\medskip\noindent
{\bf (i)} {\it They are smooth and quasi-projective.}

\smallskip\noindent
{\bf (ii)} {\it If $N_1\in N\mathbb N\setminus p\mathbb N$, then the natural level-reduction morphism $\mathcal A_{r,1,N_1}\rightarrow\mathcal A_{r,1,N}$ is finite, \'etale, and surjective. Thus the projective limit 
$$\mathcal M_r:=\projlim_{N\ge 3,(N,p)=1} \mathcal A_{r,1,N}$$ 
exists and is a regular, quasi-compact, formally smooth $\mathbb Z_{(p)}$-scheme.}

\smallskip\noindent
{\bf (iii)} {\it If $Z$ is a regular, formally smooth scheme over $\mathbb Z_{(p)}$, then each morphism $Z_{\mathbb Q}\rightarrow \mathcal M_{r,\mathbb Q}$ extends uniquely to a morphism $Z\rightarrow \mathcal M_r$.}

\smallskip\noindent
{\bf (iv)} {\it Their geometric special fibers have level $m$ {\it stratifications} ($m\in\mathbb N^*$) enjoying many properties: strata are regular, quasi-affine, equidimensional of dimensions given by explicit formulas, etc.}

\medskip
Property (i) is checked in loc. cit., cf. also Serre's Lemma of [42, Ch. IV, Sect. 21, Thm. 5]. Property (ii) is well-known. Property (iii) is implied by the fact that each abelian scheme over $Z_{\mathbb Q}$ that has level-$N$ structure for all $N\in\mathbb N\setminus (p\mathbb N\cup\{1,2\})$, extends to an abelian scheme over $Z$ (cf. the N\'eron--Ogg--Shafarevich criterion of good reduction and the purity result [67, Cor. 5]); such an extension is unique up to a unique isomorphism (cf. [49, Ch. IX, Cor. 1.4]). Property (iv) is an application of the deformation theories for abelian varieties and for Barsotti--Tate groups (i.e., $p$-divisible groups): see [45] for the case $m=1$ and see [59, Ex. 4.5] for all $m\in\mathbb N^*$.

From Yoneda Lemma we get that the regular, formally smooth  $\mathbb Z_{(p)}$-scheme $\mathcal M_r$ is uniquely determined by its generic fibre $\mathcal M_{r,\mathbb Q}$ and by the {\it universal property} expressed by the property (iii). Thus one can view $\mathcal A_{r,1,N}$ as the {\it best} smooth integral model of $\mathcal A_{r,1,N,\mathbb Q}$ over $\mathbb Z_{(p)}$. The main goal of this paper is to generalize properties (i) to (iv) to the context of Shimura varieties of Hodge type. Thus in this paper we prove the existence and the uniqueness of {\it good} smooth integral models of Shimura varieties of Hodge type in unramified mixed characteristic $(0,p)$ and we list several main properties of them, including identifying cases when the theoretical non-smooth loci are actually empty. We emphasize from the very beginning that this paper brings no new contribution to either the study of non-smooth loci (when non-empty) or to ramified mixed characteristic $(0,p)$ situations. We will begin with notation and with a review of Shimura varieties.  

\subsection{Notation} 
Let $\mathbb S:=\text{Res}_{\mathbb C/\mathbb R} \mathbb G_{m,\mathbb C}$ be the two dimensional torus over $\mathbb R$ such that we have identifications $\mathbb S(\mathbb R)=\mathbb G_{m,\mathbb C}(\mathbb C)$ and $\mathbb S(\mathbb C)=\mathbb G_{m,\mathbb C}(\mathbb C)\times\mathbb G_{m,\mathbb C}(\mathbb C)$ with the property that the monomorphism $\mathbb R\hookrightarrow\mathbb C$ induces the map $z\rightarrow (z,\bar z)$; here $z\in\mathbb S(\mathbb R)=\mathbb G_{m,\mathbb C}(\mathbb C)$ is a non-zero complex number. 

Let $R$ be a commutative $\mathbb Z$-algebra. We recall that a group scheme $F$ over $R$ is called {\it reductive} if it is smooth and affine and its fibres are connected and have trivial unipotent radicals. Let $\text{Lie}(\natural)$ be the {\it Lie algebra} over $R$ of a smooth, closed subgroup scheme $\natural$ of $F$. The group schemes $\mathbb G_{m,R}$ and $\mathbb G_{a,R}$ are over $R$. For a free module $M$ of finite rank over $R$, let $M^{\vee}:=\text{Hom}(M,R)$ be its dual, and let $\pmb{\text{GL}}_M$ be the reductive group scheme over $R$ of linear automorphisms of $M$. A bilinear form $\psi:M\times M\rightarrow R$ on $M$ is called perfect if it defines naturally an $R$-linear  isomorphism $M\rightarrow M^{\vee}$. If $\psi$ is a perfect, alternating bilinear form on $M$ (thus the rank of $M$ is even), then $\pmb{\text{Sp}}(M,\psi)$ and $\pmb{\text{GSp}}(M,\psi)$ are viewed as reductive group schemes over $R$. 

Let $k$ be a perfect field of characteristic $p$. Let $W(k)$ be the ring of $p$-typical Witt vectors with coefficients in $k$. Always $n\in\mathbb N^*$. Let $\mathbb A_f:=\widehat{\mathbb Z}\otimes_{\mathbb Z}\mathbb Q$ be the ring of finite ad\`eles of $\mathbb Q$. Let $\mathbb A_f^{(p)}$ be the ring of finite ad\`eles of $\mathbb Q$ with the $p$-component omitted; we have $\mathbb A_f=\mathbb Q_p\times \mathbb A_f^{(p)}$. If $R\in\{\mathbb A_f,\mathbb A_f^{(p)},\mathbb Q_p\}$, then the group $F(R)$ is endowed with the coarsest topology that makes all maps $R=\mathbb G_{a,R}(R)\rightarrow F(R)$ associated to morphisms $\mathbb G_{a,R}\rightarrow F$ of $R$-schemes to be continuous; thus $F(R)$ is a totally discontinuous locally compact group. Each continuous action of a totally discontinuous locally compact group on a scheme will be in the sense of [13, Subsubsect. 2.7.1] and it will be a right action. 

\subsection{Shimura varieties} 
A {\it Shimura pair} $(G,\mathcal X)$ consists of a reductive group $G$ over $\mathbb Q$ and a $G(\mathbb R)$-conjugacy class $\mathcal X$ of homomorphisms $\mathbb S\rightarrow G_{\mathbb R}$ that satisfy the three axioms [13, (2.1.1.1) to (2.1.1.3)]: the Hodge $\mathbb Q$--structure on $\text{Lie}(G)$ defined by each $h\in \mathcal X$ is of type $\{(-1,1),(0,0),(1,-1)\}$, no simple factor of the adjoint group $G^{\text{ad}}$ of $G$ becomes compact over $\mathbb R$, and $\text{Ad}(\mathbb R)(h(i))$ is a Cartan involution of $\text{Lie}(G^{\text{ad}}_{\mathbb R})$. Here $\text{Ad}:G_{\mathbb R}\rightarrow\pmb{\text{GL}}_{\text{Lie}(G^{\text{ad}}_{\mathbb R})}$ is the adjoint representation. These axioms imply that $\mathcal X$ has a natural structure of a hermitian symmetric domain, cf. [13, Cor. 1.1.17]. For $h\in \mathcal X$ we consider the Hodge cocharacter
$$\mu_h:\mathbb G_{m,\mathbb C}\rightarrow G_{\mathbb C}$$ 
which maps $z\in\mathbb G_{m,\mathbb C}(\mathbb C)$ to $\mu_h(\mathbb C)(z)=h_{\mathbb C}(\mathbb C)(z,1)\in G_{\mathbb C}(\mathbb C)$. 

The most studied Shimura pairs are constructed as follows. Let $W$ be a vector space over $\mathbb Q$ of even dimension $2r$. Let $\psi$ be a non-degenerate, alternating bilinear form on $W$. Let $\mathcal Y$ be the set of all monomorphisms $\mathbb S\hookrightarrow \pmb{\text{GSp}}(W\otimes_{\mathbb Q} {\mathbb R},\psi)$ that define Hodge $\mathbb Q$--structures on $W$ of type $\{(-1,0),(0,-1)\}$ and that have either $2\pi i\psi$ or $-2\pi i\psi$ as polarizations. The pair $(\pmb{\text{GSp}}(W,\psi),\mathcal Y)$ is a Shimura pair that defines a {\it Siegel modular variety}. Let $L$ be a $\mathbb Z$-lattice of $W$ such that $\psi$ induces a perfect bilinear form $\psi:L\times L\rightarrow\mathbb Z$. Let 
$$K(N):=\{g\in \pmb{\text{GSp}}(L,\psi)(\widehat{\mathbb Z})|g\;\;\text{mod}\;\; N\widehat{\mathbb Z}\;\;\text{is}\;\;\text{identity}\}\;\;\text{and}\;\;K_p:=\pmb{\text{GSp}}(L,\psi)(\mathbb Z_p).$$
\indent
Let $E(G,\mathcal X)\hookrightarrow\mathbb C$ be the number subfield of $\mathbb C$ that is the field of definition of the $G(\mathbb C)$-conjugacy class of the cocharacters $\mu_h$'s of $G_{\mathbb C}$, cf. [37, p. 163]. We recall that $E(G,\mathcal X)$ is called the {\it reflex field} of $(G,\mathcal X)$. The {\it Shimura variety} $\text{Sh}(G,\mathcal X)$ is identified with the canonical model over $E(G,\mathcal X)$ of the complex Shimura variety 
$$\text{Sh}(G,\mathcal X)_{\mathbb C}:={\text{proj}.}{\text{lim}.}_{K\in \Sigma(G)} G(\mathbb Q)\backslash [\mathcal X\times (G(\mathbb A_f)/K)],$$ 
where $\Sigma(G)$ is the set of compact, open subgroups of $G(\mathbb A_f)$ endowed with the inclusion relation (see [12], [13], [36], [37], [38], and [39]). Thus $\text{Sh}(G,\mathcal X)$ is an $E(G,\mathcal X)$-scheme together with a continuous $G(\mathbb A_f)$-action. For $C$ a compact subgroup of $G(\mathbb A_f)$, let 
$$\text{Sh}_C(G,\mathcal X):=\text{Sh}(G,\mathcal X)/C.$$ 

Let $K\in \Sigma(G)$. We recall that a classical result of Baily and Borel allows us to view 
$$\text{Sh}_K(G,\mathcal X)({\mathbb C}):=G(\mathbb Q)\backslash [\mathcal X\times (G(\mathbb A_f)/K)]$$ 
as a finite, disjoint union of normal, quasi-projective varieties over $\mathbb C$ and not only of normal complex analytic spaces (see [1, Thm. 10.11]) and this makes $\text{Sh}_K(G,\mathcal X)$ to be a normal, quasi-projective $E(G,\mathcal X)$-scheme. If $K$ is small enough, then $\text{Sh}_K(G,\mathcal X)$ is in fact a smooth, quasi-projective $E(G,\mathcal X)$-scheme. We also recall that $\text{Sh}_K(G,\mathcal X)$ is a projective $E(G,\mathcal X)$-scheme if and only if the $\mathbb Q$--rank of $G^{\text{ad}}$ is $0$ (i.e., the Shimura pair $(G,\mathcal X)$ is compact), cf. [6, Thm. 12.3 and Cor. 12.4].

Let $H$ be a compact, open subgroup of $G_{\mathbb Q_p}(\mathbb Q_p)$. 

We recall that the group $G_{\mathbb Q_p}$ is called {\it unramified} if and only if it has a Borel subgroup and splits over an unramified, finite field extension of $\mathbb Q_p$.

See [52] for {\it hyperspecial subgroups} of $G_{\mathbb Q_p}(\mathbb Q_p)$. In what follows we will only use the following three properties of them:

\medskip
-- the group $G_{\mathbb Q_p}(\mathbb Q_p)$ has hyperspecial subgroups if and only if $G_{\mathbb Q_p}$ is unramified,

\smallskip
-- a subgroup of $G_{\mathbb Q_p}(\mathbb Q_p)$  is hyperspecial if and only if it is the group of $\mathbb Z_p$-valued points of a reductive group scheme over $\mathbb Z_p$ whose generic fibre is $G_{\mathbb Q_p}$, and

\smallskip
-- each hyperspecial subgroup of $G_{\mathbb Q_p}(\mathbb Q_p)$ is a maximal compact, open subgroup of $G_{\mathbb Q_p}(\mathbb Q_p)$. 

\medskip
Let $v$ be a prime of $E(G,\mathcal X)$ that divides $p$. Let $k(v)$ be the residue field of $v$. Let $e(v)\in\mathbb N^*$ be the absolute ramification index of $v$.  Let $O_{(v)}$ be the localization of the ring of integers of $E(G,\mathcal X)$ with respect to $v$. 

\begin{definitions}\label{D1}

{\bf (a)} By an {\it integral model} of $\text{Sh}_K(G,\mathcal X)$ over $O_{(v)}$ we mean a faithfully flat $O_{(v)}$-scheme whose generic fibre is $\text{Sh}_K(G,\mathcal X)$. 

\smallskip
{\bf (b)} By an {\it integral model} of $\text{Sh}_H(G,\mathcal X)$ over $O_{(v)}$ we mean a faithfully flat $O_{(v)}$-scheme equipped with a continuous
 $G(\mathbb A_f^{(p)})$-action whose generic fibre is the $E(G,\mathcal X)$-scheme $\text{Sh}_H(G,\mathcal X)$ equipped with its natural continuous $G(\mathbb A_f^{(p)})$-action. 
\end{definitions}

\medskip
In this paper we study integral models of $\text{Sh}_K(G,\mathcal X)$ and $\text{Sh}_H(G,\mathcal X)$ over $O_{(v)}$. The subject has a long history, the first main result being the existence of the moduli schemes $\mathcal A_{r,1,N}$ and $\mathcal M_r$. This is so as we have natural identifications 
$$\mathcal A_{r,1,N,\mathbb Q}=\text{Sh}_{K(N)}(\pmb{\text{GSp}}(W,\psi),\mathcal Y)\;\;\text{and}\;\;\mathcal M_{r,\mathbb Q}=\text{Sh}_{K_p}(\pmb{\text{GSp}}(W,\psi),\mathcal Y)$$ 
(see [12], [37], [54], etc.). In particular, see [54, Ex. 3.2.9 and Subsect. 4.1] and [12, Thm. 4.21] for the natural continuous action of $\pmb{\text{GSp}}(W,\psi)(\mathbb A_f^{(p)})$ on $\mathcal M_r$. 

In 1976 Langlands conjectured the existence of a good integral model of $\text{Sh}_H(G,\mathcal X)$ over $O_{(v)}$, provided $H$ is a hyperspecial subgroup of $G_{\mathbb Q_p}(\mathbb Q_p)$ (see [31, p. 411]); unfortunately, Langlands did not explain what good is supposed to stand for. Only in 1992, an idea of Milne made it significantly clearer how to characterize and identify the good integral models. Milne's philosophy can be roughly summarized as follows (cf. [37]): under certain conditions, the good regular, formally smooth integral models should be uniquely determined by (N\'eron type) universal properties that are similar to the property (iii) of Subsection 1.1. 

\begin{definitions}\label{D2}

{\bf (a)} We assume that $e(v)=1$. A flat, affine group scheme $G_{\mathbb Z_{(p)}}$ over $\mathbb Z_{(p)}$ that extends $G$ (i.e., whose generic fibre is $G$) is called a {\it quasi-reductive group scheme for $(G,\mathcal X,v)$} if there exists a reductive, normal, closed subgroup scheme $G^{\text{r}}_{\mathbb Z_{p}}$ of $G_{\mathbb Z_p}$ equipped with a cocharacter $\mu_v:\mathbb G_{m,W(k(v))}\rightarrow G^{\text{r}}_{\mathbb Z_p}\times_{\Spec(\mathbb Z_p)} \Spec(W(k(v)))$ such that the extension of $\mu_v$ to $\mathbb C$ via an (any) $O_{(v)}$-monomorphism $W(k(v))\hookrightarrow\mathbb C$ defines a cocharacter of $G_{\mathbb C}$ that is $G(\mathbb C)$-conjugate to the cocharacters $\mu_h$ of $G_{\mathbb C}$ introduced above ($h\in \mathcal X$).

\smallskip
{\bf (b)} We say that a smooth $O_{(v)}$-scheme $Y$ of finite type is a {\it N\'eron model} of its generic fibre $Y_{E(G,\mathcal X)}$ over $O_{(v)}$, if for each smooth $O_{(v)}$-scheme $Z$, every morphism $Z_{E(G,\mathcal X)}\rightarrow Y_{E(G,\mathcal X)}$ of $E(G,\mathcal X)$-schemes extends uniquely to a morphism $Z\rightarrow Y$ of $O_{(v)}$-schemes. 
\end{definitions}

\medskip
Definition \ref{D2} (a) is a variation of [51, Def. 1.5]; more precisely, the group $G_{\mathbb Z_{(p)}}(\mathbb Z_p)$ is an {\it h--hyperspecial} subgroup of $G_{\mathbb Q_p}(\mathbb Q_p)$ in the sense of loc. cit. Definition \ref{D2} (b) is well-known, cf. [7, Ch. 1, Sect. 1.2, Def. 1].  

The notion $G_{\mathbb Z_{(p)}}$ is a quasi-reductive group scheme for $(G,\mathcal X,v)$ is far more general then the notion $G_{\mathbb Z_{(p)}}$ is a reductive group scheme over $\mathbb Z_{(p)}$ that extends $G$. For instance, if $G_{\mathbb Z_{(p)}}$ is a reductive group scheme over $\mathbb Z_{(p)}$, then $G_{\mathbb Q_p}$ splits over an unramified finite field extension of $\mathbb Q_p$ but if $G_{\mathbb Z_{(p)}}$ is a quasi-reductive group scheme for $(G,\mathcal X,v)$, then nothing one can say in general about a finite field extension $\sharp$ of $\mathbb Q_p$ over which the group $G_{\mathbb Q_p}$ splits (and there exist plenty of examples in which $\sharp$ must contain an arbitrary a priori given finite field extension of $\mathbb Q_p$). To exemplify this, we will assume for the remaining part of this paragraph that $G$ is a $\mathbb Q$--simple adjoint group. It is a Weil restriction $\text{Res}_{\mathbb E/\mathbb Q} J$, where $\mathbb E$ is a totally real number field and where $J$ is an absolutely simple  adjoint group over $\mathbb E$, cf. [13, Subsect. 2.3.4 (a)]. If there exists a reductive group scheme $G_{\mathbb Z_{(p)}}$ over $\mathbb Z_{(p)}$ that extends $G$ (i.e., if $G_{\mathbb Q_p}$ is unramified), then $\mathbb E$ is unramified over $p$. We assume now that the set of primes of $\mathbb E$ above $p$ is the disjoint union of two non-empty sets $S_1$ and $S_2$ with the following two properties:

\medskip\noindent
{\bf (i)}  each $w\in S_1$ is unramified over $p$ and $J_{\mathbb E_w}$ is unramified (here $\mathbb E_w$ is the completion  of $\mathbb E$ at $w$); 

\smallskip\noindent
{\bf (ii)} each $w\in S_2$ corresponds naturally to compact simple factors of $G_{\mathbb R}=\prod_{\tau:\mathbb E\rightarrow \mathbb R} J_{\tau,\mathbb R}$ (via the identification of embeddings of $\mathbb E$ in $\mathbb R$ with embeddings of $\mathbb E$ in an algebraically closed field that contains both $\mathbb Q_p$ and $\mathbb R$).

\medskip\noindent
Then there exist finite, flat group schemes $G_{\mathbb Z_{(p)}}$ over $\mathbb Z_{(p)}$ which are quasi-reductive group scheme for $(G,\mathcal X,v)$. For instance, we can choose $G_{\mathbb Z_{(p)}}$ (see [54, Cl. 3.1.3.1]) such that its extension $G_{\mathbb Z_p}$ to $\mathbb Z_p$ is a product of the form $\mathcal G_1\times_{\Spec(\mathbb Z_p)} \mathcal G_2$ with 
$$\mathcal G_1=\prod_{w\in S_1} \text{Res}_{O_w/\mathbb Z_p} J_{O_w},$$ 
where $O_w$ is the ring of integers of $\mathbb E_w$ and where $J_{O_w}$ is a reductive group scheme over $O_w$ that extends $J_{\mathbb E_w}$. There exist no additional requirements from either $S_2$ or $\mathcal G_2$ and in particular $\mathbb E$ can be arbitrarily ramified at a prime $w\in S_2$.

\subsection{Constructing integral models} 
Until the end we will assume that the Shimura pair $(G,\mathcal X)$ is of {\it Hodge type}, i.e., there exists an injective map 
$$f:(G,\mathcal X)\hookrightarrow (\pmb{\text{GSp}}(W,\psi),\mathcal Y)$$ 
for some symplectic space $(W,\psi)$ over $\mathbb Q$; thus $f:G\hookrightarrow \pmb{\text{GSp}}(W,\psi)$ is a monomorphism such that we have $f_{\mathbb R}\circ h\in \mathcal Y$ for all elements $h\in \mathcal X$. 

We recall that we identity $\mathcal M_{r,\mathbb Q}=\text{Sh}_{K_p}(\pmb{\text{GSp}}(W,\psi),\mathcal Y)$. Let $L_{(p)}:=L\otimes_{\mathbb Z} \mathbb Z_{(p)}$. The schematic closure $G_{\mathbb Z_{(p)}}$ of $G$ in $\pmb{\text{GL}}_{L_{(p)}}$ is a flat, affine group scheme over $\mathbb Z_{(p)}$. Until the end we will also assume that we have an identity $H=K_p\cap G_{\mathbb Q_p}(\mathbb Q_p)$; thus $H=G_{\mathbb Z_{(p)}}(\mathbb Z_p)$. 

The functorial morphism $f_0:\text{Sh}(G,\mathcal X)\rightarrow\text{Sh}(\pmb{\text{GSp}}(W,\psi),\mathcal Y)_{E(G,\mathcal X)}$ defined by $f$ (see [12, Cor. 5.4]) is a closed embedding as it is so over $\mathbb C$ (cf. [12, Prop. 1.15]). The morphism $f_0$ induces naturally a morphism of $E(G,\mathcal X)$-schemes 
$$f_p:\text{Sh}_H(G,\mathcal X)\rightarrow \text{Sh}_{K_p}(\pmb{\text{GSp}}(W,\psi),\mathcal Y)_{E(G,\mathcal X)}$$ 
which is a closed embedding (cf. Fact \ref{F1}). Thus we can speak about the normalization $\mathcal N$ 
of the schematic closure of $\text{Sh}_H(G,\mathcal X)$ in $\mathcal M_{r,O_{(v)}}$. As $G(\mathbb A_f^{(p)})$ acts continuously on $\text{Sh}_H(G,\mathcal X)$ and $\mathcal M_r$, it is easy to see that we get naturally an induced continuous action of $G(\mathbb A_f^{(p)})$ on $\mathcal N$ (to be compared with [54, Prop. 3.4]). Let $\mathcal N^{\text{s}}$
be the formally smooth locus of $\mathcal N$ over $O_{(v)}$; it is a $G(\mathbb A_f^{(p)})$-invariant, open subscheme of $\mathcal N$ such that we have identities $\mathcal N^{\text{s}}_{E(G,\mathcal X)}=\mathcal N_{E(G,\mathcal X)}=\text{Sh}_H(G,\mathcal X)$ (cf. Lemma \ref{L1}). Let
$$(\mathcal A,\lambda_{\mathcal A})$$
be the principally polarized abelian scheme over $\mathcal N$ which is the natural pullback of the universal principally polarized abelian scheme over $\mathcal M_r$. 

If $p>2$ and $e(v)=1$, let $\mathcal N^{\text{m}}:=\mathcal N^{\text{s}}$. If $p=2$ and $e(v)=1$, let $\mathcal N^{\text{m}}$ be the $G(\mathbb A_f^{(p)})$-invariant, open subscheme of $\mathcal N^{\text{s}}$ defined in Subsection 3.5. In this paper we study when $e(v)=1$ the following sequence
$$\mathcal N^{\text{m}}\hookrightarrow\mathcal N^{\text{s}}\hookrightarrow\mathcal N\rightarrow\mathcal M_{r,O_{(v)}}$$
of morphisms of $O_{(v)}$-schemes in order to prove the following four basic results that pertain to $\mathcal N^{\text{s}}$. 

\begin{theorem} [Basic Theorem] \label{T1}
{\it We assume that $e(v)=1$ (i.e., $v$ is unramified over $p$) and that the $k(v)$-scheme $\mathcal N^{\text{s}}_{k(v)}$ is non-empty. Then the following three properties hold:

\smallskip
{\bf (a)} The $O_{(v)}$-scheme $\mathcal N^{\text{s}}$ is the unique regular, formally smooth integral model of $\text{Sh}_H(G,\mathcal X)$ over $O_{(v)}$ that satisfies the following smooth extension property: if $Z$ is a regular, formally smooth scheme over a discrete valuation ring $O$ which is of absolute ramification index $1$ and an $O_{(v)}$-algebra, then each morphism $Z_{E(G,\mathcal X)}\rightarrow\text{Sh}_H(G,\mathcal X)$ of $E(G,\mathcal X)$-schemes extends uniquely to a morphism $Z\rightarrow \mathcal N^{\text{s}}$ of $O_{(v)}$-schemes. 

\smallskip
{\bf (b)} For each algebraically closed field $k$ of characteristic $p$, the natural morphism $\mathcal N^{\text{s}}_{W(k)}\rightarrow\mathcal M_{r,W(k)}$ induces $W(k)$-epimorphisms at the level of complete, local rings of residue field $k$ (i.e., it is a formally closed embedding at all $k$-valued point of $\mathcal N^{\text{s}}_{W(k)}$).

\smallskip
{\bf (c)} We also assume that the $\mathbb Q$--rank of the adjoint group $G^{\text{ad}}$ is $0$. Let $H^{(p)}$ be a compact, open subgroup of $G(\mathbb A_f^{(p)})$ such that $\mathcal N$ is a pro-finite pro-\'etale cover of $\mathcal N/H^{(p)}$. Then $\mathcal N^{\text{s}}/H^{(p)}$ is a N\'eron model of its generic fibre $\text{Sh}_{H\times H^{(p)}} (G,\mathcal X)$ over $O_{(v)}$.} 
\end{theorem}

\begin{corollary}\label{C1}
{\it We assume that $e(v)=1$, that $G_{\mathbb Z_{(p)}}$ is smooth, that the $k(v)$-scheme $\mathcal N^{\text{m}}_{k(v)}$ is non-empty, and that $k$ is algebraically closed. Let $H^{(p)}$ be a compact, open subgroup of $G(\mathbb A_f^{(p)})$ such that $\mathcal N$ is a pro-finite pro-\'etale cover of $\mathcal N/H^{(p)}$ and there exists a natural number $N\ge 3$ relatively prime to $p$ such that $H\times H^{(p)}$ is a subgroup of $K(N)$ (thus we have a natural morphism  $\mathcal N^{\text{s}}_{k}/H^{(p)}\rightarrow\mathcal{A}_{r,1,N,k}=\mathcal M_{r,k}/K(N)^{(p)}$, where $K(N)^{(p)}$ is the compact, open subgroup of $\pmb{\text{GSp}}(L,\psi)(\mathbb A_f^{(p)})$ such that we have $K(N)=K_p\times K(N)^{(p)}$). Then each connected component of $\mathcal N^{\text{m}}_{k}/H^{(p)}$ is a quasi Shimura $p$-variety of Hodge type in the sense of [59, Def. 4.2.1] (and therefore for all $m\in\mathbb N^*$ it has a level $m$ stratification that enjoys all the nice properties listed in [59, Cor. 4.3]).}
\end{corollary}

\begin{remark}\label{R1}
{\bf (a)} We assume that $e(v)=1$. From many points of view (such as zeta functions) one would like to have a very good understanding of $\mathcal N$ itself. However, as of today, its formally smooth locus $\mathcal N^{\text{s}}$ is the only open subscheme of $\mathcal N$ which is uniquely determined by a universal property and for which the connected components of its geometric special fibers have under natural assumptions (see Corollary \ref{C1}) level $m$ stratifications for each $m\in\mathbb N^*$ that generalize the classical Ekedahl--Oort stratifications for $m=1$ and that have all the desired good properties (strata are regular, quasi-affine, equidimensional of dimensions given by concrete formulas, etc.). 

\smallskip
{\bf (b)} Under an additional condition satisfied for instance if $G_{\mathbb Z_{(p)}}$ is also a quasi-reductive group scheme for $(G,\mathcal X,v)$ and with the notation of Corollary \ref{C1}, $\mathcal N^{\text{m}}_{k}$ itself is a quasi Shimura $p$-variety of Hodge type in the sense of [59, Def. 4.2.1] (see Subsubsection 3.5.2 for details).
\end{remark}

\begin{proposition} \label{P1}
{\it We assume that $e(v)=1$ and that $G_{\mathbb Z_{(p)}}$ is a quasi-reductive group scheme for $(G,\mathcal X,v)$. Then all ordinary points of $\mathcal N_{k(v)}$ (i.e., all points $y:\Spec(k)\rightarrow\mathcal N_{k(v)}$ with values in perfect fields such that the abelian variety $y^*(\mathcal A)$ over $k$ is ordinary) belong to  $\mathcal N^{\text{m}}_{k(v)}$.} 
\end{proposition}

\begin{theorem} [Main Theorem] \label{T2} 
We assume that $e(v)=1$ and that $G_{\mathbb Z_{(p)}}$ is a quasi-reductive group scheme for $(G,\mathcal X,v)$. 

\medskip
{\bf (a)} Then $\mathcal N^{\text{m}}_{k(v)}$ is a non-empty, open closed subscheme of $\mathcal N_{k(v)}$. 

\smallskip
{\bf (b)} If the ordinary locus of $\mathcal N_{k(v)}$ is Zariski dense in $\mathcal N_{k(v)}$, then we have $\mathcal N^{\text{m}}=\mathcal N^{\text{s}}=\mathcal N$. 

\smallskip
{\bf (c)} If the $\mathbb Q$--rank of the adjoint group $G^{\text{ad}}$ is $0$, then the following two properties hold:

\medskip\noindent
{\bf (c.i)} We have $\mathcal N^{\text{m}}=\mathcal N^{\text{s}}=\mathcal N$ and moreover $\mathcal N$ is the integral canonical model of $\text{Sh}_H(G,\mathcal X)$ over $O_{(v)}$ as defined in [54, Def. 3.2.3 6)].

\smallskip\noindent
{\bf (c.ii)} Let $H^{(p)}$ be a compact, open subgroup of $G(\mathbb A_f^{(p)})$ such that $H\times H^{(p)}$ is contained in $K(N)$ for some $N\in\mathbb N\setminus (p\mathbb N\cup\{1,2\})$; thus we have a natural finite morphism 
$$f(N):\text{Sh}_{H\times H^{(p)}}(G,\mathcal X)\rightarrow \mathcal A_{r,1,N,E(G,\mathcal X)}=\text{Sh}_{K(N)}(\pmb{\text{GSp}}(W,\psi),\mathcal Y)_{E(G,\mathcal X)}.$$ 
Then the normalization $\mathcal Q$ of $\mathcal A_{r,1,N,O_{(v)}}$ in the ring of fractions of $\text{Sh}_{H\times H^{(p)}}(G,\mathcal X)$ is a smooth, projective $O_{(v)}$-scheme that can be identified with $\mathcal N/H^{(p)}$ and that is the N\'eron model of $\text{Sh}_{H\times H^{(p)}}(G,\mathcal X)$ over $O_{(v)}$.
\end{theorem} 

\subsection{On contents and proofs} 
We detail on the contents of this Part I. Section 2 lists conventions, notation, and few basic properties that pertain to the injective map $f:(G,\mathcal X)\hookrightarrow (\pmb{\text{GSp}}(W,\psi),\mathcal Y)$ and to Hodge cycles on abelian schemes over $\mathbb Q$--schemes. In connection to Sections 3 to 5 we assume that $e(v)=1$.

Section 3 includes crystalline applications. Until Subsection 3.3 we introduce basic notation and review three relatively recent results that pertain to Barsotti--Tate groups and that play a central role in Subsections 3.2 to 3.6, Sections 4 and 5, and Appendix B. The results are:

\medskip\noindent
{\bf (i)} {\it de Jong extension theorem} (see [11] and Theorem \ref{T3}),

\smallskip\noindent
{\bf (ii)} a variant of {\it Faltings deformation theory} (see Subsection 3.2), and

\smallskip\noindent
{\bf (iii)} a refinement of {\it a motivic conjecture of Milne} proved in [64, Thm. 1.2].

\medskip
Our first main new idea is to use (ii) in order to show directly that each $W(k)$-valued point of $\mathcal N$ factors through $\mathcal N^{\text{s}}$. Based on this and [67, Cor. 5], in Subsection 3.3 we prove the Basic Theorem \ref{T1}. Subsections 3.3 and 3.4 gather extra crystalline properties required in Sections 4 and 5 and required to prove Corollary \ref{C1} and a variant of it in Subsubsections 3.5.1 and 3.5.2; these two subsections can be viewed as an enlarged version with details of [59, Ex. 4.6]. Proposition \ref{P1} is proved in Subsection 3.6 based on [44] and on (iii). 

See Lemma \ref{L8} (a) for a simple criterion on when the $k(v)$-scheme $\mathcal N^{\text{m}}_{k(v)}$ is non-empty. In Subsection 4.1 we apply Theorem \ref{T1} (a) and Lemma \ref{L8} (a) to prove the existence and the uniqueness of good regular, formally smooth integral models of $\text{Sh}_{\tilde H}(G,\mathcal X)$ over $O_{(v)}$ for a large class of compact, open subgroups $\tilde H$ of $G_{\mathbb Q_p}(\mathbb Q_p)$ (the class includes all {\it parahoric} subgroups of $G_{\mathbb Q_p}(\mathbb Q_p)$) provided $G_{\mathbb Q_p}$ splits over an unramified extension of $\mathbb Q_p$ (see Theorem \ref{T7}). In particular, Corollary \ref{C3} can be viewed as a smooth solution (answer) to the conjecture (question) of Langlands (mentioned in the paragraph before Definitions \ref{D2}) for Shimura varieties of Hodge type. 

In Section 5 we use (i), Lemma \ref{L2} (i.e., [57, Cor. 4.3]), [66], and Subsection 3.3 to prove the Main Theorem \ref{T2} (see Subsections 5.1 to 5.6). Our second main new idea is to use (i) and purity results for reductive groups as in [66] in order to get that the open subscheme $\mathcal N_{k(v)}^{\text{m}}$ of $\mathcal N_{k(v)}$ is as well stable under specializations. 

Appendices A and B review basic properties of affine group schemes and of Barsotti--Tate groups. Their subsections are numbered as A1, A2, and B1 to B6. The reader ought to refer to these subsections only when they are quoted in the main text. Modulo few parts of the notation of Subsection 2.1, Appendices A and B are entirely independent of the main text.

\subsection{On literature} 
Referring to Theorem \ref{T1} (a), all ordinary points of $\mathcal N_{k(v)}$ belong to $\mathcal N^{\text{s}}_{k(v)}$ (cf. [44, Cor. 3.8]). Thus the only new part of Proposition \ref{P1} is the case $p=2$. If the $\mathbb Q$--rank of the adjoint group $G^{\text{ad}}$ is $0$ and $\mathcal N^{\text{s}}\neq \mathcal N$, then Theorem \ref{T1} (c) provides N\'eron models over $O_{(v)}$ which are not projective and thus which are not among the N\'eron models obtained in either [57, Prop. 4.4.1] or [67, Thm. 31]. Besides their applications to the conjecture of Langlands, Theorems \ref{T1} and \ref{T2} are also key steps in proving the deep conjectures [50, Conjs. B 3.7 and B 3.12] and [51, Conj. 1.6].

The uniqueness of an integral canonical model of $\text{Sh}_H(G,\mathcal X)$ over $O_{(v)}$ for $e(v)<p-1$ was proved in [54, Subsubsect. 3.2.17] (cf. also [54, Fact of Subsubsect. 3.2.12 or Rm. 3.2.4 stated for $e(v)=1<p-1$] and [55, Prop. 4.1], the last reference being a correction to the last part of [54, Step B of Subsubsect. 3.2.17]). The uniqueness of an integral canonical model of $\text{Sh}_H(G,\mathcal X)$ over $O_{(v)}$ for $e(v)\le p-1$ is also a particular case of [67, Cor. 30]. Moreover, a second proof of [67, Cor. 30] can be obtained based on [20, Thm. 1], which also corrects [40, Subsect. 3.6.1]. 

If $p\ge 5$ and $G_{\mathbb Z_{(p)}}$ is a reductive group scheme, then Theorem \ref{T2} (c.i) was first obtained in [54, Rm. 3.2.12, Thms. 5.1 and 6.4.1], [58, App.], and [67, Thm. 31]. If the Shimura pair $(G,\mathcal X)$ is unitary (i.e., $G^{\text{ad}}_{\mathbb C}$ is a non-trivial product of $\pmb{\text{PGL}}$ groups) and $G_{\mathbb Z_{(p)}}$ is a reductive group scheme, then Theorem \ref{T2} (c.i) follows also from [58, Thm. 5.1], [54, Subsubsect. 3.2.12], and [67, Thm. 31]. If $p\ge 3$ of if $p=2$ and the $2$-rank of each geometric fibre of the abelian scheme $\mathcal A_{k(v)}$ over $\mathcal N_{k(v)}$ is $0$ and if moreover $G_{\mathbb Z_{(p)}}$ is a reductive group scheme, then Theorem \ref{T2} (c.i) has been also claimed in [27] which relies on [40]. Similarly, if $p=2$ and $G_{\mathbb Z_{(p)}}$ is a reductive group scheme, then Theorem \ref{T2} (c.i) has been also claimed in [26].

Theorem \ref{T2} (c.i) represents progress towards the proof of a conjecture of Milne (see [37, Conj. 2.7] and [54, Conj. 3.2.5]) that pertains to the existence and the uniqueness of integral canonical models of arbitrary Shimura varieties.

The published works [43], [16], [41], [68], [32], [29], [54], [55], [57], [58], [62], [63], [64], [65], [27], and [26] are the most relevant ones for the existence of good smooth integral models of Shimura varieties of Hodge type. The construction of all integral models of Shimura varieties of Hodge type (such as $\mathcal N$ in Subsection 1.4) via normalizations of schematic closures in integral models of Siegel moduli varieties (i.e., in Mumford moduli schemes) used in all these references follows entirely [53] and [54], and thus are based on an original idea of Faltings shared with us in 1993. See also [23, Sect. 5] for a translation of part of [16] in terms of the existence of good smooth integral models in arbitrary ramified mixed characteristic $(0,p)$ of very simple unitary Shimura varieties.  

Part II will complete the proof of the conjecture of Milne on integral canonical models for the case of Shimura varieties of abelian type\break (see http://arxiv.org/abs/0712.1572). Part of Part II is also claimed in [27] and [26]. 

Part I brings completely new ideas in order to: 

\medskip\noindent
$\bullet$ shorten and simplify [54];

\smallskip\noindent
$\bullet$ extend many parts of [54] that were worked out only for $p\ge 5$ to the case of small primes $p\in\{2,3\}$;

\smallskip\noindent
$\bullet$ achieve progress towards the proofs of conjectures of Langlands, Milne, and Reimann; 

\smallskip\noindent
$\bullet$ work with large classes of subgroups of $G_{\mathbb Q_p}(\mathbb Q_p)$ which in the case when $G_{\mathbb Q_p}$ splits over an unramified extension of $\mathbb Q_p$ include as a very particular case the class of parahoric subgroups of $G_{\mathbb Q_p}(\mathbb Q_p)$ (and therefore also the class of hyperspecial subgroups of $G_{\mathbb Q_p}(\mathbb Q_p)$). 

\medskip
Theorem \ref{T2} (c.ii) for $p\ge 5$ corrects an error in the proof of [54, Prop. 3.2.3.2 ii)] that invalidated [54, Rm. 6.4.1.1 2) and most of Subsubsect. 6.4.11]. This correction was acknowledged and started in [57, Rm. 4.6 (b)] and [58, Thm. 5.1 (c) and App. E.8]. We recall that [58, App.] is the published erratum to [54].
 
The theory of local models aims to construct a projective scheme $\mathcal N_{\text{local}}$ over the completion $O_v$ of $O_{(v)}$ which among other things is expected to model the singularities of the complement $\mathcal N\setminus\mathcal N^{\text{s}}$ (in the pro-\'etale topology); for instance, see [48] and [46]. To our best knowledge, so far this theory has not been able to prove either the existence or the uniqueness of integral canonical models of Shimura varieties of Hodge type which are not of PEL type. But it has been able to say a lot about the nature of the singularities of the complement $\mathcal N\setminus\mathcal N^{\text{s}}$ in many cases. The most advanced work in this direction is the recent paper [28] which works in the case when $H$ is a parahoric subgroup of $G_{\mathbb Q_p}(\mathbb Q_p)$ and $G_{\mathbb Q_p}$ splits over a tamely ramified extension of $\mathbb Q_p$. 

If $G_{\mathbb Q_p}$ splits over an unramified extension of $\mathbb Q_p$ (thus $e(v)=1$), then the class of subgroups $H$ of $G_{\mathbb Q_p}(\mathbb Q_p)$ for which our results work is a lot more general than the class of parahoric subgroups of $G_{\mathbb Q_p}(\mathbb Q_p)$ considered in [28] (see Theorem \ref{T7}). But even when $e(v)=1$ neither this paper nor its Part II says anything about the complement $\mathcal N\setminus\mathcal N^{\text{s}}$ in the case when it is non-empty.

\section{Preliminaries}
In Subsection 2.1 we include some conventions and notation to be used throughout the paper. In Subsection 2.2 we study the injective map $f:(G,\mathcal X)\hookrightarrow (\pmb{\text{GSp}}(W,\psi),\mathcal Y)$. In Subsection 2.3 we consider $\mathbb C$-valued points of $\text{Sh}(G,\mathcal X)$ and different realizations of Hodge cycles on abelian schemes over reduced $\mathbb Q$--schemes. 

\subsection{Conventions and notation} 
We recall that $p$ is a prime and that $k$ is a perfect field of characteristic $p$. Let $\sigma:=\sigma_k$ be the Frobenius automorphism of $k$, $W(k)$, and of the field of fractions $B(k):=W(k)[\frac{1}{p}]$ of $W(k)$. For a Barsotti--Tate group $D$ over $W(k)$, let $H^1(D)$ be the dual of the Tate-module of $D_{B(k)}$. 

Let $R$, $M$, and $F$ be as in the beginning of Section 1. If $*$ or $*_R$ is either a morphism or an object of the category of $R$-schemes and if $S$ is a commutative $R$-algebra, let $*_S$ be the pullback of $*$ or $*_S$ to the category of $S$-schemes. Let $Z(F)$, $F^{\text{ad}}$, and $F^{\text{der}}$ denote the center, the adjoint group scheme, and the derived group scheme  (respectively) of $F$. We have $F^{\text{ad}}=F/Z(F)$. The group schemes $\pmb{\text{SL}}_{n,R}$, etc., are over $R$. If $F_1\hookrightarrow F$ is a closed embedding monomorphism of group schemes over $R$, then we identify $F_1$ with its image in $F$ and we consider intersections of subgroups of $F_1(R)$ with subgroups of $F(R)$. By the {\it essential tensor algebra} of $M\oplus M^{\vee}$ we mean the $R$-module
$$\mathcal T(M):=\oplus_{s,t\in\mathbb N} M^{\otimes s}\otimes_R M^{\vee\otimes t}.$$
\indent
Let $F^1(M)$ be a direct summand of $M$. Let $F^0(M):=M$ and $F^2(M):=0$. Let $F^1(M^{\vee}):=0$, $F^0(M^{\vee}):=\{y\in M^{\vee}|y(F^1(M))=0\}$, and $F^{-1}(M^{\vee}):=M^{\vee}$. Let $(F^i(\mathcal T(M)))_{i\in\mathbb Z}$ be the tensor product filtration of $\mathcal T(M)$ defined by the resulting exhaustive, separated filtrations $(F^i(M))_{i\in\mathbb Z}$ and $(F^i(M^{\vee}))_{i\in\mathbb Z}$ of $M$ and $M^{\vee}$ (respectively). We refer to $(F^i(\mathcal T(M)))_{i\in\mathbb Z}$ as the filtration of $\mathcal T(M)$ defined by $F^1(M)$. 

We identify naturally $\text{End}(M)=M\otimes_R M^{\vee}$ and $\text{End}(\text{End}(M))=M^{\otimes 2}\otimes_R M^{\vee\otimes 2}$. Let $x\in R$ be a non-divisor of $0$. A family of tensors of $\mathcal T(M[\frac{1}{x}])=\mathcal T(M)[\frac{1}{x}]$ is denoted $(u_{\alpha})_{\alpha\in\mathcal J}$, with $\mathcal J$ as the set of indexes. Let $M_1$ be another free $R$-module of finite rank. Let $(u_{1,\alpha})_{\alpha\in\mathcal J}$ be a family of tensors of $\mathcal T(M_1[\frac{1}{x}])$ indexed by the same set $\mathcal J$. By an isomorphism $(M,(u_{\alpha})_{\alpha\in\mathcal J})\rightarrow (M_1,(u_{1,\alpha})_{\alpha\in\mathcal J})$ we mean an $R$-linear isomorphism $M\rightarrow M_1$ that extends naturally to an $R[\frac{1}{x}]$-linear isomorphism $\mathcal T(M[\frac{1}{x}])\rightarrow\mathcal T(M_1[\frac{1}{x}])$ which takes $u_{\alpha}$ to $u_{1,\alpha}$ for all $\alpha\in\mathcal J$. We denote two tensors or bilinear forms in the same way, provided they are obtained one from another via either a reduction modulo some ideal or a scalar extension.

The notation $r$, $N$, $\mathcal A_{r,1,N}$, $\mathcal M_r$, $\mu_h:\mathbb G_{m,\mathbb C}\rightarrow G_{\mathbb C}$, $(\pmb{\text{GSp}}(W,\psi),\mathcal Y)$, $L$, $K(N)$, $K_p$, $E(G,\mathcal X)\hookrightarrow\mathbb C$, $\text{Sh}(G,\mathcal X)$, $\text{Sh}_C(G,\mathcal X)=\text{Sh}(G,\mathcal X)/C$, $v$, $k(v)$, $e(v)$, $O_{(v)}$, $f:(G,\mathcal X)\hookrightarrow (\pmb{\text{GSp}}(W,\psi),\mathcal Y)$, $L_{(p)}:=L\otimes_{\mathbb Z} \mathbb Z_{(p)}$, $G_{\mathbb Z_{(p)}}$, $H=K_p\cap G_{\mathbb Q_p}(\mathbb Q_p)=G_{\mathbb Z_{(p)}}(\mathbb Z_p)$, $f_0:\text{Sh}(G,\mathcal X)\rightarrow\text{Sh}(\pmb{\text{GSp}}(W,\psi),\mathcal Y)_{E(G,\mathcal X)}$, $f_p:\text{Sh}_H(G,\mathcal X)\rightarrow \text{Sh}_{K_p}(\pmb{\text{GSp}}(W,\psi),\mathcal Y)_{E(G,\mathcal X)}$, $\mathcal N$, $\mathcal N^{\text{s}}$, and $(\mathcal A,\lambda_{\mathcal A})$ will be as in Subsections 1.1, 1.3, and 1.4. Let $d:=\dim_{\mathbb C}(\mathcal X)\in\mathbb N$ and $l:=\dim(G)\in\mathbb N$. Let $G^0:=G\cap \pmb{\text{Sp}}(W,\psi)$. As $G_{\mathbb C}$ is the semidirect product of $G^0_{\mathbb C}$ and the image of any cocharacter $\mu_h:\mathbb G_{m,\mathbb C}\rightarrow G_{\mathbb C}$ with $h\in \mathcal X$, $G^0$ is a connected, reductive, normal subgroup of $G$. 

\subsection{On the injective map $f$} 
Let $H^{(p)}$ be an arbitrary compact, open subgroup of $G(\mathbb A_f^{(p)})$ such that $H\times H^{(p)}\leqslant K(N)$. As the morphism $f_0:\text{Sh}(G,\mathcal X)\rightarrow\text{Sh}(\pmb{\text{GSp}}(W,\psi),\mathcal Y)_{E(G,\mathcal X)}$ is a closed embedding, the induced morphisms 
$$f_p:\text{Sh}_H(G,\mathcal X)\rightarrow \text{Sh}_{K_p}(\pmb{\text{GSp}}(W,\psi),\mathcal Y)_{E(G,\mathcal X)}$$ 
and 
$$f_{H^{(p)}}:\text{Sh}_{H\times H^{(p)}}(G,\mathcal X)\rightarrow \text{Sh}_{K(N)}(\pmb{\text{GSp}}(W,\psi),\mathcal Y)_{E(G,\mathcal X)}$$ 
are pro-finite and finite (respectively). Thus we can speak about the normalization $\mathcal Q$ of $\mathcal A_{r,1,N,O_{(v)}}$ (equivalently, of the schematic closure in $\mathcal A_{r,1,N,O_{(v)}}$ of the image of $f_{H^{(p)}}$) in the ring of fractions of $\text{Sh}_{H\times H^{(p)}}(G,\mathcal X)$. We recall that every $O_{(v)}$-scheme of finite type is excellent (for instance, cf. [34, (34.A) and (34.B)]). The $O_{(v)}$-scheme $\mathcal A_{r,1,N,O_{(v)}}$ is quasi-projective (cf. property (i) of Subsection 1.1) and thus it is also excellent. Therefore the $O_{(v)}$-scheme $\mathcal Q$ is normal, quasi-projective, flat, has a relative dimension equal to $\dim(\text{Sh}_{H\times H^{(p)}}(G,\mathcal X))=\dim_{\mathbb C}(\mathcal X)=d$, and is finite over the $O_{(v)}$-scheme $\mathcal A_{r,1,N,O_{(v)}}$. 

Let $\mathcal Q^{\text{s}}$ be the smooth locus of $\mathcal Q$ over $O_{(v)}$; it is an open subscheme of $\mathcal Q$. As $\text{Sh}(\pmb{\text{GSp}}(W,\psi),\mathcal Y)$ is a pro-finite pro-\'etale cover of $\mathcal A_{r,1,N,\mathbb Q}=\text{Sh}_{K(N)}(\pmb{\text{GSp}}(W,\psi),\mathcal Y)$, the group $K(N)$ acts freely on $\text{Sh}(\pmb{\text{GSp}}(W,\psi),\mathcal Y)$. Thus the subgroup $H\times H^{(p)}$ of $K(N)$ acts freely on $\text{Sh}(\pmb{\text{GSp}}(W,\psi),\mathcal Y)$ and therefore also on $\text{Sh}(G,\mathcal X)$. Thus $\mathcal Q_{E(G,\mathcal X)}=\text{Sh}_{H\times H^{(p)}}(G,\mathcal X)$ is a smooth $E(G,\mathcal X)$-scheme and therefore it is the open subscheme $\mathcal Q^{\text{s}}_{E(G,\mathcal X)}$ of $\mathcal Q^{\text{s}}$.  

\begin{fact}\label{F1}
The finite morphism $f_p:\text{Sh}_H(G,\mathcal X)\rightarrow \text{Sh}_{K_p}(\pmb{\text{GSp}}(W,\psi),\mathcal Y)_{E(G,\mathcal X)}$ is in fact a closed embedding.
 \end{fact}

\medskip\proof
As $f_0$ is a closed embedding, it suffices to show that the map 
$$f_p(\mathbb C):\text{Sh}_H(G,\mathcal X)(\mathbb C)\rightarrow \text{Sh}_{K_p}(\pmb{\text{GSp}}(W,\psi),\mathcal Y)_{E(G,\mathcal X)}(\mathbb C)$$ 
is injective. But we have canonical identifications
$$\text{Sh}_H(G,\mathcal X)(\mathbb C)=G(\mathbb Q)\backslash [\mathcal X\times G(\mathbb A_f)]/H$$
and
$$\text{Sh}_{K_p}(\pmb{\text{GSp}}(W,\psi),\mathcal Y)_{E(G,\mathcal X)}(\mathbb C)=\pmb{\text{GSp}}(W,\psi)(\mathbb Q)\backslash [\mathcal Y\times  \pmb{\text{GSp}}(W,\psi)(\mathbb A_f)]/K_p$$
(cf. [13, Cor. 2.1.11]) and based on these and the fact that the intersections $G(\mathbb A_f^{(p)})\cap H$ and $\pmb{\text{GSp}}(W,\psi)(\mathbb A_f^{(p)})\cap K_p$ are the trivial subgroups of $G(\mathbb A_f)$ and $\pmb{\text{GSp}}(W,\psi)(\mathbb A_f)$ (respectively), one easily gets that $f_p(\mathbb C)$ is an injective map.\endproof

\begin{proposition}\label{P2}
The following three properties hold:

\medskip
{\bf (a)} The $O_{(v)}$-scheme $\mathcal N$ is a pro-finite pro-\'etale cover of $\mathcal Q$ and $\mathcal Q$ is the quotient of $\mathcal N$ by $H^{(p)}$.

\smallskip
{\bf (b)} The morphism $\mathcal N\rightarrow\mathcal M_{r,O_{(v)}}$ is finite.

\smallskip
{\bf (c)} We assume that $e(v)\le p-1$. If $Z$ is a regular, formally smooth scheme over a discrete valuation ring $O$ which is of absolute ramification index at most $p-1$ and an $O_{(v)}$-algebra, then each morphism $Z_{E(G,\mathcal X)}\rightarrow \mathcal N_{E(G,\mathcal X)}$ extends uniquely to a morphism $Z\rightarrow \mathcal N$ of $O_{(v)}$-schemes.
\end{proposition}

\medskip\proof
Let $N_1\in N\mathbb N\setminus p\mathbb N$. Let $N_2:=N$. For $i\in\{1,2\}$ we write $K(N_i)=K_p\times K(N_i)^{(p)}$, where the group $K(N_i)^{(p)}$ is a compact, open subgroup of $\pmb{\text{GSp}}(W,\psi)(\mathbb A_f^{(p)})$. The scheme $\mathcal M_r$ is a pro-finite pro-\'etale cover of $\mathcal M_r/K(N_i)^{(p)}=\mathcal A_{r,1,N_i}$. Let $H_i$ be a compact, open subgroup of $G(\mathbb A_f^{(p)})\cap K(N_i)^{(p)}$; thus $\text{Sh}(G,\mathcal X)$ is a pro-finite pro-\'etale cover of $\text{Sh}_{H\times H_i}(G,\mathcal X)$. The morphism $\text{Sh}_{H\times H_i}(G,\mathcal X)_{\mathbb C}\rightarrow \mathcal A_{r,1,N_i,\mathbb C}$ is of finite type and a formally closed embedding at each $\mathbb C$-valued point of $\text{Sh}_{H\times H_i}(G,\mathcal X)_{\mathbb C}$. Let $\mathcal Q_i$ be the normalization of $\mathcal A_{r,1,N_i,O_{(v)}}$ in the ring of fractions of $\text{Sh}_{H\times H_i}(G,\mathcal X)$; it is a finite $\mathcal A_{r,1,N_i,O_{(v)}}$-scheme and a normal, quasi-projective, flat $O_{(v)}$-scheme of relative dimension $d$. 

As $N_1\in N_2\mathbb N^*$, we have $K(N_1)^{(p)}\leqslant K(N_2)^{(p)}$. We assume that $H_1$ is a normal subgroup of $H_2$. The natural morphism $q_{12}:\mathcal Q_1\rightarrow\mathcal Q_2\times_{\mathcal A_{r,1,N_2,O_{(v)}}} \mathcal A_{r,1,N_1,O_{(v)}}$ of normal schemes is finite. We check that $q_{12,E(G,\mathcal X)}$ is an open closed embedding. As $q_{12,E(G,\mathcal X)}$ is a finite, \'etale morphism between normal $E(G,\mathcal X)$-schemes of finite type, it is enough to check that the map $q_{12}(\mathbb C):\mathcal Q_1(\mathbb C)\rightarrow\mathcal Q_2(\mathbb C)\times_{\mathcal A_{r,1,N_2,O_{(v)}}(\mathbb C)} \mathcal A_{r,1,N_1,O_{(v)}}(\mathbb C)$ is injective. We have 
$$\text{Sh}_{K_p\times H_i}(\pmb{\text{GSp}}(W,\psi),\mathcal Y)(\mathbb C)=\pmb{\text{GSp}}(L,\psi)(\mathbb Z_{(p)})\backslash [\mathcal Y\times (\pmb{\text{GSp}}(W,\psi)(\mathbb A_f^{(p)})/H_i)]$$
(for instance, cf. [38, Prop. 4.11]). Also we have a natural disjoint union decomposition
$$\text{Sh}_{H\times H_i}(G,\mathcal X)(\mathbb C)=\sqcup_{[g_j]\in G(\mathbb Q)\backslash G(\mathbb Q_p)/H} C_j\backslash [\mathcal X\times (G(\mathbb A_f^{(p)})/H_i)],\leqno (1)$$
where $g_j\in G(\mathbb Q_p)$ is a representative of the class $[g_j]\in G(\mathbb Q)\backslash G(\mathbb Q_p)/H$ and where the group $C_j:=G(\mathbb Q)\cap g_jHg_j^{-1}$ does not depend on $i\in\{1,2\}$. As we have an identity $\pmb{\text{GSp}}(W,\psi)(\mathbb Q_p)=\pmb{\text{GSp}}(W,\psi)(\mathbb Q)K_p$ (cf. [38, Lem. 4.9]), we can write $g_j=a_jh_j$, where $a_j\in \pmb{\text{GSp}}(W,\psi)(\mathbb Q)$ and $h_j\in K_p$. Thus
$$C_j\leqslant \pmb{\text{GSp}}(W,\psi)(\mathbb Q)\cap g_jK_pg_j^{-1}=\pmb{\text{GSp}}(W,\psi)(\mathbb Q)\cap a_jK_pa_j^{-1}=a_j\pmb{\text{GSp}}(L,\psi)(\mathbb Z_{(p)})a_j^{-1}=:C_j^{\text{bigg}}.$$ 
We have $C_j=G(\mathbb Q)\cap C_j^{\text{bigg}}$. This is so as $g_jHg_j^{-1}$ is the group of $\mathbb Z_p$-valued points of the schematic closure of $G$ in $a_j\pmb{\text{GSp}}(L,\psi)_{\mathbb Z_{(p)}}a_j^{-1}$. 

To show that the map $q_{12}(\mathbb C)$ is injective, it suffices to show that each one of the following commutative diagrams indexed by $j$
\[\xymatrix{
C_j\backslash [\mathcal X\times (G(\mathbb A_f^{(p)})/H_1)]  \ar[r]^<<<<<<{s_1} \ar[d]^{\pi_{12}} & \pmb{\text{GSp}}(L,\psi)(\mathbb Z_{(p)})\backslash [\mathcal Y\times (\pmb{\text{GSp}}(W,\psi)(\mathbb A_f^{(p)})/H_1)] \ar[d]^{\pi_{12}^{\text{bigg}}} \\
C_j\backslash [\mathcal X\times (G(\mathbb A_f^{(p)})/H_2))] \ar[r]^<<<<<{s_2} & \pmb{\text{GSp}}(L,\psi)(\mathbb Z_{(p)})\backslash [\mathcal Y\times (\pmb{\text{GSp}}(W,\psi)(\mathbb A_f^{(p)})/H_2)]
}\]
is such that the maps $\pi_{12}$ and $s_1$ define an injective map of $C_j\backslash [\mathcal X\times (G(\mathbb A_f^{(p)})/H_1)]$ into the fibre product of $s_2$ and $\pi_{12}^{\text{bigg}}$. Here the maps $\pi_{12}$ and $\pi_{12}^{\text{bigg}}$ are the natural projections. The maps $s_1$ and $s_2$ are defined by the rule: the equivalence class $[h,g]$, where $h\in \mathcal X$ and $g\in G(\mathbb A_f^{(p)})$, is mapped to the equivalence class $[a_j^{-1}h,a_j^{-1}g]$. Thus the fact that $\pi_{12}$ and $s_1$ define an injective map of $C_j\backslash [\mathcal X\times (G(\mathbb A_f^{(p)})/H_1)]$ into the fibre product of $s_2$ and $\pi_{12}^{\text{bigg}}$ is a direct consequence of the identity $C_j=G(\mathbb Q)\cap C_j^{\text{bigg}}$. Thus $q_{12}(\mathbb C)$ is injective. 

Therefore  $q_{12,E(G,\mathcal X)}$ is an open closed embedding. As $q_{12}$ is also a finite morphism of normal, flat $O_{(v)}$-schemes of finite type, $q_{12}$ itself is an open closed embedding. Thus $\mathcal Q_1$ is a finite \'etale cover of $\mathcal Q_2$ that in characteristic $0$ is a finite \'etale cover which (as $H_1\vartriangleleft H_2$) induces finite Galois covers between connected components. Therefore $\mathcal Q_1$ is a finite \'etale cover of $\mathcal Q_2$ which induces finite Galois covers between connected components. This implies that $Q_2$ is the quotient of $Q_1$ under the natural action of $H_2/H_1$ on it. 

By allowing $H_1$ to vary among the normal, open subgroups of $H_2$ and by a natural passage to limits, we get that $\mathcal N$ is a pro-finite pro-\'etale cover of $\mathcal Q_2$ and that $\mathcal Q_2=\mathcal N/H_2$. Thus by taking $H_2=H^{(p)}$, we get that $\mathcal Q=\mathcal Q_2$ and  that part (a) holds. 

As each morphism $q_{12}:\mathcal Q_1\rightarrow\mathcal Q_2\times_{\mathcal A_{r,1,N_2,O_{(v)}}} \mathcal A_{r,1,N_1,O_{(v)}}$ is an open closed embedding, by allowing $H_1$ to vary through all normal, open subgroups of $H_2$ we get that $\mathcal N$ is an open closed subscheme of $\mathcal Q_2\times_{\mathcal A_{r,1,N_2,O_{(v)}}} \mathcal M_{r,O_{(v)}}$ and thus part (b) holds. 

To prove part (c), we recall that $Z$ is a healthy regular scheme in the sense of either [54, Def. 3.2.1 2)] or [55] (cf. [67, Cor. 5]). Thus part (c) is implied by [54, Ex. 3.2.9 and Prop. 3.4.1], cf. definitions [54, Def. 3.2.3 2), 3), and 6)] (to be compared with the argument for the property (iii) of Subsection 1.1).\endproof 

\begin{remark}\label{R2}
Similar arguments to the ones that checked that $\mathcal N$ is a pro-finite pro-\'etale cover of $\mathcal N/H_2$ can be used to check that the right action of $G(\mathbb A_f^{(p)})$ on $\mathcal N$ is indeed a continuous action in the sense of [13, Subsubsect. 2.7.1] and in what follows we will use this property without any extra comment.
\end{remark}

\begin{lemma}\label{L1}
The scheme $\mathcal N^{\text{s}}$ is an open subscheme of $\mathcal N$ and $\mathcal N^{\text{s}}_{E(G,\mathcal X)}=\mathcal N_{E(G,\mathcal X)}$. Moreover, if $\mathcal N^{\text{s}}_{k(v)}$ is a non-empty scheme, then $\mathcal N^{\text{s}}$ together with the resulting action of $G(\mathbb A_f^{(p)})$ on it is a regular, formally smooth integral model of $\text{Sh}_H(G,\mathcal X)$ over $O_{(v)}$.
\end{lemma}

\medskip\proof
As $\mathcal N$ is a pro-finite pro-\'etale cover of the excellent, quasi-projective $O_{(v)}$-scheme $\mathcal Q$ (see Proposition \ref{P2} (a)), $\mathcal N^{\text{s}}=\mathcal N\times_{\mathcal Q} \mathcal Q^{\text{s}}$ is an open subscheme of $\mathcal N$. As $\mathcal Q_{E(G,\mathcal X)}=\mathcal Q^{\text{s}}_{E(G,\mathcal X)}$, we have $\mathcal N^{\text{s}}_{E(G,\mathcal X)}=\mathcal N_{E(G,\mathcal X)}$. The open subscheme $\mathcal N^{\text{s}}$ of $\mathcal N$ is $G(\mathbb A_f^{(p)})$-invariant. As $G(\mathbb A_f^{(p)})$ acts continuously on $\mathcal N$, it also acts continuously on $\mathcal N^{\text{s}}$. Thus, if the scheme $\mathcal N^{\text{s}}_{k(v)}$ is non-empty, then $\mathcal N^{\text{s}}$ together with the resulting continuous action of $G(\mathbb A_f^{(p)})$ on it is a regular, formally smooth integral model of $\text{Sh}_H(G,\mathcal X)$ over $O_{(v)}$.\endproof

\begin{fact}\label{F2}
We assume that there exists a simple factor $G_1$ of $G^{\text{ad}}_{\overline{\mathbb Q}}$ which is an $\pmb{\text{SO}}_{2n+1}$ group for some $n\in\mathbb N^*$. Let $G_2$ be the semisimple, normal subgroup of $G_{\overline{\mathbb Q}}$ whose adjoint is naturally identified with $G_1$. Then $G_2$ is a $\pmb{\text{Spin}}_{2n+1}$ group.
\end{fact}

\medskip\proof
The $\text{Lie}(G_2)$-module $W\otimes_{\mathbb Q} \overline{\mathbb Q}$ is non-trivial and its irreducible $\text{Lie}(G_2)$-submodules are associated to the weight $\varpi_n$ of the $B_n$ Lie type, cf. [38, p. 456]. Thus $G_2$ is a $\pmb{\text{Spin}}_{2n+1}$ group.\endproof

\begin{lemma}\label{L2}
If the $\mathbb Q$--rank of the adjoint group $G^{\text{ad}}$ is $0$, then $\mathcal Q$ is a projective $O_{(v)}$-scheme.
\end{lemma}

\medskip\proof
Let $G^\prime$ be the smallest subgroup of $G$ such that all elements $h\in \mathcal X$ factor through $G^\prime_{\mathbb R}$. It is a normal, reductive subgroup of $G$ that contains $G^{\text{der}}$; thus $G^{\prime\text{ad}}=G^{\text{ad}}$. Let $h^\prime\in \mathcal X$ be an element such that $G^\prime$ is the smallest subgroup of $\pmb{\text{GL}}_W$ with the property that $h^\prime$ factors through $G^\prime_{\mathbb R}$. We can assume that the $\mathbb C$-valued point $[h^\prime,1_W]\in\text{Sh}_{H\times H^{(p)}}(G,\mathcal X)$ is definable over a number field (here $1_W$ is the identity element of $G(\mathbb A_f)$ modulo $H^{(p)}$) and that $2\pi i\psi$ is a principal polarization of the Hodge $\mathbb Z$-structure on $L$ defined by $h^\prime$. Thus $G^\prime$ is the Mumford--Tate group of the principally polarized Hodge $\mathbb Z$--structure on $L$ defined by $h^\prime$ and $\psi$ and this principally polarized Hodge $\mathbb Z$--structure is associated naturally to a principally polarized abelian scheme over a number field. 

Let $\mathcal X^\prime$ be the $G^\prime(\mathbb R)$-conjugacy class of $h^\prime$. The pair $(G^\prime,\mathcal X^\prime)$ is a Shimura pair whose reflex field and dimension are  $E(G,\mathcal X)$ and $d$ (respectively). Let $H^\prime:=H\cap G^\prime(\mathbb Q_p)$ and $H^{\prime(p)}:=H^{(p)}\cap G^\prime(\mathbb A_f^{(p)})$. As the $\mathbb Q$--rank of $G^{\prime\text{ad}}=G^{\text{ad}}$ is $0$, as in [57, Prop. 2.7] we argue that the normalization $\mathcal Q^\prime$ of $\mathcal A_{r,1,N,O_{(v)}}$ in $\text{Sh}_{H^\prime\times H^{\prime(p)}}(G^\prime,\mathcal X^\prime)$ is a projective $O_{(v)}$-scheme provided the Morita conjecture holds for all abelian varieties over number fields. We recall from [47], [57], and [33], that the Morita conjecture predicts that each abelian variety over a number field with the property that a pullback of it over $\mathbb C$ has a Mumford--Tate group whose adjoint has $\mathbb Q$--rank $0$, has potentially good reduction everywhere. As the Morita conjecture holds (see [33]), we get that $\mathcal Q^\prime$ is a projective $O_{(v)}$-scheme. 

The Shimura variety $\text{Sh}(G^\prime,\mathcal X^\prime)$ is a closed subscheme of $\text{Sh}(G,\mathcal X)$ of the same dimension $d$ and therefore it is an open closed subscheme of $\text{Sh}(G,\mathcal X)$. Thus each connected component of the normalization of $\mathcal A_{r,1,N,O_{(v)}}$ (equivalently of $\mathcal Q$) in the ring of fractions of $\text{Sh}(G,\mathcal X)$ is a $G(\mathbb A_f)$-translation of a connected component of the normalization of $\mathcal A_{r,1,N,O_{(v)}}$ (equivalently of $\mathcal Q^\prime$) in the ring of fractions of $\text{Sh}(G^\prime,\mathcal X^\prime)$. Thus, as $\mathcal Q^\prime$ is a projective $O_{(v)}$-scheme, we get directly that $\mathcal Q$ is a projective $O_{(v)}$-scheme.\endproof 

\subsection{Tensors} 
The image of each $h\in\mathcal X$ contains $Z(\pmb{\text{GL}}_{W\otimes_{\mathbb Q} \mathbb R})$. Thus $Z(\pmb{\text{GL}}_W)\leqslant G$ and therefore each tensor of $\mathcal T(W^{\vee})$ fixed by $G$ belongs to the direct summand $\oplus_{u\in\mathbb N} W^{\vee\otimes u}\otimes_{\mathbb Q} W^{\otimes u}$ of $\mathcal T(W^{\vee})$. We consider a family of tensors $(v_{\alpha})_{\alpha\in\mathcal J}$
in $\bigsqcup_{u=0}^{\infty} W^{\vee\otimes u}\otimes_{\mathbb Q} W^{\otimes u}\subset\mathcal T(W^{\vee})$ such that $G$ is the subgroup of $\pmb{\text{GL}}_W$ that fixes $v_{\alpha}$ for all $\alpha\in\mathcal J$, cf. [14, Prop. 3.1 c)]. 

Let $\mathfrak{T}:\text{End}(W\otimes_{\mathbb Q} \mathbb Q_p)\times \text{End}(W\otimes_{\mathbb Q} \mathbb Q_p)\rightarrow\mathbb Q_p$ be the trace bilinear form on $\text{End}(W\otimes_{\mathbb Q} \mathbb Q_p)$. If $\flat$ is a reductive subgroup of $\pmb{\text{GL}}_{W\otimes_{\mathbb Q} \mathbb Q_p}$, then the restriction of $\mathfrak{T}$ to $\text{Lie}(\flat)$ is non-degenerate (cf. Lemma \ref{L11} (b)). Let $\pi_{\flat}$ be the projector of $\text{End}(W\otimes_{\mathbb Q} \mathbb Q_p)$ on $\text{Lie}(\flat)$ along the perpendicular on $\text{Lie}(\flat)$ with respect to $\mathfrak{T}$. If $G_{\mathbb Q_p}$ normalizes $\flat$, then $G_{\mathbb Q_p}$ fixes $\pi_{\flat}$. 

\subsubsection{Complex manifolds}
For a smooth $\mathbb C$-scheme $Y$, let $Y^{\text{an}}$ be the complex manifold associated naturally to $Y$. It is well-known that for each $u\in\mathbb N^*$ and for every abelian scheme $\Pi:C\rightarrow Y$ and its associated morphism $\Pi^{\text{an}}:C^{\text{an}}\rightarrow Y^{\text{an}}$ of complex manifolds, we have a natural isomorphism
$$R^u\Pi^{\text{an}}_*(\mathbb C)\rightarrow R^u \Pi^{\text{an}}_*(\Omega^{*}_{C^{\text{an}}/Y^{\text{an}}})^{\nabla_C^{\text{an}}}\leqno (2)$$ 
of complex sheaves on $Y^{\text{an}}$, where $\nabla^{\text{an}}_C$ is the connection on $R^u \Pi^{\text{an}}_*(\Omega^{*}_{C^{\text{an}}/Y^{\text{an}}})$ induced by the Gauss--Manin connection on $R^u \Pi_{*}(\Omega^{*}_{C/Y})$. 

\subsubsection{Hodge cycles} 
We will use the terminology of [14] on Hodge cycles on an abelian scheme $B_X$ over a reduced $\mathbb Q$--scheme $X$. Thus we write each Hodge cycle $v$ on $B_X$ as a pair $(v_{\text{dR}},v_{\acute et})$, where $v_{\text{dR}}$ and $v_{\acute et}$ are the de Rham and the \'etale component of $v$ (respectively). The \'etale component $v_{\acute et}$ at its turn has an $l$-component $v_{\acute et}^l$, for each rational prime $l$. 

In what follows we will be interested only in Hodge cycles on $B_X$ that involve no Tate twists and that are tensors of different essential tensor algebras. Accordingly, if $X$ is the spectrum of a field $E$, then in applications $v_{\acute et}^p$ will be a suitable $\text{Gal}(\overline{E}/E)$-invariant tensor of $\mathcal T(H^1_{\acute et}(B_{\overline{X}},\mathbb Q_p))$, where $\overline{X}:=\Spec(\overline{E})$. If moreover $\overline{E}$ is a subfield of $\mathbb C$, then we will also use the Betti realization $v_B$ of $v$: it is a tensor of $\mathcal T(H^1((B_X\times_X \Spec(\mathbb C))^{\text{an}},\mathbb Q))$ that corresponds to $v_{\text{dR}}$ (resp. to $v_{\acute et}^l$) via the canonical isomorphism that relates the Betti cohomology of $(B_X\times_X \Spec(\mathbb C))^{\text{an}}$ with $\mathbb Q$--coefficients with the de Rham (resp. the $\mathbb Q_l$ \'etale) cohomology of $B_{\overline{X}}$ (see [14, Sect. 2]). We recall that $v_B$ is also a tensor of the $F^0$-filtration of the Hodge filtration of $\mathcal T(H^1((B_X\times_X \Spec(\mathbb C))^{\text{an}},\mathbb C))$. 

\subsubsection{On $\mathcal A_{E(G,\mathcal X)}$} 
The choice of the $\mathbb Z$-lattice $L$ of $W$ and of the family of tensors $(v_\alpha)_{\alpha\in\mathcal J}$ allows a moduli interpretation of $\text{Sh}(G,\mathcal X)$ (see  [12], [13], [38], and [54, Subsect. 4.1, Lem. 4.1.3]). For instance, $\text{Sh}(G,\mathcal X)(\mathbb C)=G(\mathbb Q)\backslash (\mathcal X\times G(\mathbb A_f))$ is the set of isomorphism classes of principally polarized abelian
varieties over $\mathbb C$ of dimension $r$, that carry a family of Hodge
cycles indexed by $\mathcal J$, that have compatible level-$N$ symplectic similitude structures for each $N\in\mathbb N^*$, and that satisfy few axioms. Thus the abelian scheme $\mathcal A_{E(G,\mathcal X)}$ over $\mathcal N_{E(G,\mathcal X)}$ is endowed with a family $(w_{\alpha}^{\mathcal A})_{\alpha\in\mathcal J}$ 
of Hodge cycles; all realizations of pullbacks of $w_{\alpha}^{\mathcal A}$ via $\mathbb C$-valued points of $\mathcal N^{\text{s}}_{E(G,\mathcal X)}$ correspond naturally to $v_{\alpha}$.

\begin{lemma}\label{L3} 
Let $w\in\text{Sh}(G,\mathcal X)(\mathbb C)$. We denote also by $w$ the $\mathbb C$-valued point of $\mathcal N$ defined by $w$; thus we can define $(A_w,\lambda_{A_w}):=w^*((\mathcal A,\lambda_{\mathcal A}))$. Let $u^w_{\alpha}$ (resp. $t^w_{\alpha}$) be the $p$-component of the \'etale component (resp. be the de Rham component) of the Hodge cycle $w^*(w_{\alpha}^{\mathcal A})$ on $A_w$. We have:

\medskip
{\bf (a)} There exist isomorphisms $(H^1_{\acute et}(A_w,\mathbb Z_p),(u_{\alpha}^w)_{\alpha\in\mathcal J})\rightarrow (L^{\vee}_{(p)}\otimes_{\mathbb Z_{(p)}} \mathbb Z_p,(v_{\alpha})_{\alpha\in\mathcal J})$ that take the perfect bilinear form on $H^1_{\acute et}(A_w,\mathbb Z_p)$ defined by $\lambda_{A_w}$ to a $\mathbb G_{m,\mathbb Z_p}(\mathbb Z_p)$-multiple of the perfect bilinear form $\psi^{\vee}$ on $L^{\vee}_{(p)}\otimes_{\mathbb Z_{(p)}} \mathbb Z_p$ defined by $\psi$.

\smallskip
{\bf (b)} There exist isomorphisms $(H^1_{\text{dR}}(A_w,\mathbb C),(t_{\alpha}^w)_{\alpha\in\mathcal J},\psi_{H^1_{\text{dR}}(A_w,\mathbb C)})\rightarrow (W^{\vee}\otimes_{\mathbb Q} \mathbb C,(v_{\alpha})_{\alpha\in\mathcal J},\psi^{\vee})$, where $\psi_{H^1_{\text{dR}}(A_w,\mathbb C)}$ is the perfect bilinear form on $H^1_{\text{dR}}(A_w,\mathbb C)$ defined by $\lambda_{A_w}$.
\end{lemma}

\medskip\proof
We write $w=[h_w,g_w]\in\text{Sh}(G,\mathcal X)(\mathbb C)=G(\mathbb Q)\backslash (\mathcal X\times G(\mathbb A_f))$, where $h_w\in\mathcal X$ and $g_w\in G(\mathbb A_f)$. From the standard moduli interpretation of $\text{Sh}(G,\mathcal X)(\mathbb C)$ applied to $w\in\text{Sh}(G,\mathcal X)(\mathbb C)$ we get (see [12], [37], [38], and [54, p. 454]) that the complex manifold $A_w^{\text{an}}$ associated to $A_w$ is $L_w\backslash W\otimes_{\mathbb Q} \mathbb C/F^{0,-1}_w$, where:

\medskip\noindent
{\bf (i)} $L_w$ is the $\mathbb Z$-lattice of $W$ defined uniquely by the identity $L_w\otimes_{\mathbb Z} \widehat{\mathbb Z}=g_w(L\otimes_{\mathbb Z} \widehat{\mathbb Z})$; 

\smallskip\noindent
{\bf (ii)} $W\otimes_{\mathbb Q} \mathbb C=F^{0,-1}_w\oplus F^{-1,0}_w$ is the usual Hodge decomposition of the Hodge $\mathbb Q$--structure on $W$ defined by $h_w\in \mathcal X$;

\smallskip\noindent
{\bf (iii)}  the principal polarization $\lambda_{A_w}$ of $A_w$ is defined naturally by a uniquely determined (non-zero) rational multiple of $\psi$;

\smallskip\noindent
{\bf (iv)}  under the canonical identifications $H^1_{\text{dR}}(A_w/\mathbb C)=H^1_{\text{dR}}(A_w^{\text{an}}/\mathbb C)=W^{\vee}\otimes_{\mathbb Q} \mathbb C=L_w^{\vee}\otimes_{\mathbb Z} \mathbb C$, each tensor $t_{\alpha}^w$ gets identified with $v_{\alpha}$ for all $\alpha\in\mathcal J$ and $\psi_{H^1_{\text{dR}}(A_w,\mathbb C)}$ gets identified with a (non-zero) complex multiple of $\psi^{\vee}$.

\medskip
Thus $(H^1_{\acute et}(A_w,\mathbb Z_p),(u_{\alpha}^w)_{\alpha\in\mathcal J})$ is identified naturally with $(L_w^{\vee}\otimes_{\mathbb Z} \mathbb Z_p,(v_{\alpha})_{\alpha\in\mathcal J})$  (cf. (iv)) and therefore also with a $G_{\mathbb Q_p}(\mathbb Q_p)$-conjugate of $(L^{\vee}_{(p)}\otimes_{\mathbb Z_{(p)}} \mathbb Z_p,(v_{\alpha})_{\alpha\in\mathcal J})$ (cf. (i)). Part (a) follows from this and from the existence of the rational multiple of $\psi$ mentioned in the property (iii).

For each non-zero complex number $\epsilon$, the automorphism of $(W^{\vee}\otimes_{\mathbb Q} \mathbb C,(v_{\alpha})_{\alpha\in\mathcal J})$ defined by $\mu_{h_w}(\mathbb C)(\epsilon^{-1})$ acts on the $\mathbb C$-span of $\psi^{\vee}$ as the multiplication by $\epsilon$. From this and the property (iv) we get that part (b) holds.\endproof 

\begin{lemma}\label{L4} 
Let $m\in\mathbb N$. Let $\mathcal R_1:=\mathbb C[[x_1,\ldots,x_m]]$, where $x_1,\ldots,x_m$ are independent variables. Let $\mathcal I_1:=(x_1,\ldots,x_m)$ be the maximal ideal of $\mathcal R_1$. Let $s\in\mathbb N^*$. Let $A_{w,s}$ be an abelian scheme over $\mathcal R_1/\mathcal I_1^s$ that is a deformation of $A_w$ (i.e., we have $A_w= A_{w,s}\times_{\Spec(\mathcal R_1/\mathcal I_1^s)} \Spec(\mathcal R_1/\mathcal I_1)$) and which has a principal polarization. Then there exists a unique isomorphism
$$I_{w,s}:H^1_{\text{dR}}(A_{w,s}/(\mathcal R_1/\mathcal I_1^s))\rightarrow H^1_{\text{dR}}(A_w/\mathbb C)\otimes_{\mathbb C} \mathcal R_1/\mathcal I_1^s$$
that has the following two properties:

\medskip\noindent
{\bf (i)} it lifts (i.e., modulo $\mathcal I_1/\mathcal I_1^s$ is) the identity automorphism of $H^1_{\text{dR}}(A_w/\mathbb C)$;

\smallskip\noindent
{\bf (ii)} under it, the Gauss--Manin connection on $H^1_{\text{dR}}(A_{w,s}/(\mathcal R_1/\mathcal I_1^s))$ becomes isomorphic to the flat connection $\delta$ on the $\mathcal R_1/\mathcal I_1^s$-module $H^1_{\text{dR}}(A_w/\mathbb C)\otimes_{\mathbb C} \mathcal R_1/\mathcal I_1^s$ that annihilates $H^1_{\text{dR}}(A_w/\mathbb C)\otimes 1$.
\end{lemma}

\medskip\proof
The uniqueness of $I_{w,s}$ is implied by the fact that $H^1_{\text{dR}}(A_w/\mathbb C)\otimes 1$ is the set of all elements of $H^1_{\text{dR}}(A_w/\mathbb C)\otimes_{\mathbb C} \mathcal R_1/\mathcal I_1^s$ that are annihilated by $\delta$. We consider an abelian scheme $\varPi:A_Y\rightarrow Y$ over a smooth $\mathbb C$-scheme $Y$ of dimension $m$ which is a global deformation of $A_{w,s}\rightarrow\Spec(\mathcal R_1/\mathcal I_1^s)$ and which has a principal polarization. Let $Z^{\text{an}}$ be a simply connected open complex submanifold of $Y^{\text{an}}$ that contains the $\mathbb C$-valued of $Y$ point defined naturally by $A_w$. We identify naturally $\Spec(\mathcal R_1/\mathcal I_1^s)$ with a complex analytic subspace of $Y^{\text{an}}$ and thus also of $Z^{\text{an}}$. We apply Formula (2) with $u=1$ and with $\Pi:C\rightarrow Y$ replaced by $\varPi:A_Y\rightarrow Y$. The pullback of $R^1\varPi_*^{\text{an}}(\mathbb C)$ to $Z^{\text{an}}$ is a constant sheaf on $Z^{\text{an}}$. Thus by pulling back Formula (2) to the complex analytic subspace $\Spec(\mathcal R_1/\mathcal I_1^s)$ of $Z^{\text{an}}$, we get directly the existence of $I_{w,s}$ for which  properties (i) and (ii) hold.\endproof

\begin{corollary}\label{C2}
{\it Let $m$, $\mathcal R_1$, and $\mathcal I_1$ be as in Lemma \ref{L4}. Let $A_{w,\infty}$ be an abelian scheme over $\mathcal R_1$ that is a deformation of $A_w$ and that has a principal polarization. Then there exists a unique isomorphism
$$I_{w,\infty}:H^1_{\text{dR}}(A_{w,\infty}/\mathcal R_1)\rightarrow H^1_{\text{dR}}(A_w/\mathbb C)\otimes_{\mathbb C} \mathcal R_1$$
that has the following two properties:

\medskip\noindent
{\bf (i)} it lifts (i.e., modulo $\mathcal I_1$ is) the identity automorphism of $H^1_{\text{dR}}(A_w/\mathbb C)$;

\smallskip\noindent
{\bf (ii)} under it, the $\mathcal I_1$-completion of the Gauss--Manin connection on $H^1_{\text{dR}}(A_{w,\infty}/\mathcal R_1)$ becomes isomorphic to the $\mathcal I_1$-completion of the flat connection $\delta$ on the $\mathcal R_1$-module $H^1_{\text{dR}}(A_w/\mathbb C)\otimes_{\mathbb C} \mathcal R_1$ that annihilates $H^1_{\text{dR}}(A_w/\mathbb C)\otimes 1$.

\medskip
If $w_{\alpha}^{\mathcal R_1}$ (resp. $\lambda_{A_{w,\infty}}$) is a Hodge cycle on (resp. a principal polarization of) $A_{w,\infty}$ that lifts the Hodge cycle $w^*(w_{\alpha}^{\mathcal A})$ on $A_w$ (resp. lifts the principal polarization $\lambda_{A_w}$ of $A_w$), then the isomorphism 
$$I_{w,\infty}:\mathcal T(H^1_{\text{dR}}(A_{w,\infty}/\mathcal R_1))\rightarrow \mathcal T(H^1_{\text{dR}}(A_w/\mathbb C))\otimes_{\mathbb C} \mathcal R_1$$ 
induced naturally by $I_{w,\infty}$ (and denoted in the same way) takes the de Rham realization of $w_{\alpha}^{\mathcal R_1}$ (resp. of $\lambda_{A_{w,\infty}}$) to $t_{\alpha}^w$ (resp. to the de Rham realization of $\lambda_{A_w}$).}
\end{corollary}

\medskip\proof
The existence and the uniqueness of $I_{w,\infty}$ follow from Lemma \ref{L4}  by taking $s\rightarrow\infty$. It is well-known that each de Rham component of a Hodge cycle on $A_{w,\infty}$ is annihilated by the Gauss--Manin connection on $\mathcal T(H^1_{\text{dR}}(A_{w,\infty}/\mathcal R_1))$. For instance, this follows from [14, Prop. 2.5] via a natural algebraization process. Thus $I_{w,\infty}(w_{\alpha}^{\mathcal R_1})$ and $t_{\alpha}^w$ are tensors of $\mathcal T(H^1_{\text{dR}}(A_w/\mathbb C))\otimes_{\mathbb C} \mathcal R_1$ which are annihilated by the $\mathcal I_1$-completion of the flat connection on $\mathcal T(H^1_{\text{dR}}(A_w/\mathbb C))\otimes_{\mathbb C} \mathcal R_1$ induced by $\delta$ and which modulo $\mathcal I_1$ coincide. Therefore the two tensors coincide, i.e., we have $I_{w,\infty}(w_{\alpha}^{\mathcal R_1})=t_{\alpha}^w$. A similar argument shows that $I_{w,\infty}$ takes $\lambda_{A_{w,\infty}}$ to the de Rham realization of $\lambda_{A_w}$.\endproof

\section{Crystalline applications}
Theorem \ref{T3} recalls a variant of the main result of [11]. In Subsection 3.1 we first introduce notation required to prove Theorems \ref{T1} and \ref{T2} and then we apply the main result of [64] in the form recalled in Theorem \ref{T9}. In Subsection 3.2 we apply the deformation theory of [17, Sect. 7]. Subsection 3.3 proves the Basic Theorem \ref{T1}. Subsection 3.4 introduces de Rham realizations of certain Hodge cycles. Subsection 3.5 defines the open subscheme $\mathcal N^{\text{m}}$ of $\mathcal N^{\text{s}}$, proves Corollary \ref{C1} and a variant of it, and lists few simple crystalline properties that are required in Sections 4 and 5. Subsection 3.6 proves Proposition \ref{P1} based also on Lemma \ref{L7}. Throughout this section we assume that $e(v)=1$.

For (crystalline or de Rham) Fontaine comparison theory we refer to [19], [17, Sect. 5], and [64]; see also Subsections B2, B3 and B6. We recall that $k$ is a perfect field of characteristic $p$. As the Verschiebung maps of Barsotti--Tate groups are not mentioned at all in what follows, we use the terminology $F$-crystals (resp. filtered $F$-crystals) associated to Barsotti--Tate groups over $k$, $k[[x]]$, or $k((x))$ (resp. over $W(k)$ or $W(k)[[x]]$) instead of the terminology Dieudonn\'e $F$-crystals (resp. filtered Dieudonn\'e $F$-crystals) of [4, Ch. 3] and [3, Chs. 2 and 3]. 

Let $x$ be an independent variable. The simplest form of [11, Thm. 1.1] says:

\begin{theorem} [de Jong] \label{T3}
{\it The natural functor from the category of $F$-crystals over $k[[x]]$ to the category of $F$-crystals over $k((x))$ is fully faithful.}
\end{theorem}

\subsection{Basic setting} 
From now on until the end, the field $k$ will be assumed to be algebraically closed and we will use the notation of Subsection 2.1. Let $z\in\mathcal N(W(k))$. Let 
$$(A,(w_{\alpha})_{\alpha\in\mathcal J},\lambda_A):=z^*(\mathcal A,(w_{\alpha}^{\mathcal A})_{\alpha\in\mathcal J},\lambda_{\mathcal A}).$$ 
Let 
$$(M,F^1,\phi,\psi_M)$$
be the principally quasi-polarized filtered $F$-crystal over $k$ of the principally quasi-polarized Barsotti--Tate group $(D,\lambda_D)$ of $(A,\lambda_A)$. Thus $\psi_M:M\times M\rightarrow W(k)
 $ is a perfect, alternating bilinear form on the free $W(k)$-module $M$ of rank $2r$, $F^1$ is a maximal isotropic submodule of $M$ with respect to $\psi_M$, and $\phi:M\rightarrow M$ is a $\sigma$-linear endomorphism such that we have an inclusion $pM\subset \phi(M)$ as well as identities $\psi_M(\phi(a),\phi(b))=p\sigma(\psi_M(a,b))$ for all $a,b\in M$. The $\sigma$-linear automorphism $\phi$ of $M[\frac{1}{p}]$ acts on $M^{\vee}[\frac{1}{p}]$ by mapping $e\in M^{\vee}[\frac{1}{p}]$ to $\sigma\circ e\circ\phi^{-1}\in M^{\vee}[\frac{1}{p}]$ and it acts on $\mathcal T(M)[\frac{1}{p}]$ in the natural tensor product way. Let $\psi_{H^1(D)}$ and $\psi_{H^1_{\acute{et}}}$ be the perfect, alternating bilinear forms on $H^1(D)=H^1_{\acute{et}}(A_{\overline{B(k)}},\mathbb Z_p)$ and $H^1_{\acute{et}}(A_{\overline{B(k)}},\mathbb Z_p)$ (respectively) defined by $\lambda_D$. We have a canonical identification of $\text{Gal}(B(k))$-modules (cf. property (ii) of Subsection B5):
$$(H^1(D),\psi_{H^1(D)})=(H^1_{\acute{et}}(A_{\overline{B(k)}},\mathbb Z_p),\psi_{H^1_{\acute{et}}}).\leqno (3)$$

Let $t_{\alpha}$ and $u_{\alpha}$ be the de Rham component of $w_{\alpha}$ and the $p$-component of the \'etale component of $w_{\alpha}$ (respectively). We have $u_{\alpha}\in \mathcal T(H^1(D))[\frac{1}{p}]$ for all $\alpha\in\mathcal J$. If $(F^i(\mathcal T(M)))_{i\in\mathbb Z}$ is the filtration of $\mathcal T(M)$ defined by $F^1$, then we have $t_{\alpha}\in F^0(\mathcal T(M))[\frac{1}{p}]$ for all $\alpha\in\mathcal J$. Let $\mathcal G$
be the schematic closure in $\pmb{\text{GL}}_M$ of the subgroup of $\pmb{\text{GL}}_{M[\frac{1}{p}]}$ that fixes $t_{\alpha}$ for all $\alpha\in\mathcal J$; it is a flat, affine group scheme over $W(k)$. It is known that $w_{\alpha}$ is a de Rham cycle, i.e., $t_{\alpha}$ and $u_{\alpha}$ correspond to each other via de Rham and thus also crystalline Fontaine comparison theory (see [61, Thm. 5.1.6 and Cor 5.1.7]). Thus $\phi(t_{\alpha})=t_{\alpha}$ for all $\alpha\in\mathcal J$.
 
Let $\mu:\mathbb G_{m,W(k)}\rightarrow \pmb{\text{GL}}_M$ be the inverse of the canonical split cocharacter of $(M,F^1,\phi)$ defined in [69, p. 512]. The cocharacter $\mu$ acts on $F^1$ via the weight $-1$ and fixes a direct supplement $F^0$ of $F^1$ in $M$; therefore we have $M=F^1\oplus F^0$. Moreover, $\mu$ fixes each tensor $t_{\alpha}$ (cf. the functorial aspects of [69, p. 513]). Thus $\mu$ factors through $\mathcal G$. Let
$$\mu:\mathbb G_{m,W(k)}\rightarrow\mathcal G$$ 
be the resulting factorization. We emphasize that in connection to different Kodaira--Spencer maps, in what follows we will identity naturally $\text{Hom}(F^1,F^0)$ with the direct summand $\{e\in\text{End}(M)|e(F^0)=0,\;\;e(F^1)\subset F^0\}\simeq \text{Hom}(M/F^0,F^0)$ of $\text{End}(M)$. 

Let $\mathcal G^0=\mathcal G\cap  \pmb{\text{Sp}}(M,\psi_M)$. As $\mathcal G$ is the semidirect product of $\mathcal G^0$ and of the image of $\mu:\mathbb G_{m,W(k)}\rightarrow\mathcal G$, we get that $\mathcal G^0$ is a flat, closed subgroup scheme of $\pmb{\text{GSp}}(M,\psi_M)$ which is smooth or reductive if and only if $\mathcal G$ is so. We consider the following family of principally quasi-polarized Dieudonn\'e modules with a group over $k$ associated to $z$:
$$\mathfrak{F}=\{(M,g\phi,\psi_M,\mathcal G^0)|g\in\mathcal G^0(W(k))\}.$$

\begin{lemma}\label{L5} 
The direct summand $\text{Lie}(\mathcal G_{B(k)})\cap\text{Hom}(F^1,F^0)$ of $\text{End}(M)$ has rank $d$. Moreover $\mathcal G_{B(k)}$ is a form of $G_{B(k)}$ and thus a reductive group.
\end{lemma}

\medskip\proof
To prove the lemma we can assume that $k$ has countable transcendental degree; thus there exists an $O_{(v)}$-monomorphism $W(k)\hookrightarrow\mathbb C$. Let $\mathcal F_{B(k)}$ be the normalizer of $F^1[\frac{1}{p}]$ in $\mathcal G_{B(k)}$. Its Lie algebra is equal to $\text{Lie}(\mathcal G_{B(k)})\cap \{e\in\text{End}(M)[\frac{1}{p}]|e(F^1[\frac{1}{p}])\subset F^1[\frac{1}{p}]\}$. As $\mu$ factors through $\mathcal G$, we have a direct sum decomposition $\text{Lie}(\mathcal G_{B(k)})=\text{Lie}(\mathcal F_{B(k)})\oplus (\text{Lie}(\mathcal G_{B(k)})\cap\text{Hom}(F^1[\frac{1}{p}],F^0[\frac{1}{p}]))$ of $B(k)$-vector spaces. Thus the rank of $\text{Lie}(\mathcal G_{B(k)})\cap\text{Hom}(F^1,F^0)$ is $\dim_{B(k)}(\text{Lie}(\mathcal G_{B(k)}))-\dim_{B(k)}(\text{Lie}(\mathcal F_{B(k)}))$ and therefore it is also equal to $\dim(\mathcal G_{B(k)}/\mathcal F_{B(k)})$. 

We will use the notation of the proof of Lemma \ref{L3} for a point $w\in\text{Sh}(G,\mathcal X)(\mathbb C)$ that lifts the $\mathbb C$-valued point of $\mathcal N_{E(G,\mathcal X)}$ defined naturally by the generic fibre of $z$ and the $O_{(v)}$-monomorphism $W(k)\hookrightarrow\mathbb C$. Let $W^{\vee}\otimes_{\mathbb Q} \mathbb C=F^{1,0}_w\oplus F^{0,1}_w$ be the Hodge decomposition  defined by $h_w\in\mathcal X$ (it is the dual of the Hodge decomposition of the property (ii) of the proof of Lemma \ref{L3}). We have a natural isomorphism $(M\otimes_{W(k)} \mathbb C,(t_{\alpha})_{\alpha\in\mathcal J})\rightarrow (W^{\vee}\otimes_{\mathbb Q} \mathbb C,(v_{\alpha})_{\alpha\in\mathcal J})$ that takes $F^1\otimes_{W(k)} \mathbb C$ to $F^{1,0}_w$, cf. Subsection B6 and Lemma \ref{L3} (b). Thus $\mathcal G_{B(k)}$ is a form of $G_{B(k)}$, $\mathcal F_{B(k)}$ is a parabolic subgroup of $\mathcal G_{B(k)}$, and we have $\dim(\mathcal G_{B(k)}/\mathcal F_{B(k)})=\dim(G_{\mathbb C}/P_w)$, where $P_w$ is the parabolic subgroup of $G_{\mathbb C}$ which is the normalizer of $F^{1,0}_w$ in $G_{\mathbb C}$. But $G_{\mathbb C}/P_w$ is the compact dual of any connected component of $\mathcal X$. Thus $\dim(G_{\mathbb C}/P_w)=d$ and therefore $\text{Lie}(\mathcal G_{B(k)})\cap\text{Hom}(F^1,F^0)$ has rank $d$.\endproof 
 
\begin{theorem} [Key Theorem] \label{T4} 
If $p=2$, then we assume that $D$ is a direct sum of connected and \'etale Barsotti--Tate groups (e.g., this holds if $G_{\mathbb Z_{(p)}}$ is a torus). We have:

\medskip
{\bf (a)} There exist isomorphisms
$$(M,(t_{\alpha})_{\alpha\in\mathcal J},\psi_M)\rightarrow (H^1_{\acute et}(A_{B(k)},\mathbb Z_p)\otimes_{\mathbb Z_p} W(k),(u_{\alpha})_{\alpha\in\mathcal J},\psi_{H^1_{\acute{et}}})\rightarrow (L^{\vee}_{(p)}\otimes_{\mathbb Z_{(p)}} W(k),(v_{\alpha})_{\alpha\in\mathcal J},\psi^{\vee}).$$
\indent
{\bf (b)} The group scheme $\mathcal G$ is isomorphic to $G_{W(k)}=G_{\mathbb Z_{(p)}}\times_{\Spec(\mathbb Z_{(p)})} \Spec(W(k))$.
\end{theorem}

\medskip\proof
From Theorem \ref{T9} applied to the triple $(D,\lambda_D,(t_{\alpha})_{\alpha\in\mathcal J})$ and from Formula (3) we get the existence of an isomorphism 
$$(M,(t_{\alpha})_{\alpha\in\mathcal J},\psi_M)\rightarrow (H^1_{\acute et}(A_{B(k)},\mathbb Z_p)\otimes_{\mathbb Z_p} W(k),(u_{\alpha})_{\alpha\in\mathcal J},\psi_{H^1_{\acute{et}}}).$$ 
Thus it suffices to prove part (a) under the extra assumption that $k$ has a countable transcendental degree, i.e., there exists an $E(G,\mathcal X)$-monomorphism $B(k)\hookrightarrow\mathbb C$. Let $w\in\mathcal N_{E(G,\mathcal X)}(\mathbb C)$ be the composite of the resulting morphism $\Spec(\mathbb C)\rightarrow\Spec(B(k))$ with the generic fibre of $z$. There exists a unit $\epsilon$ of $W(k)$ such that we have isomorphisms $(H^1_{\acute et}(A_{B(k)},\mathbb Z_p)\otimes_{\mathbb Z_p} W(k),(u_{\alpha})_{\alpha\in\mathcal J},\psi_{H^1_{\acute{et}}})\rightarrow (L^{\vee}_{(p)}\otimes_{\mathbb Z_{(p)}} W(k),(v_{\alpha})_{\alpha\in\mathcal J},\epsilon\psi^{\vee})$ (cf. Lemma \ref{L3} (a)). Thus part (a) follows once we remark that 
$\mu(W(k))(\epsilon)$ defines an isomorphism $(M,(t_{\alpha})_{\alpha\in\mathcal J},\psi_M)\rightarrow (M,(t_{\alpha})_{\alpha\in\mathcal J},\epsilon^{-1}\psi_M)$. 

Part (b) follows from part (a). \endproof

\begin{lemma}\label{L6} 
Let $G_{\mathbb Q_p}^{\text{v}}$ be a normal, reductive subgroup of $G_{\mathbb Q_p}$ such that there exists a cocharacter $\mathbb G_{m,B(k(v))}\rightarrow G_{B(k(v))}^{\text{v}}$ whose extension to $\mathbb C$ via an $O_{(v)}$-monomorphism $B(k(v))\hookrightarrow \mathbb C$ is $G(\mathbb C)$-conjugate to the cocharacters $\mu_h$ of $G_{\mathbb C}$ introduced in Subsection 1.3 ($h\in \mathcal X$). Let $\mathcal G_{B(k)}^{\text{v}}$ be the normal, reductive subgroup of $\mathcal G_{B(k)}$ which corresponds to $G_{\mathbb Q_p}^{\text{v}}$ via Fontaine comparison theory, cf. Lemmas \ref{L15} (a) and \ref{L3} (a). Then $\mu$ factors through the schematic closure $\mathcal G^{\text{v}}$ of $\mathcal G_{B(k)}^{\text{v}}$ in $\mathcal G$.
\end{lemma}

\medskip\proof
To prove this we can assume there exists a $W(k(v))$-monomorphism $W(k)\hookrightarrow\mathbb C$. We have canonical isomorphisms $(M\otimes_{W(k)} \mathbb C,(t_{\alpha})_{\alpha\in\mathcal J})\rightarrow (W^{\vee}\otimes_{\mathbb Q} \mathbb C,(v_{\alpha})_{\alpha\in\mathcal J})$
 such that $F^1\otimes_{W(k)} \mathbb C$ is mapped to the Hodge filtration of $W^{\vee}\otimes_{\mathbb Q} \mathbb C$ defined by a cocharacter $\mu_h:\mathbb G_{m,\mathbb C}\rightarrow G_{\mathbb C}$ with $h\in\mathcal X$ (see Subsection B6 and Lemma \ref{L3} (b)). We know that $\mu_{\mathbb C}$ is $\mathcal G(\mathbb C)$-conjugate to some (any) $\mu_h$, cf. Lemma \ref{L18}. From this and the very definition of $G_{\mathbb Q_p}^{\text{v}}$ we get that $\mu$ factors through $\mathcal G^{\text{v}}$.\endproof 

\subsection{Local deformation theory} Let $\mathcal G^\prime$ be the universal smoothening of $\mathcal G$, cf. Subsection A1. As $\mathcal G_{B(k)}=\mathcal G^\prime_{B(k)}$ is a form of $G_{B(k)}$, it is a reductive group over $B(k)$ of dimension $l$. Thus the relative dimension of $\mathcal G^\prime$ over $\Spec(W(k))$ is also $l$. Let $R$ be the completion of the local ring of $\mathcal G^\prime$ at the identity element of $\mathcal G^\prime_k$. Let $g_{\text{univ}}\in \mathcal G^\prime(R)$ be the natural (universal) element. Let $U$ be the connected, unipotent, smooth, closed subgroup scheme of either $\mathcal G$ of $\mathcal G^\prime$ whose Lie algebra is $\text{Lie}(\mathcal G_{B(k)})\cap\text{Hom}(F^1,F^0)$ (cf. Subsection B2 and Subsubsection B4.1). As the rank of $\text{Lie}(\mathcal G_{B(k)})\cap\text{Hom}(F^1,F^0)$ is $d$ (cf. Lemma \ref{L5}), $U$ is isomorphic to $\mathbb G_{a,W(k)}^d$. 

We fix an identification $R=W(k)[[x_1,\ldots,x_l]]$ such that the identity section of $\mathcal G^\prime$ is defined by the ideal $\mathfrak{I}:=(x_1,\ldots,x_l)$ of $R$. Let $\Phi_R$ be the Frobenius lift of $R$ that is compatible with $\sigma$ and we have $\Phi_R(x_i)=x_i^p$ for all $i\in\{1,\ldots,l\}$. The $\mathfrak{I}$-adic completion $\hat\Omega_{R/W(k)}$ of $\Omega_{R/W(k)}$ is a free $R$-module that has $\{dx_1,\ldots,dx_l\}$ as an $R$-basis. Let $d\Phi_R:\hat\Omega_{R/W(k)}\rightarrow \hat\Omega_{R/W(k)}$ be the differential map of $\Phi_R$. Let $M_R:=M\otimes_{W(k)} R$ and $F^1_R:=F^1\otimes_{W(k)} R$. We consider the $\Phi_R$-linear endomorphism
$$\Phi:=g_{\text{univ}}(\phi\otimes\Phi_R):M_R\rightarrow M_R.$$ 
Let $\nabla:M_R\rightarrow M_R\otimes_R \hat\Omega_{R/W(k)}$ be the unique connection on $M_R$ such that we have $\nabla\circ\Phi=(\Phi\otimes d\Phi_R)\circ\nabla$; it is integrable and nilpotent modulo $p$ (see Subsection B6).  See properties (i) of (iii) of Subsubsection B4.1 for three main properties of $\nabla$ and for the fact that there exists a unique $\text{Ker}(\mathbb G_{m,W(k)}(R)\rightarrow\mathbb G_{m,W(k)}(R/\mathfrak{I}))$-multiple $\psi_{M_R}$ of the perfect, alternating bilinear form $\psi_M$ on $M_R$ such that we have an identity $\psi_{M_R}(\Phi(a),\Phi(b))=p\Phi_R(\psi_{M_R}(a,b))$ for all $a, b\in M_R$. 

There exists a unique principally quasi-polarized Barsotti--Tate group $(D_R,\lambda_{D_R})$ over $R$ which modulo $\mathfrak{I}$ is $(D,\lambda_D)$ and whose principally quasi-polarized filtered $F$-crystal over $R/pR$ is the quintuple $(M_R,F^1_R,\Phi,\nabla,\psi_{M_R})$, cf. Lemmas \ref{L16} and \ref{L17}. 

Let $(B_R,\lambda_{B_R})$ be the principally polarized abelian scheme over $R$ which modulo $\mathfrak{I}$ is $(A,\lambda_A)$ and whose principally quasi-polarized Barsotti--Tate group is $(D_R,\lambda_{D_R})$, cf. Serre--Tate deformation theory and Grothendieck algebraization theorem. Let 
$$\tau_R:\Spec(R)\rightarrow\mathcal M_r$$ 
be the natural morphism that corresponds to $(B_R,\lambda_{B_R})$ and its level-$N$ symplectic similitude structures which lift those of $(A,\lambda_A)$ (here $N\in\mathbb N\setminus (p\mathbb N\cup\{1,2\})$). We have a canonical identification $H^1_{\text{dR}}(B_R/R)=M_R=M\otimes_{W(k)} R$, cf. [2, Ch. V, Subsect. 2.3] and [4, Prop. 2.5.8]. Under this identification, the following two properties hold:

\medskip\noindent
{\bf (i)} the perfect form on $M_R$ defined by the principal polarization $\lambda_{B_R}$ of $B_R$ is identified with $\psi_{M_R}$;

\smallskip\noindent
{\bf (ii)} for all $s\in\mathbb N^*$, the connection on $H^1_{\text{dR}}(B_R/R)/\mathfrak{I}^sH^1_{\text{dR}}(B_R/R)=M_R/\mathfrak{I}^sM_R$ induced by $\nabla$ is the Gauss--Manin connection of $B_R\times_{\Spec(R)} \Spec(R/\mathfrak{I}^s)$ (cf. [2, Ch. V, Prop. 3.6.4] and the fact that $R/\mathfrak{I}^s$ is $p$-adically complete).

\begin{theorem} [Faltings] \label{T5}
For each $\alpha\in\mathcal J$, the tensor $t_{\alpha}\in\mathcal T(M)\otimes_{W(k)} R[\frac{1}{p}]=\mathcal T(M_R)[\frac{1}{p}]$ is the de Rham component of a Hodge cycle on $B_{R[\frac{1}{p}]}$.
\end{theorem}

\medskip\proof
We recall that $B_R$ is a deformation of $A$ over $R$. As $t_{\alpha}\in\mathcal T(M)[\frac{1}{p}]$ is the de Rham component of the Hodge cycle $w_{\alpha}$ on $A_{B(k)}$ and due to the property (i) of Subsubsection B4.1, the theorem is a result of Faltings whose essence is outlined in [54, Rm. 4.1.5] and whose proof is presented here. 

As $\mathcal A_{r,1,N}$ is a quasi-projective $\mathbb Z_{(p)}$-scheme and as the set $\mathcal J$ is countable, it suffices to prove the theorem in the case when there exists a morphism $e_k:\Spec(\mathbb C)\rightarrow\Spec(W(k))$ of $\Spec(W(k(v)))$-schemes. We will view $\mathbb C$ as a $W(k)$-algebra via $e_k$. Let $\mathcal R:=\mathbb C[[x_1,\ldots,x_l]]$ and $\mathcal S:=\mathbb C[[x_1,\ldots,x_d]]$. Let $\mathcal I:=\mathfrak{I}\mathcal R=(x_1,\ldots,x_l)\mathcal R$ and $\mathcal I_0$ be the maximal ideals of $\mathcal R$ and $\mathcal S$ (respectively). 

Let $(B_{\mathcal R},(t_{\alpha})_{\alpha\in\mathcal J},\lambda_{B_{\mathcal R}})$ be the pullback of $(B_R,(t_{\alpha})_{\alpha\in\mathcal J},\lambda_{B_R})$ via the morphism of schemes defined by the natural $W(k)$-monomorphism 
$$R=W(k)[[x_1,\ldots,x_l]]\hookrightarrow \mathbb C[[x_1,\ldots,x_l]]=\mathcal R.$$ To prove the theorem it suffices to show that the tensor $t_{\alpha}\in \mathcal T(M)\otimes_{W(k)} \mathcal R=\mathcal T(M_R\otimes_R \mathcal R)=\mathcal T(H^1_{\text{dR}}(B_{\mathcal R}/\mathcal R))$ is the de Rham component of a Hodge cycle on $B_{\mathcal R}$.

Let $(C_{\mathcal S},(w_{\alpha}^{\mathcal S})_{\alpha\in\mathcal J},\lambda_{C_{\mathcal S}})$ be the pullback of $(\mathcal A,(w_{\alpha}^{\mathcal A})_{\alpha\in\mathcal J},\lambda_{\mathcal A})$ via a formally \'etale morphism $\Spec(\mathcal S)\rightarrow\mathcal N^{\text{s}}_{{\mathbb C}}$ whose composite with the closed embedding $\Spec(\mathbb C)\hookrightarrow\Spec(\mathcal S)$ is the point $z\circ e_k\in\mathcal N_{\mathbb C}(\mathbb C)=\mathcal N^{\text{s}}(\mathbb C)$. Let $\mathcal W:=H^1_{\text{dR}}(C_{\mathcal S}/\mathcal S)$. Let $\psi_{\mathcal W}$ be the perfect, alternating bilinear form on $\mathcal W$ defined by $\lambda_{C_{\mathcal S}}$. Let $t_{\alpha}^{\mathcal S}\in\mathcal T(\mathcal W)$ be the de Rham component of $w_{\alpha}^{\mathcal S}$. Let $\Delta$ be the Gauss--Manin connection on $\mathcal W$ defined by $C_{\mathcal S}$. 
 We recall that $\psi^{\vee}$ is the alternating bilinear form on $W^{\vee}$ (or on $L_{(p)}^{\vee}$) defined naturally by $\psi$. 
 
From Corollary \ref{C2} and (the proof of) Lemma \ref{L3} (b) we get that there exists $\epsilon\in\mathbb Q\setminus\{0\}$ for which there exist an isomorphism
$$I:(\mathcal W,(t_{\alpha}^{\mathcal S})_{\alpha\in\mathcal J},\psi_{\mathcal W})\rightarrow (W^{\vee}\otimes_{\mathbb Q} \mathcal S,(v_{\alpha})_{\alpha\in\mathcal J},\epsilon\psi^{\vee})$$ 
under which the $\mathfrak{I}_0$-completion of $\Delta$ becomes the $\mathfrak{I}_0$-completion of the flat connection on $W^{\vee}\otimes_{\mathbb Q} \mathcal S$ that annihilates $W^{\vee}\otimes 1$. As there exist isomorphisms of $(W^{\vee}\otimes_{\mathbb Q} \mathbb C,(v_{\alpha})_{\alpha\in\mathcal J})$ that take $\psi^{\vee}$ to $\epsilon\psi^{\vee}$, we can assume that $\epsilon=1$. We fix an isomorphism $I$ with $\epsilon=1$ and we view it as an identification. For each $\beta\in\mathbb G_{m,\mathbb C}(\mathcal R)$, there exist isomorphisms of $(W^{\vee}\otimes_{\mathbb Q} \mathcal R,(v_{\alpha})_{\alpha\in\mathcal J})$ that take $\psi^{\vee}$ to $\beta\psi^{\vee}$. Thus, based on the construction of $M_R$ and on either Lemma \ref{L3} (b) or the proof of Lemma \ref{L5}, we get that there exist isomorphisms
$$I_A:(M_R\otimes_R \mathcal R,(t_{\alpha})_{\alpha\in\mathcal J},\psi_{M_R})\rightarrow (W^{\vee}\otimes_{\mathbb Q} \mathcal R,(v_{\alpha})_{\alpha\in\mathcal J},\psi^{\vee}).$$
\indent
By induction on $s\in\mathbb N^*$ we show that there exists a unique morphism of $\mathbb C$-schemes
$$J_s:\Spec(\mathcal R/\mathcal I^s)\rightarrow\Spec(\mathcal S)$$ 
that has the following property:

\medskip\noindent
{\bf (i)} {\it There exists an isomorphism $\xi_s$ between the reduction of $(B_{\mathcal R},(t_{\alpha})_{\alpha\in\mathcal J},\lambda_{B_{\mathcal R}})$ modulo $\mathcal I^s$ and the pullback $J_s^*((C_{\mathcal S},(t_{\alpha}^{\mathcal S})_{\alpha\in\mathcal J},\lambda_{C_{\mathcal S}}))$  which modulo $\mathcal I/\mathcal I^s$ is defined by $1_{A_{\mathbb C}}=1_{C_{\mathcal S}\times_{\Spec(\mathcal S)} \Spec(\mathbb C)}=1_{B_{\mathcal R}\times_{\Spec(\mathcal R)} \Spec(\mathbb C)}$.} 

\medskip
As $\mathcal N^{\text{s}}_{E(G,\mathcal X)}$ is a closed subscheme of $\mathcal M_{r,E(G,\mathcal X)}$ (cf. Fact \ref{F1}) and as $\Spec(\mathcal S)\rightarrow \mathcal N_{{\mathbb C}}^{\text{s}}$ is formally \'etale, the deformation $(C_{\mathcal S},\lambda_{C_{\mathcal S}})$ of the principally polarized abelian variety $(A,\lambda_A)_{\mathbb C}$ is versal, i.e., the Kodaira--Spencer map $\mathfrak{K}$ of $\Delta$ is injective and its image is a free $\mathcal S$-module of rank $d$ which is a direct summand of its codomain. This implies the uniqueness of $J_s$. 

The existence of $J_1$ is obvious. For $s\ge 2$ the passage from the existence of $J_{s-1}$ to the existence of $J_s$ goes as follows. Let $J_s^\prime:\Spec(\mathcal R/\mathcal I^s)\rightarrow\Spec(\mathcal S)$ be an arbitrary morphism of $\mathbb C$-schemes that lifts $J_{s-1}$. Let $\Delta_s$ be the connection on $\mathcal W\otimes_{\mathcal S} \mathcal R/\mathcal I^s=W^{\vee}\otimes_{\mathbb Q} \mathcal R/\mathcal I^s$ which is the extension of the connection $\Delta$ on $\mathcal W$ via $J_s^\prime$ (the last identification is defined naturally by $I$). Let $\nabla_s$ be the Gauss--Manin connection on $H^1_{\text{dR}}(B_R/R)\otimes_R \mathcal R/\mathcal I^s=M_R\otimes_R \mathcal R/\mathcal I^s$ defined by $B_R\times_{\Spec(R)} \Spec(\mathcal R/\mathcal I^s)$; it is the extension of the connection $\nabla$ on $M_R$ (cf. property (ii) of Subsection 3.2) and thus it annihilates each tensor $t_{\alpha}\in \mathcal T(M_R)\otimes_R \mathcal R/\mathcal I^s$ (cf. property (i) of Subsubsection B4.1). From Lemma \ref{L4} we get:

\medskip\noindent
{\bf (ii)} {\it There exists a unique isomorphism $I_{A,s}:M_R\otimes_R \mathcal R/\mathcal I^s\rightarrow W^{\vee}\otimes_{\mathbb Q} \mathcal R/\mathcal I^s$ which lifts a fixed isomorphism between $(M_R\otimes_R \mathcal R\otimes_{\mathcal R} \mathcal R/\mathcal I,(t_{\alpha})_{\alpha\in\mathcal J})=(H^1_{\text{dR}}(A_{\mathbb C}/\mathbb C),(t_{\alpha})_{\alpha\in\mathcal J})$ and $(W^{\vee}\otimes_{\mathbb Q} \mathbb C,(v_{\alpha})_{\alpha\in\mathcal J})$ obtained as in Lemma \ref{L3} (b) and such that under it $\nabla_s$ becomes the flat connection $\delta_s$ on $W^{\vee}\otimes_{\mathbb C} \mathcal R/\mathcal I^s$ that annihilates $W^{\vee}\otimes 1$}. 

\medskip
We denote also by $I_{A,s}$ the isomorphism $\mathcal T(M_R\otimes_R \mathcal R/\mathcal I^s)\rightarrow \mathcal T(W^{\vee}\otimes_{\mathbb Q} \mathcal R/\mathcal I^s)$ induced by $I_{A,s}$. As $I_{A,s}(t_{\alpha})$ and $v_{\alpha}$ are two tensors of $W^{\vee}\otimes_{\mathbb C} \mathcal R/\mathcal I^s$ that are annihilated by $\delta_s$ and that coincide modulo $\mathcal I/\mathcal I^s$, we get that we have $I_{A,s}(t_{\alpha})=v_{\alpha}$ for all $\alpha\in\mathcal J$. An argument similar to the one above involving $\epsilon\in\mathbb Q\setminus\{0\}$ shows that we can assume that $I_{A,s}$ takes $\psi_{M_R}$ to $\psi^{\vee}$. Thus we can choose $I_A$ such that it lifts $I_{A,s}$.  We will view the reduction $I_{A,s}$ of $I_A$ modulo $\mathcal I^s$ as an identification. Therefore we will also identify $\nabla_s=\delta_s$.

From the existence of $I$ and the fact that $I_{A,s}$ is the reduction of $I_A$ modulo $\mathcal I^s$, we get that there exists an isomorphism 
$$\zeta_s:J_s^{\prime*}((\mathcal W,(t_{\alpha}^{\mathcal S})_{\alpha\in\mathcal J},\psi_{\mathcal W}))=(W^{\vee}\otimes_{\mathbb Q} \mathcal R/\mathcal I^s,(v_{\alpha})_{\alpha\in\mathcal J},\psi^{\vee})\rightarrow$$ 
$$\rightarrow (M_R\otimes_R \mathcal R/\mathcal I^s,(t_{\alpha})_{\alpha\in\mathcal J},\psi_{M_R})=(W^{\vee}\otimes_{\mathbb Q} \mathcal R/\mathcal I^s,(v_{\alpha})_{\alpha\in\mathcal J},\psi^{\vee})$$ 
with the properties that it lifts the identity automorphism of $W^{\vee}\otimes_{\mathbb Q} \mathbb C$ and that:

\medskip\noindent
{\bf (iii)} {\it It respects the Gauss--Manin connections, i.e., it takes $\Delta_s$ to $\nabla_s=\delta_s$.}

\medskip
From the uniqueness part of the property (ii) we also get that:

\medskip\noindent
{\bf  (iv)} {\it The reduction of $\zeta_s$ modulo $\mathcal I^{s-1}$ is the isomorphism defined by $\xi_{s-1}$.}

\medskip
Let $F_{A,s}^1$ and $F_{C,s}^1$ be the Hodge filtrations of $W^{\vee}\otimes_{\mathbb Q} \mathcal R/\mathcal I^s$ defined naturally by $B_{\mathcal R}$ and $J_s^{\prime*}(C_{\mathcal S})$ (respectively) via the above identifications. The direct summands $F^1_{A,s}$ and $\zeta_{s}(F^1_{C,s})$ of $W^{\vee}\otimes_{\mathbb Q} \mathcal R/\mathcal I^s$ coincide modulo $\mathcal I^{s-1}/\mathcal I^s$, cf. property (iv). Moreover, there exist direct sum decompositions 
$$W^{\vee}\otimes_{\mathbb Q} \mathcal R/\mathcal I^s=F^1_{A,s}\oplus F^0_{A,s}=F^1_{C,s}\oplus F^0_{C,s}$$ 
defined by cocharacters $\mu_{A,s}$ and $\mu_{C,s}$ of the reductive subgroup scheme $G_{\mathcal R/\mathcal I^s}$ of $\pmb{\text{GL}}_{W^{\vee}\otimes_{\mathbb Q} \mathcal R/\mathcal I^s}$ (here $\mathbb G_{m,\mathcal R/\mathcal I^s}$ through $\mu_{*,s}$ fixes $F^0_{*,s}$ and acts via the weight $-1$ on $F^1_{*,s}$).  The existence of $\mu_{A,s}$ is a direct consequence of the existence of the cocharacter $\mu:\mathbb G_{m,W(k)}\rightarrow\mathcal G$ (see paragraph before Lemma \ref{L5}) and of the definition of $F^1_R$ (see Subsection 3.2) while the existence of $\mu_{C,s}$ is well-known. 

As $F^1_{A,s}$ and $\zeta_{s}(F^1_{C,s})$ coincide modulo $\mathcal I^{s-1}/\mathcal I^s$, we can choose $\mu_{A,s}$ and $\mu_{C,s}$ such that $\zeta_s^{-1}\mu_{A,s}\zeta_s$ and $\mu_{C,s}$ commute modulo $\mathcal I^{s-1}/\mathcal I^s$ and thus coincide modulo $\mathcal I^{s-1}/\mathcal I^s$. Thus based on [15, Vol. II, Exp. IX, Thm. 3.6], there exists an element $g_s\in\text{Ker}(G(\mathcal R/\mathcal I^s)\rightarrow G(\mathcal R/\mathcal I^{s-1}))$ such that we have $\zeta_s^{-1}\mu_{A,s}\zeta_s=g_s\mu_{C,s}g_s^{-1}$. We have $\zeta_s(g_s(F^1_{C,s}))=F^1_{A,s}$. We would like to mention that the original approach of Faltings used the strictness of filtrations of morphisms between Hodge $\mathbb R$-structures in order to get the existence of the element $g_s$.

The image of $\mathfrak{K}$ is a free $\mathcal S$-module that has rank $d$ and that is equal to the image of $\text{Lie}(G_{\mathcal S})$ into the codomain of $\mathfrak{K}$. Thus we can replace $J_s^\prime$ by another morphism $J_s:\Spec(\mathcal R/\mathcal I^s)\rightarrow\Spec(\mathcal S)$ that lifts $J_{s-1}$ and such that under it and $I_{A,s}$ the Hodge filtration $F^1_{C,s}$ gets replaced by $g_s(F^1_{C,s})$. Thus $\zeta_s$ becomes the de Rham realization of an isomorphism $\xi_s$ between the reduction of $(B_{\mathcal R},(t_{\alpha})_{\alpha\in\mathcal J},\lambda_{B_{\mathcal R}})$ modulo $\mathcal I^s$ and $J_s^*((C_{\mathcal S},(t_{\alpha}^{\mathcal S})_{\alpha\in\mathcal J},\lambda_{C_{\mathcal S}}))$ which lifts $\xi_{s-1}$, cf. deformation theory of abelian varieties. Thus the morphism $J_s$ has the desired properties. This ends the induction.  

Let $J_{\infty}:\Spec(\mathcal R)\rightarrow\Spec(\mathcal S)$ be the morphism defined by $J_s$'s ($s\in\mathbb N^*$). The isomorphism $\xi_s$ is uniquely determined by the property (i) and this implies that $\xi_{s+1}$ lifts $\xi_s$. From this and Grothendieck algebraization theorem we get the existence of an isomorphism 
$$\xi_{\infty}:(B_{\mathcal R},(t_{\alpha})_{\alpha\in\mathcal J},\lambda_{B_{\mathcal R}})\rightarrow J_{\infty}^*((C_{\mathcal S},(t_{\alpha}^{\mathcal S})_{\alpha\in\mathcal J},\lambda_{C_{\mathcal S}}))$$ 
which modulo $\mathcal I$ is defined by $1_{A_{\mathbb C}}$ and which lifts each $\xi_s$ with $s\in\mathbb N$. Thus for each $\alpha\in\mathcal J$, the tensor $t_{\alpha}\in \mathcal T(M)\otimes_{W(k)} \mathcal R$ is the de Rham component of the Hodge cycle on $B_{\mathcal R}$ which is the pullback of the Hodge cycle $J_{\infty}^*(w_{\alpha}^{\mathcal S})$ on $J_{\infty}^*(C_{\mathcal S})$ via the isomorphism $B_{\mathcal R}\rightarrow J_{\infty}^*(C_{\mathcal S})$ that defines $\xi_{\infty}$.\endproof
 
\subsection{Proof of Theorem \ref{T1}}
In this subsection we prove the Basic Theorem \ref{T1}. Let $O$ be an $O_{(v)}$-algebra which is a discrete valuation ring of absolute ramification index $1$. We choose the field $k$ such that we have a $O_{(v)}$-monomorphism $O\hookrightarrow W(k)$. Let $Z$ be a regular, formally smooth $O$-scheme equipped with a morphism $\chi:Z_{E(G,\mathcal X)}\rightarrow\text{Sh}_H(G,\mathcal X)=\mathcal N^{\text{s}}_{E(G,\mathcal X)}$. Thus $\chi$ extends uniquely to a morphism $\chi_Z:Z\rightarrow\mathcal N$, cf. Proposition \ref{P2} (c). To prove Theorem \ref{T1} (a) we have to show that $\chi_Z$ factors through $\mathcal N^{\text{s}}$. It suffices to check this under the extra assumptions that $O=W(k)$ and $Z=\Spec(O)$. We will use the notation of Subsection 3.1 for the point
$z:=\chi_Z\in\mathcal N(W(k))$.  

Let $y:\Spec(k)\hookrightarrow\mathcal N_{W(k)}$ be the closed embedding defined by the special fibre of $z\in\mathcal N(W(k))$. Let $O_y^{\text{bigg}}$ (resp. $O_{y}$) be the completion of the local ring of $y$ in $\mathcal M_{r,W(k)}$ (resp. in $\mathcal N_{W(k)}$). As $\mathcal Q$ is a normal, flat $O_{(v)}$-scheme of relative dimension $d$ and as $\mathcal N$ is a pro-finite pro-\'etale cover of $\mathcal Q$ (cf. Proposition \ref{P2} (a)), the local ring $O_{y}$ is normal and noetherian of dimension $1+d$. The natural homomorphism $n_y:O_y^{\text{bigg}}\rightarrow O_y$ (associated at $y$ to the morphism $\mathcal N_{W(k)}\rightarrow \mathcal M_{r,W(k)}$) is finite, cf. Proposition \ref{P2} (b). Let $h_z:O_y^{\text{bigg}}\rightarrow R$ be the $W(k)$-homomorphism that defines $\tau_R:\Spec(R)\rightarrow\mathcal M_r$.
 
Let $S:=W(k)[[x_1,\ldots,x_d]]$ and let $\mathfrak{I}_0:=(x_1,\ldots,x_d)$ be its ideal. We consider a closed embedding $c_R:\Spec(S)\hookrightarrow\Spec(R)$ such that the following two properties hold (cf. Subsubsection B4.2 and Lemma \ref{L5}):

\medskip\noindent
{\bf (i)} it is defined by a $W(k)$-epimorphism $e_z:R\twoheadrightarrow S$ satisfying $e_z(\mathfrak{I})\subset \mathfrak{I}_0\subset S$;

\smallskip\noindent
{\bf (ii)} the pullback of $(M_R,F^1_R,\Phi,\nabla,\psi_{M_R})$ via the closed embedding $\Spec(S/pS)\hookrightarrow\Spec(R/pR)$, is a principally quasi-polarized filtered $F$-crystal over $S/pS$ whose Kodaira--Spencer map is injective and has an image equal to the direct summand $\text{Lie}(U)\otimes_{W(k)} S$ of $\text{Hom}(F^1,F^0)\otimes_{W(k)} S\simeq\text{Hom}(F^1,M/F^1)\otimes_{W(k)} S$. 

\medskip
From the property (ii) we get that the composite morphism $\tau_S:=\tau_R\circ c_R:\Spec(S)\rightarrow\mathcal M_r$ is defined naturally by a $W(k)$-epimorphism $s_z^{\text{bigg}}:=e_z\circ h_z:O_y^{\text{bigg}}\twoheadrightarrow S$. 

The existence of the isomorphism $\xi_{\infty}$ (see the end of the proof of Theorem \ref{T5}) implies that the morphism $\tau_R:\Spec(R)\rightarrow\mathcal M$ factors through $\mathcal N$ in such a way that modulo the ideal $\mathfrak{I}$ of $R$ it defines the point $z\in\mathcal N(W(k))$. Thus there exists a $W(k)$-epimorphism $s_z:O_y\twoheadrightarrow S$ such that we have $s_z^{\text{bigg}}=s_z\circ n_y$, i.e., the following diagram is commutative 
\[\xymatrix{
O_y^{\text{bigg}} \ar[r]^{n_y} \ar[d]^{h_Z} & O_y \ar[d]^{s_z} \\
R \ar[r]^{e_Z} & S.
}\]
\noindent
By reasons of dimensions of local, noetherian, normal rings, the $W(k)$-epimorphism $s_z:O_y\twoheadrightarrow S$ is an isomorphism. Thus $\mathcal N_{W(k)}$ is formally smooth at $z$ and therefore $z$ factors through $\mathcal N^{\text{s}}$. Therefore Theorem \ref{T1} (a) holds and $y$ is a $k$-valued point of $\mathcal N^{\text{s}}_{W(k)}$.

As $s_z$ is an isomorphism, the $W(k)$-homomorphism $n_y:O_y^{\text{bigg}}\rightarrow O_y$ is onto. Therefore the natural $W(k)$-morphism $\mathcal N^{\text{s}}_{W(k)}\rightarrow\mathcal M_{r,W(k)}$ is a formally closed embedding at $y\in \mathcal N^{\text{s}}_{W(k)}(k)$. The role of $z\in\mathcal N(W(k))$ is that of an arbitrary $W(k)$-valued of $\mathcal N$ (and thus due to Theorem \ref{T1} (a)) of $\mathcal N^{\text{s}}$. Thus the $W(k)$-morphism $\mathcal N^{\text{s}}_{W(k)}\rightarrow\mathcal M_{r,W(k)}$ is a formally closed embedding at every $k$-valued point of $\mathcal N^{\text{s}}_{W(k)}$. Thus Theorem \ref{T1} (b) also holds. 

We check that the Theorem \ref{T1} (c) holds. Let $Z$ be a smooth $O_{(v)}$-scheme such that we have a morphism $\chi:Z_{E(G,\mathcal X)}\rightarrow\text{Sh}_{H\times H^{(P)}}(G,\mathcal X)$. From Proposition \ref{P2} (b) and Lemma \ref{L2} we get that $\mathcal N/H^{(p)}$ has a finite \'etale cover which is projective; thus $\mathcal N/H^{(p)}$ is a proper $O_{(v)}$-scheme. From this and the valuative criterion of properness, we get that there exists an open subscheme $U_Z$ of $Z$ such that it contains $Z_{E(G,\mathcal X)}$, the complement of $U_Z$ in $Z$ has codimension in $Z$ at least $2$, and the morphism $\chi$ extends uniquely to a morphism $\chi_{U_Z}:U_Z\rightarrow\mathcal N/H^{(p)}$. From the classical purity theorem of Zariski, Nagata and Grothendieck (see [21, Exp. X, Thm. 3.4 (i)]) we get that the pro-finite pro-\'etale cover $U_Z\times_{\mathcal N/H^{(p)}} \mathcal N\rightarrow U_Z$ extends uniquely to a pro-finite pro-\'etale cover $Z_{\infty}\rightarrow Z$. From this and Theorem \ref{T1} (a) we get that the natural morphism $U_Z\times_{\mathcal N/H^{(p)}} \mathcal N\rightarrow\mathcal N$ extends uniquely to a morphism $Z_{\infty}\rightarrow\mathcal N$. This implies that the morphism $\chi$ extends uniquely to a morphism $\chi_Z:Z\rightarrow\mathcal N/H^{(p)}$. Thus $\mathcal N/H^{(p)}$ is a N\'eron model of its generic fibre $\text{Sh}_{H\times H^{(p)}}(G,\mathcal X)$ over $O_{(v)}$, i.e., Theorem \ref{T1} (c) holds. This ends the proof of the Basic Theorem \ref{T1}. \hspace*\fill $\Box$
\subsection{de Rham realizations of Hodge cycles}
We denote also by $\tau_R$ the factorization of $\tau_R:\Spec(R)\rightarrow\mathcal M_r$ through either $\mathcal N$ or (cf. Theorem \ref{T1} (a)) $\mathcal N^{\text{s}}$ which modulo $\mathfrak{I}$ is the $W(k)$-valued point $z\in\mathcal N(W(k))=\mathcal N^{\text{s}}(W(k))$. As $s_z:O_y\rightarrow S$ is a $W(k)$-isomorphism and as we have a $W(k)$-epimorphism $e_z:R\twoheadrightarrow S$, the morphism $\tau_R:\Spec(R)\rightarrow\mathcal N^{\text{s}}$ is formally smooth. Under the canonical identification $H^1_{\text{dR}}(B_R/R)=M_R=M\otimes_{W(k)} R$, the pullback of $w_{\alpha}^{\mathcal A}$ via the morphism $\Spec(R[\frac{1}{p}])\rightarrow\mathcal N_{E(G,\mathcal X)}=\text{Sh}_H(G,\mathcal X)$ defined by $\tau_R$, is a Hodge cycle on $B_{R[\frac{1}{p}]}$ whose de Rham component $t_{\alpha}^{\prime}\in\mathcal T(M)\otimes_{W(k)} R[\frac{1}{p}]$  modulo $\mathfrak{I}[\frac{1}{p}]$ is $t_{\alpha}$ modulo $\mathfrak{I}[\frac{1}{p}]$. In fact we have $t_{\alpha}^{\prime}=t_{\alpha}$ for all $\alpha\in\mathcal J$. This follows either from the existence of $\xi_{\infty}$ (see end of the proof of Theorem \ref{T5}) or (in Faltings' approach) from the fact that there exists no non-trivial tensor of $\mathcal T(M)\otimes_{W(k)} \mathfrak{I}[\frac{1}{p}]$ fixed by $\Phi$. Similarly, the de Rham realization of the pullback of $\lambda_{\mathcal A}$ via the morphism $\Spec(R[\frac{1}{p}])\rightarrow\mathcal N_{E(G,\mathcal X)}=\text{Sh}_H(G,\mathcal X)$ defined by $\tau_R$, is $\psi_{M_R}$.

\subsection{The open subscheme $\mathcal N^{\text{m}}$ of $\mathcal N^{\text{s}}$} 

For $p>2$ let $\mathcal N^{\text{m}}:=\mathcal N^{\text{s}}$. If $p=2$, then let $\mathcal N^{\text{m}}$ be the maximal open subscheme of $\mathcal N^{\text{s}}$ with the property that for each algebraically closed field $k$ of characteristic $p$ and for every $z\in\mathcal N^{\text{m}}(W(k))$, Theorem \ref{T4} (a) (and thus also Theorem \ref{T4} (b)) holds. Thus regardless of the parity of $p$, for each such field $k$ and for every $z\in\mathcal N^{\text{m}}(W(k))$, Theorem \ref{T4} (a) holds. 

\begin{proposition}\label{P3}
The following two properties hold:

\medskip
{\bf (a)} Always $\mathcal N^{\text{m}}$ is a $G(\mathbb A_f^{(p)})$-invariant, open subscheme of $\mathcal N^{\text{s}}$.  

\smallskip
{\bf (b)} If Theorem \ref{T4} (a) holds for $z\in\mathcal N^{\text{s}}(W(k))$, then $z\in\mathcal N^{\text{m}}(W(k))$.
\end{proposition}

\medskip\proof  
The right translations of $z$ by elements of $G(\mathbb A_f^{(p)})$ correspond to passages to isogenies prime to $p$ of the abelian scheme A. Therefore the triple $(M,\phi,(t_{\alpha})_{\alpha\in\mathcal J})$ depends only on the $G(\mathbb A_f^{(p)})$-orbit of $z$. Thus, if Theorem \ref{T4} (a) holds for $z$, then Theorem \ref{T4} (a) also holds for every point in the $G(\mathbb A_f^{(p)})$-orbit of $z$. This implies that part (a) holds. 

To check part (b), let $H^{(p)}$, $\mathcal Q$, and $\mathcal Q^{\text{s}}$ be as in Subsection 2.2. By enlarging $N$ we can assume that the triple $(\mathcal A,(w^{\mathcal A}_{\alpha})_{\alpha\in\mathcal J},\lambda_{\mathcal A})$ is the pullback of an analogue triple $\mathcal T$ over $\mathcal Q$. Let $\Spec(Q)$ be an affine, open subscheme of $\mathcal Q^{\text{s}}$ such that the composite $z_{H^{(p)}}$ of $z:\Spec(W(k))\rightarrow\mathcal N^{\text{s}}$ with $\mathcal N^{\text{s}}\rightarrow\mathcal Q^{\text{s}}$ factors through $\Spec(Q)$. Let $(M_Q,(t_{\alpha}^Q)_{\alpha\in\mathcal J},\psi_{M_Q})$ be the de Rham realization of the pullback triple $\mathcal T_Q$ (of $\mathcal T$ to $\Spec(Q)$). Let $F^1_Q$ be the direct summand of $M_Q$ which is the Hodge filtration associated to $\mathcal T_Q$. By shrinking $\Spec(Q)$, we can assume that $M_Q$ and $F^1_Q$ are free $Q$-module of ranks $2r$ and $r$ (respectively). The existence of the formally smooth morphism $\tau_R:\Spec(R)\rightarrow\mathcal N^{\text{s}}$ implies that we have isomorphisms (cf. Subsection 3.4) 
$$(M_Q\otimes_{Q} R,F^1_Q\otimes_{Q} R,(t_{\alpha}^Q)_{\alpha\in\mathcal J},\psi_{M_Q})\rightarrow (M_R,F^1_R,(t_{\alpha})_{\alpha\in\mathcal J},\psi_{M_R})$$
$$=(M\otimes_{W(k)} R,F^1\otimes_{W(k)} R,(t_{\alpha})_{\alpha\in\mathcal J},\epsilon_M\psi_M),\leqno (4)$$
where $\epsilon_M\in\text{Ker}(\mathbb G_{m,W(k)}(R)\rightarrow\mathbb G_{m,W(k)}(R/\mathfrak{I}))$. We note that $\mu(R)(\epsilon_M)$ defines an isomorphism $(M\otimes_{W(k)} R,F^1\otimes_{W(k)} R,(t_{\alpha})_{\alpha\in\mathcal J},\psi_M)\rightarrow (M\otimes_{W(k)} R,F^1\otimes_{W(k)} R,(t_{\alpha})_{\alpha\in\mathcal J},\epsilon_M\psi_M)$. As Theorem \ref{T4} (a) holds for $z\in\mathcal N^{\text{s}}(W(k))$, we also have isomorphisms $(M\otimes_{W(k)} R,(t_{\alpha})_{\alpha\in\mathcal J},\psi_M)\rightarrow (L^{\vee}_{(p)}\otimes_{\mathbb Z_{(p)}} R,(v_{\alpha})_{\alpha\in\mathcal J},\psi^{\vee})$.
From the last three sentences and Artin approximation theorem ([7, Ch. 3, Sect. 3.6, Thm. 16]) we get that there exists a smooth, affine morphism $\Spec(Q^\prime)\rightarrow\Spec(Q)$ through which $z_{H^{(p)}}:\Spec(W(k))\rightarrow\Spec(Q)$ and the natural factorization $\Spec(R)\rightarrow\Spec(Q)$ of $\tau_R$ factor naturally producing morphisms $z_{H^{(p)}}^\prime:\Spec(W(k))\rightarrow\Spec(Q^\prime)$ and $\Spec(R)\rightarrow \Spec(Q^\prime)$ and such that we have an isomorphism 
$$(M_Q\otimes_{Q} Q^\prime,(t_{\alpha}^Q)_{\alpha\in\mathcal J},\psi_{M_Q})\rightarrow (L^{\vee}_{(p)}\otimes_{\mathbb Z_{(p)}} Q^\prime,(v_{\alpha})_{\alpha\in\mathcal J},\psi^{\vee})$$
whose extension to $R$ (via $\Spec(R)\rightarrow \Spec(Q^\prime)$) defines (4).
The image $\text{Im}(\Spec(Q^\prime)\rightarrow\Spec(Q))$ is an open subscheme of $\mathcal Q^{\text{s}}$ whose pullback to $\mathcal N^{\text{s}}$ is (due to the last isomorphism) an open subscheme of $\mathcal N^{\text{m}}$ that contains the point $z\in\mathcal N^{\text{s}}(W(k))$. Thus part (b) holds.\endproof

\subsubsection{Proof of Corollary \ref{C1}}
We assume that $G_{\mathbb Z_{(p)}}$ is smooth over $\mathbb Z_{(p)}$ and that the $k(v)$-scheme $\mathcal N^{\text{m}}_{k(v)}$ is non-empty. We fix a connected component $\mathcal C^{\text{m}}$ of $\mathcal N_k^{\text{m}}/H^{(p)}$ and consider a point $y_1\in \mathcal C^{\text{m}}(k)$. Not to introduce extra notation, we will assume that $z\in\mathcal N^{\text{m}}(W(k))$ is such that its image $[\mathcal N^{\text{m}}/H^{(p)}](W(k))$ lifts the point $y_1\in\mathcal C^{\text{m}}(k)$. As we will need to vary $y_1$, we will denote $Q_1=Q$ and $Q_1^{\prime}=Q^{\prime}$. We have the following obvious property:

\medskip\noindent
{\bf (i)} {\it We can assume that $\Spec(Q_1^\prime)\rightarrow \Spec(Q_1)$ is an \'etale morphism. If $\Phi_{Q_1^\prime}$ is a Frobenius lift of the $p$-adic completion $Q_1^{\prime,\wedge}$ of $Q_1^\prime$ compatible with $\sigma$, then the Frobenius of $M_{Q_1}\otimes_{Q_1} Q_1^{\prime,\wedge}$ can be identified via an isomorphism 
$$\iota_{y_1}:(M_{Q_1}\otimes_{Q_1} Q_1^{\prime,\wedge},F^1_{Q_1}\otimes_{Q_1} Q_1^{\prime,\wedge},(t_{\alpha}^{Q_1})_{\alpha\in\mathcal J},\psi_{M_{Q_1}})\rightarrow (M\otimes_{W(k)} Q_1^{\prime,\wedge},F^1\otimes_{W(k)} Q_1^{\prime,\wedge},(t_{\alpha})_{\alpha\in\mathcal J},\epsilon_1\psi_M)\leqno (5)$$ 
with the Frobenius endomorphism $g_{Q_1^\prime}(\phi\otimes\Phi_{Q_1^\prime})$ of $M\otimes_{W(k)} Q_1^{\prime,\wedge}$ for a suitable element $g_{Q_1^\prime}\in\pmb{\text{GL}}_M(Q_1^{\prime,\wedge})$ which modulo the $p$-adic completion of the ideal of $Q_1^\prime$ that defines $z_{H^{(p)}}^\prime$ is the identity. Here $\epsilon_1$ is a unit of $Q_1^{\prime,\wedge}$. From the property (v) of Subsubsection B4.1 applied in the context of the quintuple $(M_R,F^1_R,\Phi,\nabla,(t_{\alpha})_{\alpha\in\mathcal J})$ of Subsection 3.2, we get that the Frobenius of $M_{Q_1}\otimes_{Q_1} Q_1^{\prime,\wedge}$ fixes each $t_{\alpha}^{Q_1}$ with $\alpha\in\mathcal J$. Thus $g_{Q_1^\prime}(\phi\otimes\Phi_{Q_1^\prime})$ fixes each $t_{\alpha}$ with $\alpha\in\mathcal J$ and therefore we have $g_{Q_1^\prime}\in\mathcal G(Q_1^{\prime,\wedge})$.} 

As Theorem \ref{T4} (a) holds for $z\in\mathcal N^{\text{m}}(W(k))$ and as $G_{\mathbb Z_{(p)}}$ is smooth over $\mathbb Z_{(p)}$, we get that $\mathcal G$ is a smooth, closed subgroup scheme of $\pmb{\text{GSp}}(M,\psi_M)$. Thus $\mathcal G^0$ is also smooth over $W(k)$.

As $\mathcal G^0_{B(k)}$ is the subgroup of $\pmb{\text{Sp}}(M[\frac{1}{p}],\psi_M)$ that fixes $t_{\alpha}$ for each $\alpha\in\mathcal J$ and as $\phi$ fixes each $t_{\alpha}$, we get that under the natural action of $\phi$ on $\text{End}(M[\frac{1}{p}])$, we have $\phi(\text{Lie}(\mathcal{G}^0_{B(k)}))=\mathcal{G}^0_{B(k)}$. From this and the existence of the direct sum decomposition $M=F^1\oplus F^0$ defined by $\mu$ (see Subsection 3.1) we get that the two axioms of [59, Subsect. 4.1] hold for the triple $(M,\phi,\mathcal G^0)$. 

Let $\mathcal F$ be be the normalizer of $F^1$ in $\mathcal G$ and let $\mathcal F^0=\mathcal G^0\cap \mathcal F$. From [10, Lem. 2.1.5 and Prop. 2.1.8 (3)] we get that $\mathcal F$ is smooth over $W(k)$ and the product morphism $U\times_{\Spec(W(k))} \mathcal F\rightarrow\mathcal G$ is an open embedding. As $\mathcal F$ is the semiproduct of $\mathcal F^0$ and the image of $\mu$, $\mathcal F^0$ is smooth over $W(k)$ and the product morphism $U\times_{\Spec(W(k))} \mathcal F^0\rightarrow\mathcal G^0$ is an open embedding. 

To prove the Corollary \ref{C1} it suffices to show that $\mathcal C^{\text{m}}$ equipped with the morphism $\mathcal C^{\text{m}}\rightarrow \mathcal A_{r,1,N,k}$ is a quasi Shimura $p$-variety of Hodge type relative to $\mathfrak{F}$ in the sense of [59, Def. 4.2.1], i.e., the axioms (i) of (iii) of loc. cit. hold for $\mathcal C^{\text{m}}$ (more precisely, for the morphism $\mathcal C^{\text{m}}\rightarrow \mathcal A_{r,1,N,k}$).

Axiom (i) of [59, Def. 4.2.1] holds for $\mathcal C^{\text{m}}$ as $\mathcal C^{\text{m}}$ is smooth over $k$ of dimension $d$ (the role of $e_-$ of loc. cit. is played here by the rank $d$ of the $W(k)$-module $\text{Lie}(\mathcal G_{B(k)})\cap\text{Hom}(F^1,F^0)=\text{Lie}(\mathcal G^0_{B(k)})\cap\text{Hom}(F^1,F^0)=\text{Lie}(U)$). Axiom (ii) of [59, Def. 4.2.1] holds as it is just the modulo $p$ variant of Theorem \ref{T1} (b) for $\mathcal C^{\text{m}}$. 

Let $\rho_{y_1}:\Spec(Q_1^{\prime}/pQ_1^{\prime})\rightarrow\mathcal C^{\text{m}}$ be the \'etale morphism induced naturally by the \'etale morphism $\Spec(Q_1^{\prime})\rightarrow \Spec(Q_1)$. Let $W(Q_1^{\prime}/pQ_1^{\prime})$ be the ring of $p$-typical Witt vectors with coefficients in $Q_1^{\prime}/pQ_1^{\prime}$. By shrinking $Q_1^\prime$ we can assume that we have a homomorphism $Q_1^{\prime,\wedge}\rightarrow W(Q_1^{\prime}/pQ_1^{\prime})$ which lifts the identity $Q_1^{\prime,\wedge}/pQ_1^{\prime,\wedge}=Q_1^{\prime}/pQ_1^{\prime}$ and such that the following property holds (cf. the definition of $g_{\text{univ}}$ in Subsection 3.2, the fact that the product morphism $U\times_{W(k)}\mathcal F \rightarrow \mathcal G$ is an open embedding, and the property (ii) of Subsection 3.3):

\medskip\noindent
{\bf (ii)} {\it The composite morphism $\Spec(Q_1^{\prime,\wedge})\rightarrow\mathcal G^0/\mathcal F^0=\mathcal G/\mathcal F$ induced by $g_{Q_1^\prime}$ is formally \'etale.}

\medskip
We consider a finite number of points $y_1,\ldots,y_t\in\mathcal C^{\text{m}}(k)$ such that we have $\cup_{i=1}^t \text{Im}(\rho_{y_i})=\mathcal C^{\text{m}}$; here $\rho_{y_i}:\Spec(Q_1^{\prime}/pQ_1^{\prime)}\rightarrow\mathcal C^{\text{m}}$ for $i\in\{2,\ldots,t\}$ is constructed similarly to $\rho_{y_1}$.

To end the proof of Corollary \ref{C1} it suffices to show that the axioms (iii.a) to (iii.d) of [59, Def. 4.2.1] hold in the context of the family of morphisms $(\rho_{y_i})_{i\in\{1,\ldots,t\}}$. Axiom (iii.a) of [59, Def. 4.2.1] holds for $\mathcal C^{\text{m}}$ as it just says that $\cup_{i=1}^t \text{Im}(\rho_{y_i})=\mathcal C^{\text{m}}$ and that the domain of each $\rho_{y_i}$ is affine. 

As $\mathcal G$ is the semidirect product of $\mathcal G^0$ and the image of $\mu$, by replacing the isomorphism of (5) with an automorphism of $(M\otimes_{W(k)} Q_1^{\prime,\wedge},F^1\otimes_{W(k)} Q_1^{\prime,\wedge},(t_{\alpha})_{\alpha\in\mathcal J},\epsilon_1\psi_M)$ defined by a $Q_1^{\prime,\wedge}$-valued point of the image of $\mu$, we can assume that $\epsilon_1=1$. This implies that in fact we  have $g_{Q_1^\prime}\in\mathcal G^0(Q_1^{\prime,\wedge})$. Now the fact that the axiom (iii.b) of [59, Def. 4.2.1] holds for $\mathcal C^{\text{m}}$ follows from the formally \'etale part of the property (ii) via a natural extension through the homomorphism $Q_1^{\prime,\wedge}\rightarrow W(Q_1^{\prime}/pQ_1^{\prime})$.

Let $i,j\in\{1,\ldots,t\}$. Let $\nu_i:\Spec(Q_i^{\prime})\rightarrow\mathcal N^{\text{m}}_{W(k)}/H^{(p)}$ be an \'etale morphism analogues to the morphism $\nu_1:\Spec(Q_1^{\prime})\rightarrow\mathcal N^{\text{m}}_{W(k)}/H^{(p)}$ induced naturally by the composite morphism $\Spec(Q_1^{\prime})\rightarrow \Spec(Q_1)\rightarrow \mathcal Q^{\text{s}}_{W(k)}\rightarrow\mathcal N^{\text{s}}_{W(k)}/H^{(p)}$ and let $\Spec(Q_{i,j}^\prime)$ define the cartesian product of $\nu_i$ and $\nu_j$. As we have $\epsilon_i=\epsilon_j=1$, the extensions of the isomorphisms $\iota_{y_i}$ and $\iota_{y_j}$ to $W(Q_{ij}^{\prime}/pQ_{ij}^{\prime})$ via composite homomorphisms $Q_i^{\prime,\wedge}\rightarrow Q_{i,j}^{\prime}\rightarrow W(Q_{ij}^{\prime}/pQ_{ij}^{\prime})$ and $Q_j^{\prime,\wedge}\rightarrow Q_{i,j}^{\prime}\rightarrow W(Q_{ij}^{\prime}/pQ_{ij}^{\prime})$ that lift the homomorphisms $Q_i^{\prime}/pQ_i^{\prime}\rightarrow Q_{ij}^{\prime}/pQ_{ij}^{\prime}$ and $Q_j^{\prime}/pQ_j^{\prime}\rightarrow Q_{ij}^{\prime}/pQ_{ij}^{\prime}$ (respectively), when viewed without filtrations differ by an automorphism of $(M\otimes_{W(k)} W(Q_{ij}^{\prime}/pQ_{ij}^{\prime}),(t_{\alpha})_{\alpha\in\mathcal J},\psi_M)$ and thus by an element $h_{ij}\in\mathcal G^0(W(Q_{ij}^{\prime}/pQ_{ij}^{\prime}))$. Thus the axiom (iii.c) of [59, Def. 4.2.1] holds for $\mathcal C^{\text{m}}$.

Axiom (iii.d) of [59, Def. 4.2.1] holds for $\mathcal C^{\text{m}}$ as it just says that regardless of what choices are made in Subsubsection B4.2 to define the closed embedding $\Spec(S)\rightarrow\Spec(R)$ used in Subsection 3.3, the morphism $O_y/pO_y\rightarrow S/pS$ of Subsection 3.3 induced by the isomorphism $s_Z:O_y\rightarrow S$ does exist and it is an isomorphism and moreover the variant of Theorem \ref{T10} described in Subsubsection B4.2 holds.\endproof

\subsubsection{Variant of Corollary \ref{C1}}
The isomorphism class of the quintuple $(M,(t_{\alpha})_{\alpha\in\mathcal J},\psi_M,\mathcal G^0,\mathcal G)$ does not depend on $z\in\mathcal N^{\text{m}}(W(k))$ (as Theorem \ref{T4} (a) holds). But the isomorphism class of the family $\mathfrak{F}$ does depend in general on the connected component $\mathcal C^{\text{m}}_{\infty}$ of $\mathcal N^{\text{m}}_k$ which contains the $k$-valued point $y\in\mathcal N^{\text{m}}(k)$ defined by $z$. This is so as in general the $\mathcal G^0(W(k))$-conjugacy class $[\mu]$ of the cocharacter $\mu:\mathbb G_{m,W(k)}\rightarrow\mathcal G$ does depend on the connected component $\mathcal C^{\text{m}}_{\infty}$, despite the fact that it is easy to check that the $\mathcal G^0(B(k))$-conjugacy class of $\mu_{B(k)}:\mathbb G_{m,B(k)}\rightarrow\mathcal G_{B(k)}$ does not depend on $z\in\mathcal N^{\text{m}}(W(k))$.

If the $\mathcal G^0(W(k))$-conjugacy class $[\mu]$ does not depend on $\mathcal C^{\text{m}}_{\infty}$ and if the hypotheses of Corollary \ref{C1} hold, then (the isomorphism class of) the family $\mathfrak{F}$ does not depend on $z\in\mathcal N^{\text{m}}(W(k))$ and from this and Corollary \ref{C1} we get that $\mathcal N^{\text{m}}_k/H^{(p)}$ itself is a quasi Shimura $p$-variety of Hodge type relative to $\mathfrak{F}$. 

In this paragraph we check that if the hypotheses of Corollary \ref{C1} hold and if moreover $G_{\mathbb Z_{(p)}}$ is a quasi-reductive group scheme for $(G,\mathcal X,v)$, then $\mathcal N^{\text{m}}_k/H^{(p)}$ itself is a quasi Shimura $p$-variety of Hodge type relative to $\mathfrak{F}$. To check this we can assume that $k$ has a countable transcendental degree and thus that we have an $E(G,\mathcal X)$-monomorphism $B(k)\hookrightarrow\mathbb C$. Let $\mathcal G^{\text{r}}$ be the schematic closure in $\mathcal G$ of the normal, reductive subgroup $\mathcal G_{B(k)}^{\text{r}}$ of $\mathcal G_{B(k)}$ which corresponds to $G_{\mathbb Q_p}^{\text{r}}$ via Fontaine comparison theory, cf. Lemma \ref{L6} applied with $\mathcal G_{B(k)}^{\text{v}}=\mathcal G_{B(k)}^{\text{r}}$ and Definition \ref{D2} (a). As Theorem \ref{T4} (a) holds for $z\in\mathcal N^{\text{m}}(W(k))$, $\mathcal G^{\text{r}}$ is isomorphic to $G^{\text{r}}_{\mathbb Z_p}\times_{\Spec(\mathbb Z_p)}\Spec(W(k))$ and thus it is a reductive group scheme. From Lemma \ref{L6} we get that $\mu$ factors through the normal, reductive subgroup scheme of $\mathcal G^{\text{r}}$ of $\mathcal G$. Thus $\mathcal G^{\text{r}}$ is the semiproduct of $\mathcal G^{\text{r},0}=\mathcal G^{\text{r}}\cap \pmb{\text{Sp}}(M,\psi_M)$ and of the image of $\mu$ and therefore $\mathcal G^{\text{r},0}$ is a normal, reductive subgroup scheme of $\mathcal G^{\text{r}}$ (or of $\mathcal G$ or $\mathcal G^0$). The $\mathcal G(\mathbb C)$-conjugacy class of $\mu_{\mathbb C}$ does not depend on $z\in\mathcal N^{\text{m}}(W(k))$ as via isomorphisms $(M\otimes_{W(k)} \mathbb C,(t_{\alpha})_{\alpha\in\mathcal J})\rightarrow (W^{\vee}\otimes_{\mathbb Q} \mathbb C,(v_{\alpha})_{\alpha\in\mathcal J})$ it corresponds to the $G(\mathbb C)$-conjugacy class of the cocharacters $\mu_h:\mathbb G_{m,\mathbb C}\rightarrow G_{\mathbb C}$ with $h\in\mathcal X$, cf. proof of Lemma \ref{L6}. It is well-known that the last two sentences imply that $[\mu]$ equals to the $\mathcal G^{\text{r},0}(W(k))$-conjugacy class of $\mu$ and it does not depend on $z\in\mathcal N^{\text{m}}(W(k))$. Thus $\mathcal N^{\text{m}}_k/H^{(p)}$ itself is a quasi Shimura $p$-variety of Hodge type relative to $\mathfrak{F}$, cf. previous paragraph. 

\medskip
The following lemma will be used in the proof of Lemma \ref{L7} and in Subsections 5.4 and 5.6.

\begin{lemma}\label{L7} 
Let $\mu:\mathbb G_{m,W(k)}\rightarrow \mathcal G$ and $M=F^1\oplus F^0$ be as in Subsection 3.1. Let $y\in\mathcal N^{\text{s}}(k)$ be defined by $z\in\mathcal N^{\text{s}}(W(k))=\mathcal N(W(k))$. Let $\mu_1:\mathbb G_{m,W(k)}\rightarrow \mathcal G$ be a cocharacter such that we have a direct sum decomposition $M=F^1_1\oplus F^0_1$ with the properties that $\mathbb G_{m,W(k)}$ acts through $\mu_1$ on each $F^i_1$ via the weight $-i$ and we have $F^1_1/pF^1_1=F^1/pF^1$. Then the following three properties hold:

\medskip
{\bf (a)} There exists a point $z_1\in\mathcal N^{\text{s}}(W(k))$ that lifts $y\in\mathcal N^{\text{s}}(k)$ and such that the principally quasi-polarized filtered $F$-crystal over $k$ of $z_1^*(\mathcal A,\lambda_{\mathcal A})$ is $(M,F^1_1,\phi,\psi_M)$.

{\bf (b)} If $p>2$, then a point $z_1$ as in part (a) is unique. 

\smallskip
{\bf (c)} We assume that $p=2$ and that $y$ factors through $\mathcal N^{\text{m}}$. We also assume that one of the following two conditions holds:

\medskip
{\bf (c.i)} the abelian variety $A_k$ is ordinary and $\mathcal G$ is smooth;

\smallskip
{\bf (c.ii)} $G_{\mathbb Z_{(p)}}$ is a quasi-reductive group scheme for $(G,\mathcal X,v)$.

\medskip\noindent
Then there exist exactly $2^a$ such $z_1$'s, where $a$ is the multiplicity of the Newton polygon slope $-1$ for $(\text{Lie}(\mathcal G),\phi)$. Moreover, if (c.i) holds, then we can choose $\mu_1$ and $z_1$ such that the abelian scheme $z_1^*(\mathcal A)$ is the canonical lift of $A_k$.
\end{lemma}

\medskip\proof
For $n\in\mathbb N^*$ let $W_n(k):=W(k)/p^nW(k)$. By induction on $n\in\mathbb N^*$ we show that there exists a point $z(n)\in\mathcal N^{\text{s}}(W(k))=\mathcal N(W(k))$ that has the following three properties:

\medskip\noindent
{\bf (i)} it lifts $y\in\mathcal N^{\text{s}}(k)$;

\smallskip\noindent
{\bf (ii)} for $n\ge 2$ the $W_{n-1}(k)$-valued points of $\mathcal N^{\text{s}}$ defined by $z(n-1)$ and $z(n)$ coincide;

\smallskip\noindent
{\bf (iii)} the principally quasi-polarized filtered $F$-crystal over $k$ of $z(n)^*(\mathcal A,\lambda_{\mathcal A})$ is a quadruple of the form $(M,F^1_1(n),\phi,\psi_M)$, where $F^1_1(n)$ is congruent  to $F^1_1$ modulo $p^n$. 

\medskip
Let $z(1):=z$; obviously the base of the induction for $n=1$ holds. The passage from $n$ to $n+1$ goes as follows. Not to introduce extra notation by replacing $z$ with $z(n)$, we can assume that $z(n)=z$; thus we have $F^1/p^nF^1=F^1_1/p^nF^1_1$. Let $v\in p\text{Lie}(U)$ be such that for $u:=1_M+v$ we have $u(F^1)=F^1_1$, cf. Lemma \ref{L12}. As $F^1/p^nF^1=F^1_1/p^nF^1_1$, we have $v\in p^n\text{Lie}(U)$.

As the image of the Kodaira--Spencer map of $\nabla$ is $\text{Lie}(U)\otimes_{W(k)} R$ (cf. property (iii) of Subsubsection B4.1) and as the morphism $\tau_R:\Spec(R)\rightarrow\mathcal N^{\text{s}}$ is formally smooth, from the relation $v\in p^n\text{Lie}(U)$ we get that there exists a lift $z(n+1)\in\mathcal N^{\text{s}}(W(k))$ of $z(n)$ modulo $p^n$ such that the principally quasi-polarized filtered $F$-crystal over $k$ of $z(n+1)^*(\mathcal A,\lambda_{\mathcal A})$ is $(M,F^1_1(n+1),\phi,\psi_M)$, where $F^1_1(n+1)$ is congruent  to $u(F^1)=F^1_1$ modulo $p^{n+1}$ (this holds even if $p=2$; for instance, the proof of [63, Prop. 6.4.6 (b)] applies entirely in the slightly more general context of our  present lemma). This ends the induction.

From the property (ii) we get that there exists a point $z_1\in \mathcal N^{\text{s}}(W(k))$ that lifts $z(n)$ modulo $p^n$ for all $n\in\mathbb N^*$. Thus $z_1$ also lifts $y$, cf. property (i). The principally quasi-polarized filtered $F$-crystal over $k$ of $z_1^*(\mathcal A,\lambda_{\mathcal A})$ is $(M,F^1_1,\phi,\psi_M)$, cf. property (iii). Thus part (a) holds. Part (b) follows from Theorem \ref{T1} (b) and the Grothendieck--Messing deformation theory.

If the condition (c.i) holds, then part (c) follows from Theorem \ref{T1} (b) and Theorem \ref{T11} (c) (more precisely, the constructions of Subsection 3.2 needed to prove Theorem \ref{T1} (b) in Subsection 3.3 are exactly the constructions of Subsubsection B4.1 and Theorem \ref{T10}). 

We are left to check that part (c) holds if the  condition (c.ii) holds. Let $G^{\text{r}}_{\mathbb Z_p}=G^{\text{r}}_{\mathbb Z_2}$ be as in Definition \ref{D2} (a). As $y$ factors through $\mathcal N^{\text{m}}$, the normal, closed subgroup scheme $\mathcal G^{\text{r}}$ of $\mathcal G$ obtained as in Lemma \ref{L6} but for $G^{\text{v}}_{\mathbb Q_p}=G^{\text{r}}_{\mathbb Q_2}$, is a reductive group scheme over $W(k)$ through which $\mu$ and thus through which $\mu_1$'s factor as well (note that $U$ is a closed subgroup scheme of $\mathcal G^{\text{r}}$). Thus this case of part (c) also follows from Theorem \ref{T1} (b) and Theorem \ref{T11} (c) and (d) applied to $(M,\phi,\mathcal G^{\text{r}})$ instead of $(M,\phi,\mathcal G)$, with the cocharacter $\mu_1$ chosen such that we have $\phi(F^1_1)=2F^1_1$ (see the second paragraph of the proof of Theorem \ref{T11} for the existence of such a $\mu_1$).\endproof

\subsection{Proof of Proposition \ref{P1}} Let $y\in\mathcal N(k)$ be such that $A_k:=y^*(\mathcal A)$ is an ordinary abelian variety. From [44, Cor. 3.8] we get that $y$ factors through $\mathcal N^{\text{s}}$. Thus to prove the Proposition \ref{P1} we can assume that $p=2$ and we have to show that $y\in\mathcal N^{\text{m}}(k)$. We will use the previous notation for a lift $z\in\mathcal N^{\text{s}}(W(k))=\mathcal N(W(k))$ of $y$. We have a direct sum decomposition $M=F^1_1\oplus F^0_1$ such that $\phi(F^1_0)=F^1_0$ and $\phi(F^1_1)=2F^1_1$; obviously, $F^1_1/2F^1_1=F^1/2F^1$. The cocharacter $\mu_1:\mathbb G_{m,W(k)}\rightarrow\pmb{\text{GL}}_M$ associated to it factors through $\mathcal G$, cf. second paragraph of the proof of Theorem \ref{T11}. Based on Lemma \ref{L7} (a) we can assume that $F^1=F^1_1$ and $\mu=\mu_1$; thus $F^1$ is the Hodge filtration of $M$ which defines the canonical lift $A_{\text{can}}$ of $A_k$. The $\text{Gal}(B(k))$-module $H^1(D_{\text{can}})$ is canonically identified with a $\text{Gal}(B(k))$-submodule of ${1\over 2}H^1(D)$ which contains $2H^1(D)$ (see [64, Subsubsects. 2.2.1 and 2.2.3] applied as in the  second paragraph of the proof of [64, Lem. 2.2.5]). Let $H^1(D_{\text{can}})=H^1(D_{\text{can}})_1\oplus H^1(D_{\text{can}})_0$ be the direct sum decomposition that corresponds naturally to the direct sum decomposition $(M,F^1,\phi)=(F^1,F^1,\phi)\oplus (F^0,0,\phi)$. As we have a short exact sequence $0\rightarrow H^1(D_{\text{can}})_0\rightarrow H^1(D)\rightarrow H^1(D_{\text{can}})_1\rightarrow 0$, there exists $c\in {1\over 2}\text{Hom}(H^1(D_{\text{can}})_1,H^1(D_{\text{can}})_0)$ such that we have $(1_{H^1(D_{\text{can}})}+c)(H^1(D_{\text{can}}))=H^1(D)$. Let $\mu^{\acute{et}}:\mathbb G_{m,\mathbb Z_2}\rightarrow \pmb{\text{GL}}_{H^1(D_{\text{can}})}$ be the cocharacter that fixes $H^1(D_{\text{can}})_0$ and that acts on $H^1(D_{\text{can}})_1$ via the weight $-1$. We consider an isomorphism $(H^1(D),(u_{\alpha})_{\alpha\in\mathcal J})\rightarrow (L^{\vee}_{(2)}\otimes_{\mathbb Z_{(2)}} \mathbb Z_2,(v_{\alpha})_{\alpha\in\mathcal J})$ (cf. Lemma \ref{L3} (a)) to be viewed as an identification.

Let $G^{\text{r}}_{\mathbb Z_2}=G^{\text{r}}_{\mathbb Z_p}$ be as in Definition \ref{D2} (a); it is a reductive, normal, closed subgroup scheme of $G_{\mathbb Z_2}$ and thus of $\pmb{\text{GL}}_{H^1(D)}$, cf. last identification. We know that $\mu^{\acute{et}}_{\mathbb Q_2}$ is the \'etale counterpart of the cocharacter $\mu_{B(k)}$ of $\mathcal G_{B(k)}$, i.e., they correspond to each other via Fontaine comparison theory (the functorial isomorphism $i_D$ of the property (v) of Subsection B2 preserves the direct sum decompositions of the previous paragraph). Thus from Lemma \ref{L6} we get that $\mu^{\acute{et}}_{\mathbb Q_2}$ factors through $G^{\text{r}}_{\mathbb Q_2}$. Let $U^{\acute{et}}_{\text{bigg}}$ (resp. $U^{\acute{et}}$) be the unipotent radical of the parabolic subgroup scheme of $\pmb{\text{GL}}_{H^1(D)}$ (resp. of the parabolic subgroup scheme $P^{\text{r}}_{\mathbb Z_2}$ of $G^{\text{r}}_{\mathbb Z_2}$) that normalizes $H^1(D)_0$ (cf. the existence of $\mu^{\acute{et}}_{\mathbb Q_2}$ and the fact that the $\mathbb Z_2$-schemes of parabolic subgroup schemes of reductive group schemes over $\mathbb Z_2$ are projective, see [15, Vol. III, Exp. XXVI, Cor. 3.5]). As a  $P^{\text{r}}_{\mathbb Z_2}(\mathbb Q_2)$-conjugate of $\mu^{\acute{et}}_{\mathbb Q_2}$ extends to a cocharacter of $G^{\text{r}}_{\mathbb Z_2}$, it is easy to see that we have $G^{\text{r}}_{\mathbb Z_2}\cap U^{\acute{et}}_{\text{bigg}}=U^{\acute{et}}$ (this is a particular case of [10, Lem. 2.1.5 and Prop. 2.1.8 (3)]).

We claim that there exists an element $g\in G^{\text{r}}_{\mathbb Z_2}(\mathbb Q_2)$ such that we have $g(H^1(D_{\text{can}}))=H^1(D)=L^{\vee}_{(2)}\otimes_{\mathbb Z_{(2)}} \mathbb Z_2$. It suffices to show that the reduction $\bar e$ of $e:=2c\in\text{Lie}(U^{\acute{et}}_{\text{bigg}})$ modulo $2$  belongs to $\text{Lie}(U^{\acute{et}}_{\mathbb F_2})$. Let $G^{\text{r}}_{1,\mathbb Z_2}:=(1_{H^1(D_{\text{can}})}-c)G^{\text{r}}_{\mathbb Z_2}(1_{H^1(D_{\text{can}})}+c)$; it is a reductive, closed subgroup scheme $\pmb{\text{GL}}_{H^1(D_{\text{can}})}$. For $t\in W(k)$, $\mu^{\acute{et}}((1+2t)^{-1})$ normalizes $H^1(D)\otimes_{\mathbb Z_2} W(k)=(1_{H^1(D_{\text{can}})}+c)(H^1(D_{\text{can}}))\otimes_{\mathbb Z_2} W(k)$ and thus its conjugate under $1_{H^1(D_{\text{can}})}-c$ belongs to $G^{\text{r}}_{1,\mathbb Z_2}(W(k))$. Therefore $1_{H^1(D_{\text{can}})/2H^1(D_{\text{can}})}-\bar t\bar e$ belongs to $G^{\text{r}}_{1,\mathbb Z_2}(k)$ for all $\bar t\in k$ and thus $\bar e\in\text{Lie}(G^{\text{r}}_{1,\mathbb F_2})$. Conjugating via $1_{H^1(D_{\text{can}})}+c$ we get that $\bar e\in\text{Lie}(G^{\text{r}}_{\mathbb F_2})\cap \text{Lie}(U^{\acute{et}}_{\text{bigg},\mathbb F_2})=\text{Lie}(U^{\acute{et}}_{\mathbb F_2})$. Thus the claim holds. 

Let $z_{\text{can}}:\Spec(W(k))\rightarrow\mathcal M_{r,O_{(v)}}$ be the morphism which is the canonical lift of the composite morphism $\Spec(k)\rightarrow\mathcal N\rightarrow \mathcal M_{r,O_{(v)}}$ which is defined naturally by $y$ and which factors through the ordinary locus of $\mathcal M_{r,k(v)}$. From the above claim we get that the generic fibres of $z_{\text{can}}$ and $z$ define $B(k)$-valued points of $\mathcal M_{r,E(G,\mathcal X)}$ which are images of complex points of $\text{Sh}(G,\mathcal X)$ that differ by the right translation through the element $g\in G(\mathbb Q_2)\leqslant G(\mathbb A_f)$, cf. proof of Lemma \ref{L3}. Therefore we have a unique factorization $z_{\text{can}}:\Spec(W(k))\rightarrow\mathcal N$ with the property that each $t_{\alpha}$ is the crystalline realization of the Hodge cycle  $z_{\text{can},B(\mathbb F)}^*(\mathcal A_{E(G,\mathcal X)})$ on $A_{\text{can},B(\mathbb F)}$. We know that $z_{\text{can}}:\Spec(W(k))\rightarrow\mathcal N$ factors through $\mathcal N^{\text{s}}$ (cf. Theorem \ref{T1} (a)) and even through $\mathcal N^{\text{m}}$ (cf. Theorem \ref{T4} (a)). From this and the existence of $g\in G^{\text{r}}_{\mathbb Z_2}(\mathbb Q_2)$ such that $g(H^1(D_{\text{can}}))=H^1(D)$ we get that Theorem \ref{T4} (a) holds for $z\in\mathcal N^{\text{s}}(W(k))$ and therefore we have $z\in \mathcal N^{\text{m}}(W(k))$ and $y\in\mathcal N^{\text{m}}(k)$. Thus Proposition \ref{P1} holds.\endproof 

\section{Applications to integral models}
In this section we take $k$ to be an algebraic closure of $k(v)$. This implies that there exist $O_{(v)}$-monomorphisms $W(k)\hookrightarrow \mathbb C$.  

Lemma \ref{L8} presents a simple criterion on when the $k(v)$-scheme $\mathcal N^{\text{m}}_{k(v)}$ is non-empty or when the $W(k)$-valued points of $\mathcal N^{\text{m}}_{W(k)}$ are Zariski dense. In Subsection 4.1 we apply Theorem \ref{T1} (a) and Lemma \ref{L8} (a) to prove the existence of good smooth integral models of $\text{Sh}_{\tilde H}(G,\mathcal X)$ over $O_{(v)}$ for a large class of maximal compact, open subgroups $\tilde H$ of $G_{\mathbb Q_p}(\mathbb Q_p)$. Corollary \ref{C3} can be viewed as a smooth solution (answer) to the conjecture (question) of Langlands of [31, p. 411] for Shimura varieties of Hodge type. Theorem \ref{T7} shows that, in the case when $G_{\mathbb Q_p}$ splits over an unramified extension of $\mathbb Q_p$, Subsection 4.1 extends naturally to the case of parahoric subgroups of $G_{\mathbb Q_p}(\mathbb Q_p)$.

\begin{lemma}\label{L8}
We assume that one of the following two conditions holds:

\medskip\noindent
{\bf (i)} there exists a smooth, affine group scheme $G^{\text{v}}_{\mathbb Z_{(p)}}$ over $\mathbb Z_{(p)}$ that extends $G$ (i.e., it has $G$ as its generic fibre), that has a special fibre $G^{\text{v}}_{\mathbb F_p}$ of the same rank as $G$, and that has the property that there exists a homomorphism $G^{\text{v}}_{\mathbb Z_{(p)}}\rightarrow G_{\mathbb Z_{(p)}}$ which extends the identity automorphism of $G$;

\smallskip\noindent
{\bf (ii)} we have $e(v)=1$ and $G_{\mathbb Z_{(p)}}$ is a quasi-reductive group scheme for $(G,\mathcal X,v)$ in the sense of Definition \ref{D2} (a). 

\medskip
{\bf (a)} Then $e(v)=1$ and the $k(v)$-scheme $\mathcal N^{\text{m}}_{k(v)}$ (and thus also $\mathcal N^{\text{s}}_{k(v)}$) is non-empty.

\smallskip
{\bf (b)} Then the $W(k)$-valued points of $\mathcal N^{\text{m}}_{W(k)}$ are Zariski dense in $\mathcal N^{\text{m}}_{W(k)}$.
\end{lemma} 

\medskip\proof
We prove part (a). We first assume that (i) holds. Each torus of $G^{\text{v}}_{\mathbb F_p}$ lifts to a torus of $G^{\text{v}}_{\mathbb Z_p}$, cf. [15, Vol. II, Exp. XII, Cor. 1.10]. Thus $G_{\mathbb Z_p}^{\text{v}}$ has tori of rank equal to the rank of $G$. Let $T^{\text{v}}_{\mathbb Z_{(p)}}$ be a torus of $G^{\text{v}}_{\mathbb Z_{(p)}}$ of the same rank as $G$ and such that there exists $h^{\text{v}}\in \mathcal X$ which factors through $T^{\text{v}}_{\mathbb R}$. Its existence is implied by [22, Lem. 5.5.3]. The pair $(T^{\text{v}}_{\mathbb Q},\{h^{\text{v}}\})$ is a Shimura subpair of $(G,\mathcal X)$ and therefore we have an inclusion $E(G,\mathcal X)\subset E(T^{\text{v}}_{\mathbb Q},\{h^{\text{v}}\})$ of reflex fields. Each prime of $E(T^{\text{v}}_{\mathbb Q},\{h^{\text{v}}\})$ that divides $v$ is unramified over $p$ (cf. [38, Prop. 4.6 and Cor. 4.7]) and thus we have $e(v)=1$. The intersection $H^{\text{v}}:=H\cap T^{\text{v}}_{\mathbb Z_{(p)}}(\mathbb Q_p)$ is the unique hyperspecial subgroup $T^{\text{v}}_{\mathbb Z_{(p)}}(\mathbb Z_p)$ of $T^{\text{v}}_{\mathbb Z_{(p)}}(\mathbb Q_p)$. Therefore there exists an integral model $\mathcal Z^{\text{v}}$ of $\text{Sh}_{H^{\text{v}}}(T^{\text{v}}_{\mathbb Q},\{h^{\text{v}}\})$ over the spectrum of the normalization of $O_{(v)}$ in $E(T^{\text{v}}_{\mathbb Q},\{h^{\text{v}}\})$ which is a pro-finite pro-\'etale cover of $\Spec(O_{(v)})$, cf. either [37, Rm. 2.16] or [54, Ex. 3.2.8]. In particular, $\mathcal Z^{\text{v}}$ is a regular, formally \'etale, faithfully flat $O_{(v)}$-scheme. The functorial morphism $\text{Sh}_{H^{\text{v}}}(T^{\text{v}}_{\mathbb Q},\{h^{\text{v}}\})\rightarrow \text{Sh}_H(G,\mathcal X)$ of $E(G,\mathcal X)$-schemes extends uniquely to a morphism $\mathcal Z^{\text{v}}\rightarrow\mathcal N^{\text{s}}$ of $O_{(v)}$-schemes, cf. Theorem \ref{T1} (a). There exist points $z\in\mathcal Z^{\text{v}}(W(k))$. Let $(v_{\alpha})_{\alpha\in\mathcal J^{\text{v}}}$ be a family of tensors of $\mathcal T(W^{\vee})$ such that $T^{\text{v}}_{\mathbb Q}$ is the subgroup of $\pmb{\text{GL}}_{W^{\vee}}$ that fixes $v_{\alpha}$ for all $\alpha\in\mathcal J^{\text{v}}$. We can assume that $\mathcal J\subset\mathcal J^{\text{v}}$ and that for each $\alpha\in\mathcal J$, the tensor $v_{\alpha}$ is the tensor introduced in Subsection 2.3. We will use the notation of Subsection 3.1 for $z\in\mathcal Z^{\text{v}}(W(k))$. From Theorem \ref{T4} (a) applied to the point $z\in\mathcal Z^{\text{v}}(W(k))$ we get that there exists an isomorphism $(M,(t_{\alpha})_{\alpha\in\mathcal J^{\text{v}}})\rightarrow (L^{\vee}_{(p)}\otimes_{\mathbb Z_{(p)}} W(k),(v_{\alpha})_{\alpha\in\mathcal J^{\text{v}}})$ (each $t_{\alpha}$ with $\alpha\in\mathcal J^{\text{v}}$, is the de Rham realization of the Hodge cycle on $A_{B(k)}$ that corresponds naturally to $v_{\alpha}$). Thus as $\mathcal J\subset\mathcal J^{\text{v}}$, Theorem \ref{T4} (a) holds for the $W(k)$-valued point of $\mathcal N^{\text{s}}$ defined by $z$. From this and Proposition \ref{P3} (b) we get that this last point factors through $\mathcal N^{\text{m}}$. Therefore the $k(v)$-scheme $\mathcal N^{\text{m}}_{k(v)}$ is non-empty. 

We now assume that (ii) holds; thus $e(v)=1$. Let $G^{\text{r}}_{\mathbb Z_p}$ and $\mu_v$ be as in Definition \ref{D2} (a). Let $T^{\text{r}}_{\mathbb F_p}$ be a maximal torus of $G^{\text{r}}_{\mathbb F_p}$. Due to the existence of $\mu_v$, $T^{\text{r}}_{\mathbb F_p}$ has positive rank. The torus $T^{\text{r}}_{\mathbb F_p}$ lifts to a torus $T^{\text{r}}_{\mathbb Z_p}$ of $G^{\text{r}}_{\mathbb Z_p}$, cf. [15, Vol. II, Exp. XII, Cor. 1.10]. Let $T^{\text{v}}_{0,\mathbb Q_p}$ be a maximal torus of $G_{\mathbb Q_p}$ which has $T^{\text{r}}_{\mathbb Q_p}$ as a subtorus. Let $T^{\text{v}}$ be a maximal torus of $G$ such that there exists an element $h^{\text{v}}\in \mathcal X$ which factors through $T^{\text{v}}_{\mathbb R}$ and moreover $T^{\text{v}}_{\mathbb Q_p}$ is $H$-conjugate to $T^{\text{v}}_{0,\mathbb Q_p}$. Again, the existence of $T^{\text{v}}$ is implied by [22, Lem. 5.5.3]. Thus (up to $H$-conjugation) we can assume that we have $T^{\text{v}}_{0,\mathbb Q_p}=T^{\text{v}}_{\mathbb Q_p}$. 

The intersection $H^{\text{v}}:=H\cap T^{\text{v}}(\mathbb Q_p)$ is not necessarily the maximal compact, open subgroup of $T^{\text{v}}(\mathbb Q_p)$ and the subgroup $T^{\text{v}}(\mathbb Q)H^{\text{v}}$ of $T^{\text{v}}(\mathbb Q_p)$ is not necessarily $T^{\text{v}}(\mathbb Q_p)$. But the intersection $T^{\text{r}}_{\mathbb Z_p}(\mathbb Q_p)\cap H$ is the unique hyperspecial subgroup $T^{\text{r}}_{\mathbb Z_{p}}(\mathbb Z_p)$ of $T^{\text{r}}_{\mathbb Z_{p}}(\mathbb Q_p)$. We fix an $O_{(v)}$-monomorphism $W(k(v))\hookrightarrow\mathbb C$ as in Definition \ref{D2} (a). As $\mu_{h^{\text{v}}}$ and $\mu_{v,\mathbb C}$ are $G(\mathbb C)$-conjugate and as $G^{\text{r}}_{\mathbb C}$ is a normal subgroup of $G_{\mathbb C}$, $\mu_{h^{\text{v}}}$ factors through the intersection $T^{\text{v}}_{\mathbb C}\cap G^{\text{r}}_{\mathbb C}$ and therefore through $T^{\text{r}}_{\mathbb C}=T^{\text{r}}_{\mathbb Z_p}\times_{\Spec(\mathbb Z_p)} \Spec(\mathbb C)$. Thus as $T^{\text{r}}_{\mathbb Z_p}$ splits over a finite, unramified extension of $\mathbb Z_p$, we get that the field of definition $E(T^{\text{v}}_{\mathbb Q},\{h^{\text{v}}\})$ of $\mu_{h^{\text{v}}}$ is a number subfield of $\mathbb C$ that contains $E(G,\mathcal X)$ and that is unramified over $v$. From the class field theory (see [30, Th. 4 of p. 220]) and the reciprocity map of [37, pp. 163--164] we easily get that each connected component of $\text{Sh}_{H^{\text{v}}}(T^{\text{v}}_{\mathbb Q},\{h^{\text{v}}\})_{\mathbb C}$ is defined over the spectrum of an abelian extension of $E(T^{\text{v}}_{\mathbb Q},\{h^{\text{v}}\})$ unramified over all primes of $E(T^{\text{v}}_{\mathbb Q},\{h^{\text{v}}\})$ that divide $v$. Thus there exists an integral model $\mathcal Z^{\text{v}}$ of $\text{Sh}_{H^{\text{v}}}(T^{\text{v}}_{\mathbb Q},\{h^{\text{v}}\})$ over the normalization of $O_{(v)}$ in $E(T^{\text{v}}_{\mathbb Q},\{h^{\text{v}}\})$ which has the same properties as above. Let $z\in\mathcal Z^{\text{v}}(W(k))$.

Let $(v_{\alpha})_{\alpha\in\mathcal J^{\text{r}}}$ be a family of tensors of $\mathcal T(W^{\vee}\otimes_{\mathbb Q} \mathbb Q_p)$ such that $T^{\text{r}}_{\mathbb Q_p}$ is the subgroup of $\pmb{\text{GL}}_{W^{\vee}\otimes_{\mathbb Q} \mathbb Q_p}$ that fixes $v_{\alpha}$ for all $\alpha\in\mathcal J^{\text{r}}$. We can assume that $\mathcal J\subset\mathcal J^{\text{r}}$ and that for each $\alpha\in\mathcal J$, the tensor $v_{\alpha}$ is the tensor introduced in Subsection 2.3. 

We will use the notation of Subsection 3.1 for $z\in\mathcal Z^{\text{v}}(W(k))$ and for $k$ of countable transcendental degree. Let $\rho_D:\text{Gal}(B(k))\rightarrow \pmb{\text{GL}}_{H^1_{\acute et}(A_{B(k)},\mathbb Q_p)}(\mathbb Q_p)\rightarrow \pmb{\text{GL}}_{L_{(p)}^{\vee}\otimes_{\mathbb Z_{(p)}} \mathbb Q_p}(\mathbb Q_p)$ be the $p$-adic Galois representation associated to the Barsotti--Tate group $D$ of $A$. Let $\mathcal D_{\mathbb Q_p}^{\acute et}$ be the schematic closure of $\text{Im}(\rho_D)$ in $\pmb{\text{GL}}_{L_{(p)}^{\vee}\otimes_{\mathbb Z_{(p)}} \mathbb Q_p}$; it is a connected group (cf. Subsection B1) which is a subgroup of $T^{\text{v}}_{\mathbb Q_p}$. As the groups $T^{\text{v}}_{\mathbb Q_p}$, $T^{\text{r}}_{\mathbb Q_p}$, and $G^{\text{r}}_{\mathbb Q_p}$ are normalized by $\mathcal D_{\mathbb Q_p}^{\acute et}$, we can speak about the subgroups $\mathcal T^{\text{r}}_{B(k)}$, $\mathcal T^{\text{v}}_{B(k)}$, and $\mathcal G^{\text{r}}_{B(k)}$ of $\mathcal G_{B(k)}$ that correspond to $T^{\text{v}}_{\mathbb Q_p}$, $T^{\text{r}}_{\mathbb Q_p}$, and $G^{\text{r}}_{\mathbb Q_p}$ (respectively) via Fontaine comparison theory for $D$ (cf. Lemma \ref{L15} (a)).  The generic fibre of $\mu$ factors through $\mathcal T^{\text{v}}_{B(k)}$ (cf. Subsection 3.1 applied in the context of $z\in\mathcal Z^{\text{v}}(W(k))$) and through $\mathcal G^{\text{r}}_{B(k)}$ (cf. Lemma \ref{L6} applied with $G^{\text{v}}_{\mathbb Q_p}=G^{\text{r}}_{\mathbb Q_p}$ to the image of $z\in\mathcal Z^{\text{v}}(W(k))$ in $\mathcal N^{\text{s}}(W(k))$) and thus it factors through $\mathcal T^{\text{r}}_{B(k)}=\mathcal T^{\text{v}}_{B(k)}\cap\mathcal G^{\text{r}}_{B(k)}$. From this and Lemma \ref{L15} (b) we get that $\mathcal D_{\mathbb Q_p}^{\acute et}$ is a subgroup of $T^{\text{r}}_{\mathbb Q_p}$. This implies that each $v_{\alpha}$ with $\alpha\in\mathcal J^{\text{r}}$ defines naturally an \'etale Tate-cycle $u_{\alpha}$ on $D_{B(k)}$. 

As $T^{\text{r}}_{\mathbb Z_p}$ is a torus, from Theorem \ref{T9} applied to the pair $(D,(u_{\alpha})_{\alpha\in\mathcal J^{\text{r}}})$, from Formula (3), and from Lemma \ref{L3} (b) applied to $\text{Sh}(T^{\text{v}}_{\mathbb Q},\{h^{\text{v}}\})$ we get that there exist isomorphisms 
$$(M,(t_{\alpha})_{\alpha\in\mathcal J^{\text{r}}},\psi_M)\rightarrow (H^1_{\acute et}(A_{B(k)},\mathbb Z_p)\otimes_{\mathbb Z_p} W(k),(u_{\alpha})_{\alpha\in\mathcal J^{\text{r}}},\psi_{H^1_{\text{\'et}}})\rightarrow (L^{\vee}_{(p)}\otimes_{\mathbb Z_{(p)}} W(k),(v_{\alpha})_{\alpha\in\mathcal J^{\text{r}}},\psi^{\vee})$$ 
(each $t_{\alpha}\in\mathcal T(M[\frac{1}{p}])$ with $\alpha\in\mathcal J^{\text{r}}$ corresponds to $u_{\alpha}$ via Fontaine comparison theory for $D$). As $\mathcal J\subset\mathcal J^{\text{r}}$, from this and Proposition \ref{P3} (b) we get that the image of $z\in\mathcal Z^{\text{v}}(W(k))$ in $\mathcal N^{\text{s}}(W(k))$ belongs to $\mathcal N^{\text{m}}(W(k))$. Thus the $k(v)$-scheme $\mathcal N^{\text{m}}_{k(v)}$ is non-empty, i.e.,  part (a) holds.

We prove part (b). If (i) holds, let $T_{\mathbb Q_p}^{\text{r}}:=T_{\mathbb Q_p}^{\text{v}}$. Thus $T_{\mathbb Q_p}^{\text{r}}$ is well-defined regardless of which one of the conditions (i) and (ii) holds. Due to Formula (1) and the fact that $\mathcal N^{\text{m}}$ is $G(\mathbb A_f^{(p)})$-invariant, to prove that the $W(k)$-valued points of $\mathcal N^{\text{m}}_{W(k)}$ are Zariski dense, it suffices to show that for each open subset $\mathcal K$ of $\mathcal X$ and for every element of $G(\mathbb Q)\backslash G(\mathbb Q_p)/H$, we can choose a representative $g_j\in G(\mathbb Q_p)\leqslant G(\mathbb A_f)$ of this element and we can choose $(T_{\mathbb Q}^{\text{v}},\{h^{\text{v}}\})$ such that $h^{\text{v}}\in\mathcal K$ and the elements of $T_{\mathbb Q_p}^{\text{r}}(\mathbb Q_p)\cap H$ act via left translation trivially on the image of $g_j$ in $G(\mathbb Q_p)/H$ (this is so as from the class field theory and the reciprocity map of [37, pp. 163--164] we easily get that the complex point $[h^{\text{v}},g_j]$ of $\text{Sh}_{H^{\text{v}}}(T^{\text{v}}_{\mathbb Q},\{h^{\text{v}}\})$ is defined over the spectrum of an abelian extension of $E(T^{\text{v}}_{\mathbb Q},\{h^{\text{v}}\})$ unramified over all primes of $E(T^{\text{v}}_{\mathbb Q},\{h^{\text{v}}\})$ that divide $v$). 

If (i) holds, then the existence of $T^{\text{v}}_{\mathbb Z_{(p)}}$ implies that $G_{\mathbb Q_p}$ splits over a finite unramified extension of $\mathbb Q_p$ and therefore we have $G(\mathbb Q_p)=G(\mathbb Q)H$ (cf. [38, Lem. 4.10]). This implies that we can take $g_j$ to be the identity element and based on [22, Lem. 5.5.3] we can assume that $h^{\text{v}}\in\mathcal K$. 

If (ii) holds, then $g_j$ can be any representative and we choose $T^{\text{r}}_{\mathbb Z_p}$ so that it is also a maximal torus of $g_jG_{\mathbb Z_p}^{\text{r}}g_j^{-1}$ (this is argued similarly to [62, Subsect. 5.2] based on [52, pp. 43--44, Subsect. 3.4, and Subsubsect. 3.8.1]); based on [22, Lem. 5.5.3] we can assume that $h^{\text{v}}\in\mathcal K$ and that $T^{\text{v}}_{\mathbb Q_p}$ is $H\cap g_jHg_j^{-1}$-conjugate to $T^{\text{v}}_{0,\mathbb Q_p}$ and thus that the elements of $T_{\mathbb Q_p}^{\text{r}}(\mathbb Q_p)\cap H$ act via left translation trivially on the image of $g_j$ in $G(\mathbb Q_p)/H$. We conclude that part (b) holds.\endproof

\subsection{Integral models for maximal compact, open subgroups} Let $\tilde H$ be a maximal compact, open subgroup of $G_{\mathbb Q_p}(\mathbb Q_p)$. Let $\tilde G_{\mathbb Z_{p}}$ be a smooth, affine group scheme over $\mathbb Z_{p}$ that extends $G_{\mathbb Q_p}$ and such that we have $\tilde H=\tilde G_{\mathbb Z_{p}}(\mathbb Z_p)$, cf. [52, p. 52]. Let $\tilde G_{\mathbb Z_{(p)}}$ be the smooth, affine group scheme over $\mathbb Z_{(p)}$ that extends $G$ and whose extension to $\mathbb Z_p$ is $\tilde G_{\mathbb Z_{p}}$, cf. [54, Cl. 3.1.3.1]. Let $\tilde L_{(p)}$ be a $\mathbb Z_{(p)}$-lattice of $W$ such that the monomorphism $G\hookrightarrow \pmb{\text{GL}}_W$ extends to a homomorphism $\tilde G_{\mathbb Z_{(p)}}\rightarrow \pmb{\text{GL}}_{\tilde L_{(p)}}$, cf. [25, Part I, 10.9].

\begin{lemma}\label{L9} We can modify the $\mathbb Z$-lattice $L$ of $W$ and the injective map $f:(G,\mathcal X)\hookrightarrow (\pmb{\text{GSp}}(W,\psi),\mathcal Y)$, such that we have an identity $H=\tilde H$ and $L_{(p)}$ is a $\tilde G_{\mathbb Z_{(p)}}$-module (but we emphasize that the resulting homomorphism $\tilde G_{\mathbb Z_{(p)}}\rightarrow \pmb{\text{GL}}_{L_{(p)}}$ of smooth group schemes over $\mathbb Z_{(p)}$ is not necessarily a closed embedding).
\end{lemma}

\medskip\proof
Let $\tilde L$ be the $\mathbb Z$-lattice of $W$ such that we have $\tilde L[\frac{1}{p}]=L[\frac{1}{p}]$ and $\tilde L\otimes_{\mathbb Z} \mathbb Z_{(p)}=\tilde L_{(p)}$. If $\psi$ induces a perfect form on $\tilde L$, then by replacing $L$ with $\tilde L$ we get that $H=\tilde H$. This is so as the fact that $\tilde H$ is a maximal compact, open subgroup of $G_{\mathbb Q_p}(\mathbb Q_p)$ implies that the monomorphism $\tilde H\hookrightarrow G_{\mathbb Q_p}(\mathbb Q_p)\cap \pmb{\text{GL}}_{\tilde L\otimes_{\mathbb Z} \mathbb Z_p}(\mathbb Z_p)$ is an isomorphism. If $\psi$ does not induces a perfect form on $\tilde L$, then we will have to modify $f$ as follows. 

Let $L_1^\prime:=\tilde L\oplus \tilde L^{\vee}$. Let $W_1:=L_1^\prime\otimes_{\mathbb Z} \mathbb Q$ and $L^\prime_{1,(p)}:=L_1^\prime\otimes_{\mathbb Z} \mathbb Z_{(p)}$. Let $\psi_1^\prime$ be a perfect, alternating bilinear form on $L_1^\prime$ such that the group scheme $\pmb{\text{SL}}_{\tilde L}$, when viewed naturally as a closed subgroup scheme of $\pmb{\text{SL}}_{L_1^\prime}$, is in fact a subgroup scheme of $\pmb{\text{Sp}}(L_1^\prime,\psi_1^\prime)$. We can assume that $\tilde L$ and $\tilde L^{\vee}$ are both maximal isotropic $\mathbb Z$-lattices of $W_1$ with respect to $\psi_1^\prime$ (this automatically holds if $r>1$). Let $\tilde G^0_{\mathbb Z_{(p)}}$ be the schematic closure in $\tilde G_{\mathbb Z_{(p)}}$ of $G^0$; it is a flat, closed subgroup scheme of $\pmb{\text{SL}}_{\tilde L\otimes_{\mathbb Z} \mathbb Z_{(p)}}$ and thus also of $\pmb{\text{GSp}}(L_{1,(p)}^\prime,\psi_1^\prime)$. The subgroup scheme of $\pmb{\text{GSp}}(L_{1,(p)}^\prime,\psi_1^\prime)$ generated by $Z(\pmb{\text{GL}}_{L_{1,(p)}^\prime})$ and $\tilde G^0_{\mathbb Z_{(p)}}$ is a group scheme which is naturally identified with $\tilde G_{\mathbb Z_{(p)}}$ itself. 

Let $\mathfrak{A}$ be the free $\mathbb Z_{(p)}$-module of alternating bilinear forms on $L_1^\prime\otimes_{\mathbb Z} \mathbb Z_{(p)}$ fixed by $\tilde G_{\mathbb Z_{(p)}}^0$. There exist elements of $\mathfrak{A}\otimes_{\mathbb Z_{(p)}} \mathbb R$ that define polarizations of the Hodge $\mathbb Q$--structure on $W_1$ defined by a fixed element $h\in\mathcal X$, cf. [13, Cor. 2.3.3]. Thus the real vector space $\mathfrak{A}\otimes_{{\mathbb Z}_{(p)}} \mathbb R$ has a non-empty, open subset of such polarizations, cf. [13, Subsubsect. 1.1.18 (a)]. A standard application to $\mathfrak{A}$ of the approximation theory  for independent  valuations, shows the existence of an alternating bilinear form $\psi_1$ on $L_1^\prime\otimes_{\mathbb Z} \mathbb Z_{(p)}$ that is fixed by $\tilde G_{\mathbb Z_{(p)}}^0$, that is congruent to $\psi^\prime_1$ modulo $p$, and that defines a polarization of the mentioned Hodge $\mathbb Q$--structure. Thus there exists an injective map $f_1:(G,\mathcal X)\hookrightarrow (\pmb{\text{GSp}}(W_1,\psi_1),\mathcal Y_1)$ of Shimura pairs. 

As $\psi_1$ is congruent  to $\psi^\prime_1$ modulo $p$, it is a perfect, alternating bilinear form on $L_1^\prime\otimes_{\mathbb Z} \mathbb Z_{(p)}$. Let $L_1$ be a $\mathbb Z$-lattice of $W_1$ such that $\psi_1$ induces a perfect, alternating bilinear form on $L_1$ and we have $L_{1,(p)}:=L_1\otimes_{\mathbb Z} \mathbb Z_{(p)}=L_{1,(p)}^\prime$; thus $L_{1,(p)}$ is a $\tilde G_{\mathbb Z_{(p)}}$-module. As above we argue that $\tilde H=G_{\mathbb Q_p}(\mathbb Q_p)\cap \pmb{\text{GL}}_{L_1\otimes_{\mathbb Z} \mathbb Z_p}(\mathbb Z_p)$. Therefore the lemma holds.\endproof 

\begin{theorem}\label{T6} 
Let $\tilde H$ be a maximal compact, open subgroup of $G_{\mathbb Q_p}(\mathbb Q_p)$. Let $\tilde G_{\mathbb Z_{(p)}}$ be a smooth, affine group scheme over $\mathbb Z_{(p)}$ that has $G$ as its generic fibre and such that $\tilde H=\tilde G_{\mathbb Z_{(p)}}(\mathbb Z_p)$ (see beginning of Subsection 4.1). We assume that one of the following two conditions holds:

\medskip\noindent
{\bf (i)} the special fibre $\tilde G_{\mathbb F_{p}}$ of $\tilde G_{\mathbb Z_{p}}$ has a torus of the same rank as $G$ (e.g., this holds if $G_{\mathbb Q_p}$ splits over an unramified extension of $\mathbb Q_p$, cf. [52, Sects. 1.10 and 3.4]);

\smallskip\noindent
{\bf (ii)} we have $e(v)=1$ and $\tilde G_{\mathbb Z_{(p)}}$ is a quasi-reductive group scheme for $(G,\mathcal X,v)$.

\medskip
Then there exists a unique regular, formally smooth integral model $\tilde{\mathcal N}^{\text{s}}$ of $\text{Sh}_{\tilde H}(G,\mathcal X)$ over $O_{(v)}$ that satisfies the following smooth extension property: if $Z$ is a regular, formally smooth scheme over a discrete valuation ring $O$ which is of absolute ramification index $1$ and is an $O_{(v)}$-algebra, then each morphism $Z_{E(G,\mathcal X)}\rightarrow\tilde{\mathcal N}^{\text{s}}_{E(G,\mathcal X)}$ of $E(G,\mathcal X)$-schemes extends uniquely to a morphism $Z\rightarrow\tilde{\mathcal N}^{\text{s}}$ of $O_{(v)}$-schemes.
\end{theorem}

\medskip\proof
We can assume that the injective map $f:(G,\mathcal X)\rightarrow (\pmb{\text{GSp}}(W,\psi),\mathcal Y)$ of Shimura pairs is such that $\tilde H=H$ and $L_{(p)}$ is a $\tilde G_{\mathbb Z_{(p)}}$-module, cf. Lemma \ref{L9}. If (i) holds, then the condition (i) of Lemma \ref{L8}  holds. If (ii) holds, let $\tilde G^{\text{r}}_{\mathbb Z_p}$ be a reductive, normal, closed subgroup scheme of $\tilde G_{\mathbb Z_p}$ such that there exists a cocharacter $\mu_v:\mathbb G_{m,W(k(v))}\rightarrow \tilde G^{\text{r}}_{W(k(v))}$ with the property that the extension of $\mu_v$ to $\mathbb C$ via an (any) $O_{(v)}$-monomorphism $W(k(v))\hookrightarrow\mathbb C$ defines a cocharacter of $G_{\mathbb C}$ that is $G(\mathbb C)$-conjugate to the cocharacters $\mu_h$ ($h\in \mathcal X$) introduced in the beginning of Subsection 1.3. The group $G_{\mathbb C}^{\text{der}}$ has no simple factors that are $\pmb{\text{SO}}_{2n+1}$ groups for some $n\in\mathbb N^*$, cf. Fact \ref{F2}. Therefore the natural homomorphism $\tilde G^{\text{r}}_{\mathbb Z_p}\rightarrow \pmb{\text{GL}}_{L_{(p)}\otimes_{\mathbb Z_{(p)}} \mathbb Z_p}$ is a closed embedding, cf. [56, Thm. 1.1 (d)]. Thus $\tilde G^{\text{r}}_{\mathbb Z_p}$ is naturally a closed subgroup scheme of $G_{\mathbb Z_p}$. This implies that $G_{\mathbb Z_{(p)}}$ is also a quasi-reductive group scheme for $(G,\mathcal X,v)$. Thus, if (ii) holds, then the condition (ii) of Lemma \ref{L8} holds. 

As one of the two conditions (i) and (ii) of Lemma \ref{L8} holds, the $k(v)$-scheme $\mathcal N^{\text{s}}_{k(v)}$ is non-empty (cf. Lemma \ref{L8} (a)). Based on Theorem \ref{T1} (a) and the fact that $\tilde H=H$, we get that as $\tilde{\mathcal N}^{\text{s}}$ we can take $\mathcal N^{\text{s}}$ itself.\endproof 

\begin{corollary}\label{C3} 
{\it Let $(G,\mathcal X)$ be a Shimura pair of Hodge type. Let $v$ a prime of the reflex field $E(G,\mathcal X)$ that divides a prime $p$ with the property that the group $G_{\mathbb Q_p}$ is unramified. Then for each hyperspecial subgroup $\tilde H$ of $G_{\mathbb Q_p}(\mathbb Q_p)$, there exists a unique regular, formally smooth integral model $\tilde{\mathcal N}^{\text{s}}$ of $\text{Sh}_{\tilde H}(G,\mathcal X)$ over $O_{(v)}$ that satisfies the following smooth extension property: if $Z$ is a regular, formally smooth scheme over a discrete valuation ring $O$ which is of absolute ramification index $1$ and is an $O_{(v)}$-algebra, then each morphism $Z_{E(G,\mathcal X)}\rightarrow \tilde{\mathcal N}^{\text{s}}_{E(G,\mathcal X)}$ extends uniquely to a morphism $Z\rightarrow \tilde{\mathcal N}^{\text{s}}$ between $O_{(v)}$-schemes.}
\end{corollary}

\medskip\proof
As $\tilde H$ is a hyperspecial subgroup, we can assume that the group scheme $\tilde G_{\mathbb Z_p}$ is reductive. Thus $\tilde G_{\mathbb Z_{(p)}}$ is a reductive group scheme. Therefore the condition (i) (in fact even the condition (ii)) of Theorem \ref{T9}  holds and hence the corollary follows from Theorem \ref{T6}.\endproof  

\begin{theorem}\label{T7}
We assume that $G_{\mathbb Q_p}$ splits over an unramified extension of $\mathbb Q_p$ (thus $e(v)=1$). Then for each parahoric subgroup $\tilde H$ of $G_{\mathbb Q_p}(\mathbb Q_p)$ in the sense of [9, Def. 5.2.6], there exists a unique regular, formally smooth integral model $\tilde{\mathcal N}^{\text{s}}$ of $\text{Sh}_{\tilde H}(G,\mathcal X)$ over $O_{(v)}$ that satisfies the following smooth extension property: if $Z$ is a regular, formally smooth scheme over a discrete valuation ring $O$ which is of absolute ramification index $1$ and is an $O_{(v)}$-algebra, then each morphism $Z_{E(G,\mathcal X)}\rightarrow \tilde{\mathcal N}^{\text{s}}_{E(G,\mathcal X)}$ extends uniquely to a morphism $Z\rightarrow \tilde{\mathcal N}^{\text{s}}$ between $O_{(v)}$-schemes.
\end{theorem}

\medskip\proof
As $\tilde H$ is a parahoric subgroup of $G_{\mathbb Q_p}(\mathbb Q_p)$, it is the subgroup of $G_{\mathbb Q_p}(\mathbb Q_p)$ that fixes all vertices $v_1,\ldots,v_s$ of a facet $\tilde{\mathcal F}$ of an apartment $\tilde{\mathcal A}$ of the building of $G_{\mathbb Q_p}$ over $\mathbb Q_p$. For $i\in\{1,\ldots,s\}$, let $\tilde H_i$ be the maximal, compact subgroup subgroup of $G_{\mathbb Q_p}(\mathbb Q_p)$ that fixes $v_i$. We have:

\medskip\noindent
{\bf (i)} $\tilde H=\cap_{i=1}^s \tilde H_i$.

\medskip
Let $\tilde G_{i,\mathbb Z_{p}}$ be the smooth, affine group scheme over $\mathbb Z_{p}$ that extends $G_{\mathbb Q_p}$, that satisfies the identity $\tilde H_i=\tilde G_{i,\mathbb Z_{p}}(\mathbb Z_p)$, and that is constructed as in [52, p. 52]. Let $\tilde G_{\mathbb Z_{(p)}}$ be the flat group scheme over $\mathbb Z_{(p)}$ that extends $G$ and such that its extension to $\mathbb Z_p$ is $\tilde G_{i,\mathbb Z_{p}}$, cf. [54, Cl. 3.1.3.1].
Let $\tilde G_{\mathbb Z_{(p)}}$ be the schematic closure of $G$ embedded diagonally into the generic fibre of the product $\prod_{i=1}^s \tilde G_{i,\mathbb Z_{(p)}}$; it is a flat, affine group scheme over $\mathbb Z_{(p)}$ such that (due to property (i)) we have $\tilde G_{\mathbb Z_{(p)}}(\mathbb Z_p)=\tilde H$.

As $G_{\mathbb Q_p}$ splits over an unramified extension of $\mathbb Q_p$, there exists a maximal torus $\tilde{\mathcal T}$ of $G_{\mathbb Q_p}$ which splits over a finite unramified Galois extension $\mathbb E_p$ of $\mathbb Q_p$ and which contains the maximal, split torus of $G_{\mathbb Q_p}$ that is related to the apartment $\tilde{\mathcal A}$ in such a way that the apartment of the building of $G_{\mathbb E_p}$ related to $\tilde{\mathcal T}_{\mathbb E_p}$ contains $\tilde{\mathcal A}$ (see [52, Subsects. 1.10 and 2.6]). Let $\tilde{\mathcal T}_{\mathbb Z_p}$ be the torus over $\mathbb Z_p$ whose generic fibre is $\tilde T$; it is a maximal torus of each $\tilde G_{i,\mathbb Z_{p}}$ and therefore also of the pullback of $\tilde G_{\mathbb Z_{(p)}}$ to $\Spec(\mathbb Z_p)$.

Based on the last two paragraphs, it is easy to check that (cf. also [9, Subsect. 5.2]):

\medskip\noindent
{\bf (ii)} $\tilde G_{\mathbb Z_{(p)}}$ is a smooth group scheme over $\Spec(\mathbb Z_{(p)})$ whose special fibre has the same rank as $G$. 

\medskip
Based on Lemma \ref{L9}, for each $i\in\{1,\ldots,s\}$ there exists an injective map $f_i:(G,\mathcal X)\hookrightarrow (\pmb{\text{GSp}}(W_i,\psi_i),\mathcal Y_i)$ and a $\mathbb Z$-lattice $L_i$ of $W_i$ such that $\psi_i$ induces a perfect, alternating bilinear form $\psi_i:L_i\times L_i\rightarrow \mathbb Z$, $L_i\otimes_{\mathbb Z} \mathbb Z_{(p)}$ is a $\tilde G_{i,\mathbb Z_{(p)}}$-module, and we have $\tilde H_i=G(\mathbb Q_p)\cap \pmb{\text{GSp}}(L_i,\psi_i)(\mathbb Z_p)$. We fix an element $x\in\mathcal X$. By replacing each $\psi_i$ by either itself or $-\psi_i$ we can assume that $2\pi i\psi_i$ is a polarization of the $\mathbb Q$-structure on $W_i$ defined naturally by $x$ for all $i\in\{1,\ldots,s\}$. Let
$$(W,\psi,L)=(\oplus_{i=1}^s W_i,\oplus_{i=1}^s \psi_i,\oplus_{i=1}^s L_i).$$

We have a natural diagonal embedding $(G,\mathcal X)\hookrightarrow (\pmb{\text{GSp}}(W,\psi),\mathcal Y)$, where $\mathcal Y$ is the $\pmb{\text{GSp}}(W,\psi)(\mathbb R)$-conjugacy class of homomorphisms $\text{Res}_{\mathbb C/\mathbb R} \mathbb G_m\rightarrow \pmb{\text{GSp}}(W,\psi)_{\mathbb R}$ that contains all those homomorphisms that are defined naturally by elements of $\mathcal X$. The group $H=G(\mathbb Q_p)\cap \pmb{\text{GSp}}(L,\psi)(\mathbb Z_p)$ is $\cap_{i=1}^s \tilde H_i=\tilde H$, cf. property (i). From the construction of $\tilde G_{\mathbb Z_{(p)}}$ (as a schematic closure) we get that $L_{(p)}$ is a $\tilde G_{\mathbb Z_{(p)}}$-module. From this and from the rank part of the property (ii), we get that the condition (i) of Lemma \ref{L8} holds. Thus the $k(v)$-scheme $\mathcal N^{\text{s}}_{k(v)}$ is non-empty (cf. Lemma \ref{L8} (a)). Based on Theorem \ref{T1} (a) and the fact that $\tilde H=H$, we get that as $\tilde{\mathcal N}^{\text{s}}$ we can take $\mathcal N^{\text{s}}$ itself.\endproof 

\section{Proof of the Main Theorem \ref{T2}}
In this section we take $k$ to be a field extension of $k(v)$ that is algebraically closed and has a countable transcendental degree. Let $(v_{\alpha})_{\alpha\in\mathcal J}$ and $(w_{\alpha}^{\mathcal A})_{\alpha\in\mathcal J}$ be as in Subsection 2.3. For a point $z\in\mathcal N^{\text{s}}(W(k))=\mathcal N(W(k))$, the notation $(A,(w_{\alpha})_{\alpha\in\mathcal J},\lambda_A)$, $(M,F^1,\phi,(t_{\alpha})_{\alpha\in\mathcal J},\psi_M)$, $M=F^1\oplus F^0$, and $\mu:\mathbb G_{m,W(k)}\rightarrow\mathcal G$ is as in Subsection 3.1. Subsections 5.1 to 5.6 prove the Main Theorem \ref{T2}. 

Let $R_0:=W(k)[[x]]$, where $x$ is an independent variable. Let $\Phi_{R_0}$ be the Frobenius lift of $R_0$ that is compatible with $\sigma$ and that takes $x$ to $x^p$.  
 
\subsection{Basic notation and facts} 
We begin the proof of the Main Theorem \ref{T2}  by introducing notation and some basic facts. We have $e(v)=1$ and $G_{\mathbb Z_{(p)}}$ is a quasi-reductive group scheme for $(G,\mathcal X,v)$. We recall that $\mathcal N^{\text{m}}$ is an open subscheme of $\mathcal N^{\text{s}}$ (cf. Subsection 3.5) and therefore also of $\mathcal N$ (cf. Lemma \ref{L1}). Thus $\mathcal N^{\text{m}}_{k(v)}$ is also an open subscheme of $\mathcal N_{k(v)}$. Moreover, the open embedding $\mathcal N^{\text{m}}\hookrightarrow\mathcal N$ is a pro-finite pro-\'etale cover of an open embedding between quasi-projective $O_{(v)}$-schemes (cf. Propositions \ref{P2} (a) and \ref{P3} (a)) and the $k(v)$-scheme $\mathcal N^{\text{m}}_{k(v)}$ is non-empty (cf. Lemma \ref{L8} (a)). Thus to show that $\mathcal N^{\text{m}}_{k(v)}$ is a non-empty, open closed subscheme of $\mathcal N_{k(v)}$, we only have to show that for each commutative diagram of the following type  
\[\xymatrix{
\Spec(k)\ar[r]^{} \ar[d]^{y} & \Spec(k[[x]])  \ar[d]^{\tau} & \Spec(k((x))) \ar[l]^{} \ar[d]^{\tau_{k((x))}} \\
\mathcal N & \mathcal N_{k(v)} \ar[l]^{} & \mathcal N^{\text{m}}_{k(v)} \ar[l]^{} ,
}\]
the morphism $y:\Spec(k)\rightarrow\mathcal N$ factors through the open subscheme $\mathcal N^{\text{m}}$ of $\mathcal N$. All the horizontal arrows of this diagram are natural (closed or open) embeddings. Until Subsection 5.4 inclusive we study properties of this diagram that are needed to prove Theorem \ref{T2} in Subsections 5.4 to 5.6. 

We consider the principally quasi-polarized $F$-crystal
$$(M_0,\Phi_0,\nabla_0,\psi_{M_0})$$ 
over $k[[x]]$ of $\tau^*((\mathcal A,\lambda_{\mathcal A})\times_{\mathcal N} \mathcal N_{k(v)})$. Thus $M_0$ is a free $R_0$-module of rank $2r$, $\Phi_0$ is a $\Phi_{R_0}$-linear endomorphism of $M_0$, $\nabla_0$ is an integrable and nilpotent modulo $p$ connection on $M_0$ such that we have $\nabla_0\circ\Phi_0=(\Phi_0\otimes d\Phi_{R_0})\circ\nabla_0$, and $\psi_{M_0}$ is a perfect, alternating bilinear form on $M_0$ that defines a principal quasi-polarization of $(M_0,\Phi_0,\nabla_0)$. 

Let $O$ be the unique local ring of ${R_0}$ that is a discrete valuation ring of mixed characteristic $(0,p)$. Let $\mathcal O$ be the completion of $O$. Let $\Phi_{\mathcal O}$ be the Frobenius lift of $\mathcal O$ defined by $\Phi_{R_0}$ via a natural localization and completion. Let $k_1:=\overline{k((x))}$. Let $\Spec(W(k_1))\rightarrow \Spec(R_0)$ be the lift that is compatible with the Frobenius lifts $\sigma_{k_1}$ and $\Phi_{R_0}$; under it $W(k_1)$ gets naturally the structure of a $*$-algebra, where $*\in\{R_0,O,\mathcal O\}$. 

As the $O_{(v)}$-scheme $\mathcal N^{\text{m}}$ is formally smooth, there exists a lift $\tilde z_1:\Spec(\mathcal O)\rightarrow\mathcal N^{\text{m}}$ of the morphism $\tau_{k((x))}:\Spec(k((x)))\rightarrow\mathcal N^{\text{m}}$ defined naturally by $\tau_{k((x))}$ and denoted in the same way. Let $z_1:\Spec(W(k_1))\rightarrow \mathcal N^{\text{m}}$ be the composite of $\Spec(W(k_1))\rightarrow\Spec(\mathcal O)$ with $\tilde z_1$; we also view $\tilde z_1$ and $z_1$ as valued points of either $\mathcal N^{\text{s}}$ or $\mathcal N$. Let 
$$(\tilde A_1,(w_{1,\alpha})_{\alpha\in\mathcal J},\lambda_{\tilde A_1}):=\tilde z_1^*(\mathcal A,(w_{\alpha}^{\mathcal A})_{\alpha\in\mathcal J},\lambda_{\mathcal A})\; \text{and}\; (A_1,\lambda_{A_1}):=z_1^*(\mathcal A,\lambda_{\mathcal A})=(\tilde A_1,\lambda_{\tilde A_1})_{W(k_1)}.$$ 
\indent
For $\alpha\in\mathcal J$ let $n(\alpha)\in\mathbb N$ be such that we have $v_{\alpha}\in W^{\vee\otimes n(\alpha)}\otimes_{\mathbb Q} W^{\otimes n(\alpha)}\subset\mathcal T(W^{\vee})$, cf. definition of $v_{\alpha}$ in Subsection 2.3. Let $t_{1,\alpha}$ be the de Rham realization of $w_{1,\alpha}$. We identify canonically $M_0\otimes_{R_0} \mathcal O=H^1_{\text{dR}}(\tilde A_1/\mathcal O)$ (cf. [2, Ch. V, Subsect. 2.3]) and thus we view each $t_{1,\alpha}$ as a tensor of $(M_0^{\otimes n(\alpha)}\otimes_{R_0} M_0^{\vee\otimes n_{\alpha}})\otimes_{R_0} \mathcal O[\frac{1}{p}]\subset \mathcal T(M_0\otimes_{R_0} \mathcal O)[\frac{1}{p}]$. Let $n_{\alpha}\in\mathbb N$ be the smallest number such that we have $p^{n_{\alpha}}t_{1,\alpha}\in (M_0^{\otimes n(\alpha)}\otimes_{R_0} M_0^{\vee\otimes n_{\alpha}})\otimes_{R_0} \mathcal O\subset \mathcal T(M_0\otimes_{R_0} \mathcal O)$.

\begin{proposition}\label{P4}
For all $\alpha\in\mathcal J$ we have $p^{n_{\alpha}}t_{1,\alpha}\in M_0^{\otimes n(\alpha)}\otimes_{R_0} M_0^{\vee\otimes n(\alpha)}\subset\mathcal T(M_0)$.
\end{proposition}

\medskip\proof
The tensor $p^{n_{\alpha}}t_{1,\alpha}$ is fixed by the $\sigma_{k_1}$-linear automorphism of $\mathcal T(M_0\otimes_{R_0} B(k_1))$ defined by $\Phi_0$, cf. Subsection 3.1. Thus (as $\Spec(W(k_1))\rightarrow \Spec(R_0)$ is a Teichm\"uller lift) $p^{n_{\alpha}}t_{1,\alpha}$ is also fixed by the $\Phi_{\mathcal O}$-linear endomorphism of $\mathcal T(M_0\otimes_{R_0} \mathcal O)[\frac{1}{p}]$ defined by $\Phi_0$.  

The field $k((x))$ has $\{x\}$ as a $p$-basis, i.e., $\{1,x,\ldots,x^{p-1}\}$ is a basis of $k((x))$ over $k((x))^p=k((x^p))$. Thus the $p$-adic completion of the $\mathcal O$-module $\Omega_{\mathcal O/W(k)}$ of relative differentials is naturally isomorphic to $\mathcal O dx$, cf. [3, Prop. 1.3.1]. Let $\nabla_0:M_0\otimes_{R_0} \mathcal O\rightarrow M_0\otimes_{R_0} \mathcal O dx$ be the connection which is the natural extension of the connection $\nabla_0$ on $M_0$. 

The de Rham component of $w_{\alpha}^{\mathcal A}$ is annihilated by the Gauss--Manin connection of $\mathcal A$ (this is a property of Hodge cycles; for instance, it follows from [14, Prop. 2.5] applied in the context of a quotient of $\text{Sh}_H(G,\mathcal X)$ by a small compact, open subgroup of $G(\mathbb A_f^{(p)})$). Thus the tensor $p^{n_{\alpha}}t_{1,\alpha}$ is annihilated by the Gauss--Manin connection on $\mathcal T(H^1_{\text{dR}}(\tilde A_1/\mathcal O))=\mathcal T(M_0\otimes_{R_0} \mathcal O)$ of $\tilde A_1$ and thus also by the $p$-adic completion of this connection. Therefore $p^{n_{\alpha}}t_{1,\alpha}$ is annihilated by the connection $\nabla_0:M_0\otimes_{R_0} \mathcal O\rightarrow M_0\otimes_{R_0} \mathcal O dx$, cf. [2, Ch. V, Prop. 3.6.4]. 

As the field $k((x))$ has a $p$-basis, each $F$-crystal over $k((x))$ is uniquely determined by its evaluation at the thickening naturally associated to the closed embedding $\Spec(k((x)))\hookrightarrow\Spec(\mathcal O)$ (cf. [3, Prop. 1.3.3]). Thus the natural identification 
$$(M_0^{\otimes n(\alpha)}\otimes_{R_0} M_0^{\vee\otimes n_{\alpha}})\otimes_{R_0} \mathcal O=\text{End}(M_0^{\otimes n(\alpha)}\otimes_{R_0} \mathcal O)$$ 
allows us to view $p^{n_{\alpha}}t_{1,\alpha}$ as an endomorphism of the $F$-crystal over $k((x))$ defined by the tensor product of $n(\alpha)$-copies of $(M_0\otimes_{R_0} \mathcal O,\Phi_0\otimes\Phi_{\mathcal O},\nabla_0)$. From this and Theorem \ref{T3} we get that $p^{n_{\alpha}}t_{1,\alpha}$ can be viewed as an endomorphism of the $F$-crystal over $k[[x]]$ defined by the tensor product of $n(\alpha)$-copies of $(M_0,\Phi_0,\nabla_0)$ and therefore in fact we have $p^{n_{\alpha}}t_{1,\alpha}\in M_0^{\otimes n(\alpha)}\otimes_{R_0} M_0^{\vee\otimes n(\alpha)}\subset\mathcal T(M_0)$.\endproof

\subsubsection{Group schemes} Next we introduce notation that pertains to group schemes. Let $G^{\text{r}}_{\mathbb Z_p}$ be a reductive, normal, closed subgroup scheme of $G_{\mathbb Z_p}$ as in Definition \ref{D2} (a); we emphasize that in general $G^{\text{r}}_{\mathbb Z_p}$ is not the pullback to $\Spec(\mathbb Z_p)$ of a closed subgroup scheme of $G_{\mathbb Z_{(p)}}$. Let $\pi_{\text{r}}\in\text{End}(M_0\otimes_{R_0} B(k_1))$ be the projector that corresponds to the projector $\pi_{G^{\text{r}}_{\mathbb Q_p}}$ of Subsection 2.3 via Fontaine comparison theory for (the Barsotti--Tate group of) $\tilde A_{1,W(k_1)}$, cf. Subsection B3. As $\pi_{G^{\text{r}}_{\mathbb Q_p}}$ is fixed by $G_{\mathbb Q_p}$, by enlarging the family $(v_{\alpha})_{\alpha\in\mathcal J}$, we can assume that $\pi_{G^{\text{r}}_{\mathbb Q_p}}$ is a $\mathbb Q_p$-linear combinations of the $v_{\alpha}$'s with $\alpha\in\mathcal J$. Thus $\pi_{\text{r}}$ is a $\mathbb Q_p$-linear combination of the $t_{1,\alpha}$'s with $\alpha\in\mathcal J$. From this and Proposition \ref{P4} we get that in fact we have $\pi_{\text{r}}\in\text{End}(M_0[\frac{1}{p}])$. Thus there exists $n_{\text{r}}\in\mathbb N$ such that $p^{n_{\text{r}}}\pi_{\text{r}}\in\text{End}(M_0)$. 

Let $\eta$ be the field of fractions of $R_0$ (or of $O$). Let $\mathcal G_{0,\eta}$ be the subgroup of $(\pmb{\text{GL}}_{M_0})_{\eta}$ that fixes $p^{n_{\alpha}}t_{1,\alpha}$ for all $\alpha\in\mathcal J$ (this definition makes sense due to Proposition \ref{P4}). The group $\mathcal G_{0,B(k_1)}$ corresponds to $G_{\mathbb Q_p}$ via Fontaine comparison theory for (the Barsotti--Tate group of) $\tilde A_{1,B(k_1)}$. This implies that $\mathcal G_{0,\eta}$ is a reductive group.  

\begin{lemma}\label{L10}
There exists a unique reductive subgroup $\mathcal G_{0,\eta}^{\text{r}}$ of $\mathcal G_{0,\eta}$ whose Lie algebra is $\text{Im}(\pi_{\text{r}})\otimes_{R_0[\frac{1}{p}]} \eta$. The subgroup $\mathcal G_{0,\eta}^{\text{r}}$ of $\mathcal G_{0,\eta}$ is normal. Moreover each geometric pullback of $\mathcal G_{0,\eta}^{\text{r},\text{der}}$ has no normal subgroup which is an $\pmb{\text{SO}}_{2n+1}$ group for some $n\in\mathbb N^*$.
\end{lemma}

\medskip\proof 
From Fontaine comparison theory for (the Barsotti--Tate group of) $\tilde A_{1,W(k_1)}$ we get that there exists a unique reductive subgroup $\mathcal G_{0,B(k_1)}^{\text{r}}$ of $\pmb{\text{GL}}_{M_0\otimes_{R_0} B(k_1)}$ whose Lie algebra is $\text{Im}(\pi_{\text{r}})\otimes_{R_0[\frac{1}{p}]} B(k_1)$, cf. Lemma \ref{L15} (a). From Lemma \ref{L11} (a) applied with $(\mathcal W,\mathcal L,\eta,\eta_1)=(M_0\otimes_{R_0} \eta,\text{Im}(\pi_{\text{r}})\otimes_{R_0[\frac{1}{p}]} \eta,\eta,B(k_1))$, we get that there exists a unique reductive subgroup $\mathcal G_{0,\eta}^{\text{r}}$ of $\pmb{\text{GL}}_{M_0\otimes_{R_0} \eta}$ whose Lie algebra is $\text{Im}(\pi_{\text{r}})\otimes_{R_0[\frac{1}{p}]} \eta$. The group $\mathcal G_{0,\eta}^{\text{r}}$ is a subgroup of $\mathcal G_{0,\eta}$, as this holds after extension to $B(k_1)$. Thus the first part of the lemma holds.  

But $\pi_{\text{r}}$ is fixed by $\mathcal G_{0,\eta}$  (as this holds after tensorization with $B(k_1)$, cf. Subsection B3) and thus $\text{Im}(\pi_{\text{r}})\otimes_{R_0[\frac{1}{p}]} \eta$ is a $\mathcal G_{0,\eta}$-submodule of $\text{Lie}(\mathcal G_{0,\eta})$. From this and the uniqueness part of the lemma, we get that $\mathcal G_{0,\eta}^{\text{r}}$ is a subgroup of $\mathcal G_{0,\eta}$ normalized by $\mathcal G_{0,\eta}(\eta)$ and thus also by $\mathcal G_{0,\eta}$. As $\mathcal G_{0,B(k_1)}^{\text{r}}$ corresponds to the normal subgroup $G^{\text{r}}_{\mathbb Q_p}$ of $G_{\mathbb Q_p}$ via Fontaine comparison theory for (the Barsotti--Tate group of)  $\tilde A_{1,W(k_1)}$, from Fact \ref{F2} we get that each geometric pullback of $\mathcal G_{0,\eta}^{\text{r},\text{der}}$ has no normal subgroup which is an $\pmb{\text{SO}}_{2n+1}$ group for some $n\in\mathbb N^*$.\endproof 
 
\begin{theorem}\label{T8} 
{The schematic closure $\mathcal G_{0}^{\text{r}}$ of $\mathcal G_{0,\eta}^{\text{r}}$ in $\pmb{\text{GL}}_{M_0}$ is a reductive subgroup scheme over $\Spec(R_0)$.}
\end{theorem}
\begin{proof}
We check that if $\mathcal V$ is a local ring of $R_0$ which is a discrete valuation ring, then $\mathcal G_{0,\mathcal V}^{\text{r}}$ is a reductive group scheme over $\mathcal V$. 

We first assume that $\mathcal V=O$. As we know that $\tilde z_1\in\mathcal N^{\text{m}}(\mathcal O)$, there exist isomorphisms $(M_0\otimes_{R_0} W(k_1),(t_{1,\alpha})_{\alpha\in\mathcal J})\rightarrow (L_{(p)}^{\vee}\otimes_{\mathbb Z_{(p)}} W(k_1),(v_{\alpha})_{\alpha\in\mathcal J})$. Therefore the schematic closure of $\mathcal G^{\text{r}}_{0,B(k_1)}$ in $\pmb{\text{GL}}_{M_0\otimes_{R_0} W(k_1)}$ is isomorphic to $G^{\text{r}}_{W(k_1)}$ and thus it is a reductive group scheme over $W(k_1)$. As the natural morphism $\Spec(W(k_1))\rightarrow\Spec(\mathcal V)=\Spec(O)$ is faithfully flat, this schematic closure is $\mathcal G_{0,\mathcal V}^{\text{r}}\times_{\Spec(\mathcal V)} \Spec(W(k_1))$. Thus $\mathcal G_{0,\mathcal V}^{\text{r}}$ is a reductive group scheme over $\mathcal V$.

We now assume that $\mathcal V\neq O$, i.e., $\mathcal V$ is of equal characteristic $0$. Thus $\mathcal G_{0,\mathcal V}^{\text{r}}$ is a smooth, closed subgroup scheme of $\pmb{\text{GL}}_{M_0\otimes_{R_0} \mathcal V}$, cf. Cartier theorem. Its Lie algebra $\mathfrak{g}_{\mathcal V}$ is $(\text{Im}(p^{n_{\text{r}}}\pi_{\text{r}})\otimes_{R_0} \eta)\cap \text{End}(M_0\otimes_{R_0} \mathcal V)=\text{Im}(\pi_{\text{r}})\otimes_{R_0[\frac{1}{p}]} \mathcal V$ and thus the restriction of the trace bilinear form on $\text{End}(M_0\otimes_{R_0} \mathcal V)$ to $\mathfrak{g}_{\mathcal V}$ is perfect. From this and Lemma \ref{L11} (b) we get that the identity component of the special fibre of $\mathcal G_{0,\mathcal V}^{\text{r}}$ is a reductive group. Let $\mathcal G_{0,\mathcal V}^{\text{r},\text{id}}$ be the open subgroup scheme of $\mathcal G_{0,\mathcal V}^{\text{r}}$ whose special fibre is the identity component of the special fibre of $\mathcal G_{0,\mathcal V}^{\text{r}}$. As $\mathcal G_{0,\mathcal V}^{\text{r},\text{id}}$ is the complement in $\mathcal G_{0,\mathcal V}^{\text{r}}$ of a divisor of $\mathcal G_{0,\mathcal V}^{\text{r}}$, it is an affine $\mathcal G_{0,\mathcal V}^{\text{r}}$-scheme and thus it is an affine scheme. Therefore $\mathcal G_{0,\mathcal V}^{\text{r},\text{id}}$ is a reductive group scheme. Based on this and the second part of Lemma \ref{L10}, from [56, Thm. 1.1 (d)] we get that the homomorphism  $\mathcal G_{0,\mathcal V}^{\text{r},\text{id}}\rightarrow \pmb{\text{GL}}_{M_0\otimes_{R_0} \mathcal V}$ is a closed embedding. Thus $\mathcal G_{0,\mathcal V}^{\text{r},\text{id}}\rightarrow \mathcal G_{0,\mathcal V}^{\text{r}}$ is a closed embedding. Being also an open embedding, we conclude that $\mathcal G_{0,\mathcal V}^{\text{r},\text{id}}=\mathcal G_{0,\mathcal V}^{\text{r}}$ is a reductive, closed subgroup scheme of $\pmb{\text{GL}}_{M_0\otimes_{R_0} \mathcal V}$.

Let $\mathcal U:=\Spec(R_0)\setminus\Spec(k)$. As $\mathcal G_{0,\mathcal U}^{\text{r}}$ is a reductive, closed group scheme of ${\pmb{\text{GL}}}_{M_0,\mathcal U}$ (cf. last two paragraphs), it extends uniquely to a reductive group scheme $\tilde{\mathcal G}_{0}^{\text{r}}$ over $R_0$ (cf. [66, Thm. 1.4 (b)]). The closed embedding homomorphism $\mathcal G_{0,\mathcal U}^{\text{r}}\rightarrow \pmb{\text{GL}}_{M_0,\mathcal U}$ extends to a  closed embedding homomorphism $\tilde{\mathcal G}_{0}^{\text{r}}\rightarrow\pmb{\text{GL}}_{M_0}$, cf. [66, Prop. 5.1 (c)] and for $p=2$ cf. also the last property of Lemma \ref{L10}. Thus $\mathcal G_{0}^{\text{r}}=\tilde{\mathcal G}_{0}^{\text{r}}$ is a reductive, closed subgroup scheme of $\pmb{\text{GL}}_{M_0}$.
\end{proof}

\subsection{Applying Theorem \ref{T8}} Let $(M_1,F^1_1,\phi_1,\psi_{M_1})$ be the principally quasi-polarized filtered $F$-crystal over $k_1$ of $(A_1,\lambda_{A_1})$. We have $M_1=M_0\otimes_{R_0} W(k_1)$, $\phi_1=\Phi_0\otimes\sigma_{k_1}$, and each $t_{1,\alpha}\in\mathcal T(M_1)[\frac{1}{p}]$ with $\alpha\in\mathcal J$ is the de Rham realization of the Hodge cycle $z_1^*(w_{\alpha}^{\mathcal A})$ on $A_{1,B(k)}$. Let $\mu_1:\mathbb G_{m,W(k_1)}\rightarrow \mathcal G_1=\mathcal G_{0,W(k_1)}$ be the analogue of $\mu:\mathbb G_{m,W(k)}\rightarrow\mathcal G$ but obtained working with $z_1\in\mathcal N(W(k_1))$ instead of some $z\in\mathcal N(W(k))$. We know that $\mu_1$ factors through $\mathcal G^{\text{r}}_{0,W(k_1)}$, cf. Lemma \ref{L6} applied to $z_1\in\mathcal N(W(k_1))$ with $G^{\text{v}}_{\mathbb Q_p}=G_{\mathbb Q_p}^{\text{r}}$. 

Let $\bar F^1_0$ be the kernel of $\Phi_0$ modulo $p$; it is a free module over $k[[x]]=R_0/pR_0$ of rank $r$. As the cocharacter $\mu_1$ factors through $\mathcal G^{\text{r}}_{0,W(k_1)}$, the normalizer of $\bar F^1_0\otimes_{k[[x]]} k_1$ in $\mathcal G^{\text{r}}_{0,k_1}$ is a parabolic subgroup of the reductive group $\mathcal G^{\text{r}}_{0,k_1}$ which (as $\bar F^1_0\otimes_{k[[x]]} k_1$ is defined over $k((x))$) is the pullback of a parabolic subgroup $\mathcal F^{\text{r}}_{0,k((x))}$ of $\mathcal G^{\text{r}}_{0,k((x))}$. The $k[[x]]$-scheme of parabolic subgroup schemes of the reductive group scheme $\mathcal G^{\text{r}}_{0,k[[x]]}$ is projective, cf. [15, Vol. III, Exp. XXVI, Cor. 3.5].  Thus the schematic closure $\mathcal F^{\text{r}}_{0,k[[x]]}$ of $\mathcal F^{\text{r}}_{0,k((x))}$ in $\mathcal G^{\text{r}}_{0,k[[x]]}$ is a parabolic subgroup scheme of $\mathcal G^{\text{r}}_{0,k[[x]]}$. 

As $\mathcal G^{\text{r}}_0$ is a split reductive group scheme and $\mu_{1,k_1}$ factors through $\mathcal G^{\text{r}}_{0,k_1}$, there exists a cocharacter $\mu_{0,k[[x]]}:\mathbb G_{m,k[[x]]}\rightarrow\mathcal G^{\text{r}}_{0,k[[x]]}$ that factors through $\mathcal F^{\text{r}}_{0,k[[x]]}$ and that produces a direct sum decomposition $M_0/pM_0=\bar F^1_0\oplus \bar F^0_0$ such that for each $i\in\{0,1\}$, every $\beta\in\mathbb G_{m,k[[x]]}(k[[x]])$ acts via $\mu_{0,k[[x]]}$ on $\bar F^i_0$ as the multiplication by $\beta^{-i}$. We consider a cocharacter
$$\mu_0:\mathbb G_{m,R_0}\rightarrow \mathcal G^{\text{r}}_0$$ 
that lifts $\mu_{0,k[[x]]}$, cf. [15, Vol. II, Exp. IX, Thms. 3.6 and 7.1]. Let $M_0=F^1_0\oplus F^0_0$ be the direct sum decomposition such that for each $i\in\{0,1\}$, every element $\beta\in\mathbb G_{m,R_0}(R_0)$ acts via $\mu_0$ on $F^i_0$ as the multiplication by $\beta^{-i}$; the notation matches, i.e., we have $F^i_0/pF^i_0=\bar F^i_0$. 

We consider the $W(k)$-epimorphism $R_0\twoheadrightarrow W(k)$ whose kernel is the ideal $(x)$. Let 
$$(M,F^1,\phi,\mathcal G,\mathcal G^{\text{r}},(t_{\alpha})_{\alpha\in\mathcal J},\psi_{M}):=(M_0,F^1_0,\Phi_0,\mathcal G_0,\mathcal G^{\text{r}}_0,(t_{1,\alpha})_{\alpha\in\mathcal J},\psi_{M_0})\otimes_{R_0} W(k).$$

\subsection{Extra crystalline applications} Let $(D,\lambda_D)$ be an arbitrary principally quasi-polarized Barsotti--Tate group over $W(k)$ whose principally quasi-polarized filtered $F$-crystal over $k$ is $(M,F^1,\phi)$ and for which we have an isomorphism 
$$(M,(t_{\alpha})_{\alpha\in\mathcal J},\psi_M)\rightarrow (H^1(D)\otimes_{\mathbb Z_p} W(k),(u_{\alpha})_{\alpha\in\mathcal J},\psi_{H^1(D)}),\leqno (6)$$
where $\psi_{H^1(D)}$ is the perfect, alternating bilinear form on $H^1(D)$ which is the \'etale realization of $\lambda_D$ and where $u_{\alpha}\in\mathcal T(H^1(D))[\frac{1}{p}]$ corresponds to $t_{\alpha}$ via Fontaine comparison theory for $D$.  If $p=2$, then the existence of $(D,\lambda_D)$ is implied by Theorem \ref{T11} (b) applied to $(M,\phi,\mathcal G^{\text{r}},\psi_M)$ instead of $(M,\phi,\mathcal G,\psi_M)$. If $p>2$ or if $p=2$ and $(M,\phi)$ has no integral slopes, then there exists a unique Barsotti--Tate group $D$ over $W(k)$ whose filtered $F$-crystal over $k$ is $(M,F^1,\phi)$ (cf. [64, Prop. 2.2.6] for $p=2$); due to the uniqueness part, $\psi_M$ is the crystalline realization of a (unique) principal quasi-polarization $\lambda_D$ of $D$. The fact that (6) holds in this case follows from Theorem \ref{T9}. 

Let $(D_{R_0},\lambda_{D_{R_0}})$ be the principally quasi-polarized Barsotti--Tate group over $R_0$ which modulo the ideal $(x)$ is $(D,\lambda_D)$ and whose principally quasi-polarized $F$-crystal over $R_0/pR_0$ is the quintuple $(M_0,F^1_0,\Phi_0,\nabla_0,\psi_{M_0})$, cf. Lemmas \ref{L16} and \ref{L17}. Let
$$\tau_{R_0}:\Spec(R_0)\rightarrow\mathcal M_r$$ 
be the morphism that has the following two properties: 

\medskip\noindent
{\bf (i)} it lifts the composite of $y:\Spec(k)\rightarrow\mathcal N$ with the morphism $\mathcal N\rightarrow\mathcal M_r$, and 

\smallskip\noindent
{\bf (ii)} the principally quasi-polarized Barsotti--Tate group of the pullback  via $\tau_{R_0}$ of the universal principally polarized abelian scheme over $\mathcal M_r$, is $(D_{R_0},\lambda_{D_{R_0}})$. 

\medskip
Let 
$$z_2:\Spec(W(k_1))\rightarrow\mathcal M_r$$ 
be the composite of the morphism $\Spec(W(k_1))\rightarrow \Spec(R_0)$ of Subsection 5.1 with $\tau_{R_0}$. Let $(A_2,\lambda_{A_2})$ be the principally polarized abelian scheme over $W(k_1)$ that is the pullback through $z_2$ of the universal principally polarized abelian scheme over $\mathcal M_r$ and let $(D_2,\lambda_{D_2})$ be its principally quasi-polarized Barsotti--Tate group. The principally quasi-polarized filtered $F$-crystal of $(D_2,\lambda_{D_2})$ is canonically identified with $(M_1,F^1_2,\phi_1,\psi_{M_1})$, where $F^1_2$ is a direct summand of $M_1$ of rank $r$. Let $(F^i_2(\mathcal T(M_1)))_{i\in\mathbb Z}$ be the filtration of $\mathcal T(M_1)$ defined by $F^1_2$ and let $(F^i_0(\mathcal T(M_0)))_{i\in\mathbb Z}$ be the filtration of $\mathcal T(M_0)$ defined by $F^1_0$. For each $\alpha\in\mathcal J$, the tensor $t_{1,\alpha}\in\mathcal T(M_0)[\frac{1}{p}]$ is annihilated by $\nabla_0$, is fixed by $\Phi_0$, and belongs to $F^0_0(\mathcal T(M_0))[\frac{1}{p}]$. This implies that we have $t_{1,\alpha}\in F^0_2(\mathcal T(M_1))[\frac{1}{p}]$ for all $\alpha\in\mathcal J$. Thus as before Lemma \ref{L5} we argue that the inverse of the canonical split cocharacter of $(M_1,F^1_2,\phi_1)$ defined in [69, p. 512] factors through the closed subgroup scheme $\mathcal G_1=\mathcal G_{0,W(k_1)}$ of $\pmb{\text{GL}}_{M_1}$; let $\mu_2:\mathbb G_{m,W(k_1)}\rightarrow \mathcal G_1$ be the resulting factorization. 

Due to Lemma \ref{L7} (a) applied to $z_1\in\mathcal N^{\text{m}}(W(k_1))\subset\mathcal N^{\text{s}}(W(k_1))$ and to $\mu_2:\mathbb G_{m,W(k_1)}\rightarrow \mathcal G_1$, there exists a point $z_3\in\mathcal N^{\text{m}}(W(k_1))\subset \mathcal N^{\text{s}}(W(k_1))=\mathcal N(W(k_1))$ that lifts the $k_1$-valued point $y_1$ of $\mathcal N^{\text{m}}$ defined naturally by $z_1$ (or $z_2$) and such that the filtered $F$-crystal of $(A_3,\lambda_{A_3}):=z_3^*(\mathcal A,\lambda_{\mathcal A})$ is precisely $(M_1,F^1_2,\phi_1,\psi_{M_1})$. Let $(D_3,\lambda_{D_3})$ be the principally quasi-polarized Barsotti--Tate group of $(A_3,\lambda_{A_3})$.

\subsection{Proof of Theorem \ref{T2} (a), part I} In this subsection we assume that either $p>2$ or $p=2$ and the $2$-rank of $y_1^*(\mathcal A)=A_{1,k_1}=\tilde A_{1,k_1}$ is $0$. Due to this assumption, the Barsotti--Tate groups $D_2$ and $D_3$ are the same lift of the Barsotti--Tate group of $y_1^*(\mathcal A)$ (cf. [64, Prop. 2.2.6] for $p=2$). Therefore the $W(k_1)$-valued points of $\mathcal M_r$ defined by $z_2$ and $z_3$ coincide. Thus $z_2$ factors through $\mathcal N^{\text{s}}$. From this and Theorem \ref{T1} (b) we get that $\tau_{R_0}$ factors through $\mathcal N^{\text{s}}$. Let $z:\Spec(W(k))\rightarrow\mathcal N^{\text{s}}$ be the composite of the factorization $\Spec(R_0)\rightarrow\mathcal N^{\text{s}}$ of $\tau_{R_0}$ with the closed embedding $\Spec(W(k))\hookrightarrow \Spec(R_0)$ defined by the $W(k)$-epimorphism $R_0\twoheadrightarrow R_0/(x)=W(k)$; it lifts $y$. The notation matches with the one of Subsection 3.1: $(D,\lambda_D)$ is the principally quasi-polarized Barsotti--Tate group of $(A,\lambda_A):=z^*(\mathcal A,\lambda_{\mathcal A})$, the principally quasi-polarized filtered $F$-crystal of $(D,\lambda_D)$ is $(M,F^1,\phi,\psi_M)$, and $u_{\alpha}$ corresponds to $t_{\alpha}$ via Fontaine comparison theory for $D$.

There exists an isomorphism $(M\otimes_{W(k)} W(k_1),(t_{\alpha})_{\alpha\in\mathcal J},\psi_M)\rightarrow (M_1,(t_{1,\alpha})_{\alpha\in\mathcal J},\psi_{M_1})$, cf. proof of Proposition \ref{P3} (b). Thus, as Theorem \ref{T4} (a) holds for $z_1\in\mathcal N^{\text{m}}(W(k_1))$, there exist isomorphisms $(M\otimes_{W(k)} W(k_1),(t_{\alpha})_{\alpha\in\mathcal J},\psi_M)\rightarrow (L^{\vee}_{(p)}\otimes_{\mathbb Z_{(p)}} W(k_1),(v_{\alpha})_{\alpha\in\mathcal J},\psi^{\vee})$ and therefore (using an argument similar to the one used to prove Lemma \ref{L14} we get that) there exist isomorphisms $(M,(t_{\alpha})_{\alpha\in\mathcal J},\psi_M)\rightarrow (L^{\vee}_{(p)}\otimes_{\mathbb Z_{(p)}} W(k),(v_{\alpha})_{\alpha\in\mathcal J},\psi^{\vee})$. From this and Lemma \ref{L3} (a) we get that Theorem \ref{T4} (a) holds for $z\in\mathcal N^{\text{s}}(W(k))$. Thus we have $z\in\mathcal N^{\text{m}}(W(k))$ (cf. Proposition \ref{P3} (b)) and therefore the morphism $y:\Spec(k)\rightarrow\mathcal N$ factors through $\mathcal N^{\text{m}}$. This ends the proof of Theorem \ref{T2} (a) provided either $p>2$ or $p=2$ and the $2$-rank of $A_{1,k_1}$ is $0$.\endproof

\subsection{Proof of Theorem \ref{T2} (a), part II}  We prove Theorem \ref{T2} (a) for $p=2$ in the general case. Let $a\in\mathbb N$ be the multiplicity of the Newton polygon slope $-1$ for $(\text{Lie}(\mathcal G_1)[{1\over 2}],\phi_1)$. Let $u_{2,\alpha}\in\mathcal T(H^1(D_2))[{1\over 2}]=\mathcal T(H^1(D_3))[{1\over 2}]$ correspond to $t_{1,\alpha}$ and let $\psi_{H^1(D_2)}$ correspond to $\psi_{M_1}$ via Fontaine comparison theory for $D_2$. 

Let $Z_0=\Spec(R_0)$. From Theorem \ref{T10} we get that the sextuple 
$$(M_0,F^1_0,\Phi_0,\nabla_0,(t_{1,\alpha})_{\alpha\in\mathcal J},\psi_{M_0})$$ 
is the pullback of the sextuple $(M_R,F^1_R,\Phi,\nabla,(t_{\alpha})_{\alpha\in\mathcal J},\psi_{M_R})$ of Susubbsection B4.1 via a morphism $i_{Z_0}: Z_0\rightarrow\Spec(R)$ of $\Spec(W(k))$-schemes which lifts the identity $R/(x_1,\ldots,x_l)=W(k)=R_0/(x)$. This implies that there exists an isomorphism 
$$(M\otimes_{W(k)} W(k_1),(t_{\alpha})_{\alpha\in\mathcal J},\psi_M)=(M_1,(t_{1,\alpha})_{\alpha\in\mathcal J},\psi_{M_1}).\leqno (7)$$
We have a canonical identification
$$(H^1(D_2),(u_{2,\alpha})_{\alpha\in\mathcal J})=(H^1(D),(u_{\alpha})_{\alpha\in\mathcal J}),\leqno (8)$$ 
cf. [64, Subsubsect. 3.4.2] (to be compared also with either the fifth paragraph of the proof of Theorem \ref{T11} or [64,  Lem. 3.4.3]).
From (6), (7), and (8) we get that there exists an isomorphism 
$$(M_1,(t_{1,\alpha})_{\alpha\in\mathcal J},\psi_M)\rightarrow (H^1(D_2)\otimes_{\mathbb Z_2} W(k_1),(u_{2,\alpha})_{\alpha\in\mathcal J},\epsilon_2\psi_{H^1(D_2)}),\leqno (9)$$
with $\epsilon_2$ a unit of $W(k_1)$. As in the proof of Theorem \ref{T4} we argue that we can assume that $\epsilon_2=1$.

We have exactly $2^a$ possibilities for a lift $z_3\in\mathcal N^{\text{m}}(W(k_1))$ of $y_1$ as in the end of Subsection 5.3 (cf. Lemma \ref{L7} (c)) for which Theorem \ref{T4} (a) holds and each such $z_3$ is uniquely determined by $(D_3,\lambda_{D_3})$ (cf. Theorem \ref{T1} (b)). Theorem \ref{T11} (c) proves that there are exactly $2^a$ possibilities for $(D_2,\lambda_{D_2})$ such as (9) holds with $\epsilon_2=1$ and its principally quasi-polarized filtered $F$-crystal over $k$ is $(M_1,F^1_2,\phi_1,\psi_{M_1})$. From the last two sentences we get that we can choose $z_3$ such that we have $(D_3,\lambda_{D_3})=(D_2,\lambda_{D_2})$ as lifts of the principally quasi-polarized $2$-divisible group of $y_1^*(\mathcal A,\lambda_{\mathcal A})$. Therefore the $W(k_1)$-valued points of $\mathcal M_r$ defined by $z_2$ and $z_3$ coincide and thus, as in Subsection 5.4 we argue that $\tau_{R_0}$ factors through $\mathcal N^{\text{s}}$ and that $y:\Spec(k)\rightarrow\mathcal N$ factors through $\mathcal N^{\text{m}}$. Thus Theorem \ref{T2} (a) holds.\endproof
  
\subsection{Proof of Theorem \ref{T2} (b) and (c)} Theorem \ref{T2} (b) follows from Theorem \ref{T2} (a) and Proposition \ref{P1}. To prove Theorem \ref{T2} (c), let $\mathcal Q$ and $\mathcal Q^{\text{s}}$ be as in Subsection 2.2. As the $\mathbb Q$-rank of the adjoint group $G^{\text{ad}}$ is $0$, $\mathcal Q$ is a projective $O_{(v)}$-scheme (cf. Lemma \ref{L2}).  From Proposition \ref{P3} (a) we get that $\mathcal N^{\text{m}}$ is the pullback of a smooth, open subscheme $\mathcal Q^{\text{m}}$ of $\mathcal Q$. To prove Theorem \ref{T2} (c.i), it suffices to show that $\mathcal Q^{\text{m}}=\mathcal Q$, i.e., to show that if $\mathcal C$ is a connected component of $\mathcal Q_{W(k)}$, then we have $\mathcal C\subset\mathcal Q^{\text{m}}$.  As $\mathcal C_{B(k)}\subset \mathcal C\cap\mathcal Q^{\text{m}}_{W(k)}$, from Lemma \ref{L8} (b) we get that the intersection $\mathcal C_k\cap\mathcal Q^{\text{m}}_k$ is non-empty and thus (as $\mathcal Q^{\text{m}}$ is smooth) there exist $W(k)$-valued points of $\mathcal C$. Thus the ring of global functions of the connected, flat, normal, projective $W(k)$-scheme $\mathcal C$ is $W(k)$. From this and [24, Ch. III, Cor. 11.3] we get that the special fibre $\mathcal C_k$ of $\mathcal C$ is connected. But the non-empty scheme $\mathcal C_k\cap \mathcal Q^{\text{m}}_k$ is an open, closed subscheme of $\mathcal C_k$, cf. Theorem \ref{T2} (a). From the last two sentences we get that $\mathcal Q^{\text{m}}_k\cap\mathcal C_k=\mathcal C_k$. Thus $\mathcal Q^{\text{m}}_{W(k)}\cap\mathcal C=\mathcal C$. Therefore Theorem \ref{T2} (c.i) holds.

We know that $\mathcal Q=\mathcal N/H^{(p)}=\mathcal N^{\text{s}}/H^{(p)}$ is a normal, projective $O_{(v)}$-scheme and that the quotient morphism $\mathcal N^{\text{s}}=\mathcal N\rightarrow\mathcal Q$ is a pro-\'etale cover, cf. beginning of Subsection 2.2 and Proposition \ref{P2} (a). Thus the $O_{(v)}$-scheme $\mathcal Q$ is smooth, i.e., we have $\mathcal Q=\mathcal Q^{\text{s}}$. As $\mathcal Q$ is a N\'eron model of its generic fibre $\text{Sh}_{H\times H^{(p)}}(G,\mathcal X)$ over $O_{(v)}$ (cf. Theorem \ref{T1} (c)), Theorem \ref{T2} (c.ii) holds.\endproof

\medskip\noindent
{\bf Acknowledgement.}
We would like to thank Universities of Utah and Arizona, MPI, Bonn, SUNY Binghamton, and IAS, Princeton for good working conditions. This research was partially supported by the NSF grants DMF 97-05376 and DMS \#0900967.

\noindent{\large \bf APPENDIX}

\begin{appendix}
\section{On affine group schemes}

Let $p\in\mathbb N$ be a prime. Let $k$ be an algebraically field of characteristic $p$. Let $W(k)$ be the ring of $p$-typical Witt vectors with coefficients in $k$ and let $B(k):=W(k)[\frac{1}{p}]$ be its field of fractions. 

\subsection{Universal smoothenings}

Let $\mathcal G$ be a flat, affine group scheme over $W(k)$. For $a\in\mathcal G(W(k))$, the N\'eron measure of the defect of smoothness $\delta(a)\in\mathbb N$ of $\mathcal G$ at $a$ is the length of the torsion part of the finitely generated $W(k)$-module $a^*(\Omega_{\mathcal G/\Spec(W(k))})$. As $\mathcal G$ is a group scheme over $W(k)$, the value of $\delta(a)$ does not depend on $a$ and thus we denote it by $\delta(\mathcal G)$. We have $\delta(\mathcal G)\in\mathbb N^*$ if and only if $\mathcal G$ is not smooth, cf. [7, Ch. 3, Sect. 3.3, Lem. 1]. Let $\mathcal F_k$ be the schematic closure in $\mathcal G_k$ of all special fibres of $W(k)$-valued points of $\mathcal G$; it is a reduced subgroup of $\mathcal G_k$. We write $\mathcal F_k=\Spec(R_{\mathcal G}/J_{\mathcal G})$, where $\mathcal G=\Spec(R_{\mathcal G})$ and where $J_{\mathcal G}$ is the ideal of $R_{\mathcal G}$ that defines $\mathcal F_k$ and contains $p$. By the {\it canonical dilatation} of $\mathcal G$ we mean the affine $\mathcal G$-scheme $\mathcal G_1=\Spec(R_{\mathcal G_1})$, where $R_{\mathcal G_1}$ is the $R_{\mathcal G}$-subalgebra of $R_{\mathcal G}[\frac{1}{p}]$ generated by ${x\over p}$ with $x\in J_{\mathcal G}$. 

The $W(k)$-scheme $\mathcal G_1$ has a canonical group scheme structure and the morphism $\mathcal G_1\rightarrow \mathcal G$ is a homomorphism of group schemes over $W(k)$, cf. [7, Ch. 3, Sect. 3.2, Prop. 2 (d)]. Moreover the $W(k)$-morphism $\mathcal G_1\rightarrow \mathcal G$ has the following universal property:  each $W(k)$-morphism $Z\rightarrow \mathcal G$ of flat $W(k)$-schemes whose special fibre $Z_k\rightarrow \mathcal G_k$ factors through the closed embedding $\mathcal F_k\hookrightarrow \mathcal G_k$, factors uniquely through $\mathcal G_1\rightarrow \mathcal G$ (cf. [7, Ch. 3, Sect. 3.2, Prop. 1 (b)]). If $\mathcal G$ is smooth, then $\mathcal F_k=\mathcal G_k$ and therefore $\mathcal G_1=\mathcal G$. 

Either $\mathcal G_1$ is smooth or we have $0<\delta(\mathcal G_1)<\delta(\mathcal G)$, cf. [7, Ch. 3, Sect. 3.3, Prop. 5]. Thus by using a sequence of at most $\delta(\mathcal G)$ canonical dilatations (the first one of $\mathcal G$, the second one of $\mathcal G_1$, etc.), we get the existence of a unique smooth, affine group scheme $\mathcal G^\prime$ over $W(k)$ endowed with a homomorphism $\mathcal G^\prime\rightarrow \mathcal G$ whose generic fibre over $B(k)$ is an isomorphism and which has the following universal property:  each $W(k)$-morphism $Z\rightarrow \mathcal G$, with $Z$ a smooth $W(k)$-scheme, factors uniquely through $\mathcal G^\prime\rightarrow \mathcal G$. One calls $\mathcal G^\prime$ the {\it universal smoothening} of $\mathcal G$.

\subsection{On Lie algebras}

\begin{lemma} \label{L11} 
{\it Let $\mathcal W$ be a finite dimensional vector space over a field $\eta$ of characteristic $0$. Let $\mathcal L$ be a Lie subalgebra of $\text{End}(\mathcal W)$. We assume that there exists a field extension $\eta_1$ of $\eta$ such that $\mathcal L\otimes_{\eta} \eta_1$ is the Lie algebra of a connected (resp. reductive) subgroup $\mathcal F_{\eta_1}$ of $\pmb{\text{GL}}_{\mathcal W\otimes_{\eta} \eta_1}$.

\medskip
{\bf (a)} Then there exists a unique connected (resp. reductive)  subgroup $\mathcal F$ of $\pmb{\text{GL}}_\mathcal W$ whose Lie algebra is $\mathcal L$ (the notation matches, i.e., the extension of $\mathcal F$ to $\eta_1$ is $\mathcal F_{\eta_1}$).

\smallskip
{\bf (b)} The restriction $\mathfrak{t}:\mathcal L\times\mathcal L\rightarrow \eta$ of the trace bilinear form on $\text{End}(\mathcal W)$ to $\mathcal L$ is non-degenerate if and only if $\mathcal F$ is a reductive subgroup of $\pmb{\text{GL}}_\mathcal W$.}
\end{lemma}

\medskip\proof 
We prove part (a). The uniqueness part is implied by [5, Ch. I, Sect. 7.1]. Loc. cit. also implies that if $\mathcal F$ exists, then its extension to $\eta_1$ is indeed $\mathcal F_{\eta_1}$. It suffices to prove part (a) for the case when $\mathcal F_{\eta_1}$ is connected. We consider commutative $\eta$-algebras $\kappa$ such that there exists a closed subgroup scheme $\mathcal F_{\kappa}$ of $\pmb{\text{GL}}_{\mathcal W\otimes_{\eta} \kappa}$ whose Lie algebra is $\mathcal L\otimes_{\eta} \kappa$. Our hypotheses imply that as $\kappa$ we can take $\eta_1$. Thus as $\kappa$ we can also take a finitely generated $\eta$-subalgebra of $\eta_1$. By considering the reduction modulo a maximal ideal of this last $\eta$-algebra, we can assume that $\kappa$ is a finite field extension of $\eta$ which is separable as $\eta$ has characteristic $0$. Thus we can assume that $\kappa$ is a finite Galois extension of $\eta$. By replacing $\mathcal F_{\kappa}$ with its identity component, we can assume that $\mathcal F_{\kappa}$ is connected. Due to the mentioned uniqueness part, the Galois group $\text{Gal}(\kappa/\eta)$ acts naturally on the connected subgroup $\mathcal F_{\kappa}$ of $\pmb{\text{GL}}_{\mathcal W\otimes_{\eta} \kappa}$. As $\mathcal F_{\kappa}$ is an affine scheme, the resulting Galois descent datum on $\mathcal F_{\kappa}$ with respect to $\text{Gal}(\kappa/\eta)$ is effective (cf. [7, Ch. 6, Sect. 6.1, Thm. 5]). This implies the existence of a subgroup $\mathcal F$ of $\pmb{\text{GL}}_\mathcal W$ whose extension to $\kappa$ is $\mathcal F_{\kappa}$. As $\text{Lie}(\mathcal F)\otimes_{\eta} \kappa=\text{Lie}(\mathcal F_{\kappa})=\mathcal L\otimes_{\eta} \kappa$, we have $\text{Lie}(\mathcal F)=\mathcal L$. The group $\mathcal F$ is connected as $\mathcal F_{\kappa}$ is so. Thus $\mathcal F$ exists, i.e., part (a) holds.

Part (b) follows from [8, Ch. I, Sect. 6, Prop. 5 and Thm. 4]. For the sake of completeness we include here a short proof of part (b). We can assume that $\eta$ is algebraically closed. We first prove the if part. Using isogenies, we are reduced to the case when $\mathcal F$ is either $\mathbb G_{m,\eta}$ or a semisimple group whose adjoint is simple. If $\mathcal F$ is $\mathbb G_{m,\eta}$, then the $\mathcal F$-module $\mathcal W$ is a direct sum of one dimensional $\mathcal F$-modules. We easily get that there exists an element $x\in\mathcal L\setminus\{0\}$ which is a semisimple element of $\text{End}(\mathcal W)$ whose eigenvalues are integers. The trace of $x^2$ is a sum of squares of natural numbers not all zero and thus it is non-zero. If $\mathcal F$ is a semisimple group whose adjoint is simple, then $\mathcal L$ is a simple Lie algebra over $\eta$. From Cartan solvability criterion we get that $\mathfrak{t}$ is non-zero, i.e., the ideal $\Ker(\mathfrak{t})=\{x\in\mathcal L|\mathfrak{t}(x,y)=0\;\forall\;y\in\mathcal L\}$ of $\mathcal L$ is not $\mathcal L$. As $\mathcal L$ is a simple Lie algebra over $\eta$, we get that $\Ker(\mathfrak{t})=0$, i.e., $\mathfrak{t}$ is non-degenerate. 

To prove the only if part, we consider the unipotent radical $\mathcal U$ of $\mathcal F$. Let $0=\mathcal W_0\subsetneqq \mathcal W_1\subsetneqq\cdots\subsetneqq\mathcal W_s=\mathcal W$ be a strictly increasing filtration of $\mathcal W$ by $\mathcal F$-modules such that $\mathcal U$ acts trivially on $\mathcal W_i/\mathcal W_{i-1}$ for all $i\in\{1,\ldots,s\}$. Based on the existence of this filtration, it is easy to see that $\text{Lie}(\mathcal U)\subset\Ker(\mathfrak{t})=0$ and thus $\text{Lie}(\mathcal U)=0$. Therefore $\mathcal U$ is the trivial subgroup, i.e., $\mathcal F$ is reductive. Thus part (b) holds.\endproof

\medskip
See [60, Prop. 3.2] for a different approach to prove Lemma \ref{L11} (a). 

\section{Complements on Barsotti--Tate groups}

Let $p$, $k$, $W(k)$, and $B(k)$ be as in Appendix A. Let $\sigma:=\sigma_k$ be the Frobenius automorphism of $k$, $W(k)$, and $B(k)$. We fix an algebraic closure $\overline{B(k)}$ of $B(k)$. Let $\text{Gal}(B(k)):=\text{Gal}(\overline{B(k)}/B(k))$. Let $D$ be a Barsotti--Tate group over $W(k)$. Let $\flat^{\text{t}}$ be the Cartier dual of a Barsotti--Tate group $\flat$ over a $W(k)$-algebra. Let $(M,\phi)$ be the $F$-crystal of $D_k$ (i.e., the contravariant Dieudonn\'e module of $D_k$ with the Verschiebung map
suppressed). Thus $M$ is a free $W(k)$-module of rank equal to the height of $D$ and $\phi:M\rightarrow M$ is a $\sigma$-linear endomorphism such that we have $pM\subset\phi(M)$. Let $F^1$ be the direct summand of $M$ that is the Hodge filtration defined by $D$. We have $\phi(M+\frac{1}{p}F^1)=M$. The rank of $F^1$ is the dimension of $D_k$. Let $M^{\vee}:=\text{Hom}(M,W(k))$. Let $\mathcal T(M)$ and its filtration $(F^i(\mathcal T(M)))_{i\in\mathbb Z}$ defined by $F^1$, be as in Subsection 2.1. For $f\in M^{\vee}[\frac{1}{p}]$ let $\phi(f):=\sigma\circ f\circ\phi^{-1}\in M^{\vee}[\frac{1}{p}]$. Thus $\phi$ acts in the usual tensor product way on $\mathcal T(M[\frac{1}{p}])$. 

If $D$ has a principal quasi-polarization $\lambda_D$, let $\psi_M:M\times M\rightarrow W(k)$ be the perfect, alternating form defined by $\lambda_D$. For all $a,b\in M$ we have $\psi_M(\phi(a),\phi(b))=p\sigma(\psi_M(a,b))$. Moreover, we have $\psi_M(F^1,F^1)=0$. 

\subsection{Galois modules} 

Let $H^1(D):=T_p(D^{\text{t}}_{B(k)})(-1)$ be the dual of the Tate-module $T_p(D_{B(k)})$ of $D_{B(k)}$. Thus $H^1(D)$ is a free $\mathbb Z_p$-module of the same rank as $M$ on which $\text{Gal}(B(k))$ acts. If $D$ has a principal quasi-polarization $\lambda_D$, let $\psi_{H^1(D)}:H^1(D)\times H^1(D)\rightarrow\mathbb Z_p$ be the perfect, alternating form defined by $\lambda_D$. Let $F^0(H^1(D)):=H^1(D)$ and $F^1(H^1(D)):=0$. Let 
$$\rho_D:\text{Gal}(B(k))\rightarrow \pmb{\text{GL}}_{H^1(D)}(\mathbb Z_p)$$ 
be the natural Galois representation associated to $D_{B(k)}$. Let $\mathcal D^{\acute et}$ be the schematic closure in $\pmb{\text{GL}}_{H^1(D)}$ of $\text{Im}(\rho_D)$; it is a flat, affine group scheme over $\mathbb Z_p$. From [69, Prop. 4.2.3] one gets that the generic fibre $\mathcal D^{\acute et}_{\mathbb Q_p}$ is connected. See Subsection 2.1 for $\mathcal T(H^1(D))$; it is naturally a $\text{Gal}(B(k))$-module. By an {\it \'etale Tate-cycle} on $D_{B(k)}$ we mean a tensor of $\mathcal T(H^1(D[\frac{1}{p}]))=\mathcal T(H^1(D))[\frac{1}{p}]$ that is fixed by $\text{Gal}(B(k))$ (equivalently, by $\mathcal D^{\acute et}_{\mathbb Q_p}$). In what follows we will fix a family $(v_{\alpha})_{\alpha\in\mathcal J}$ of \'etale Tate-cycles on $D_{B(k)}$. Let $\mathcal G^{\acute et}$ be the schematic closure in $\pmb{\text{GL}}_{H^1(D)}$ of the subgroup of $\pmb{\text{GL}}_{H^1(D)[\frac{1}{p}]}$ that fixes $v_{\alpha}$ for all $\alpha\in\mathcal J$. The flat, affine group scheme $\mathcal D^{\acute et}$ is a closed subgroup scheme of $\mathcal G^{\acute et}$.

\subsection{Fontaine comparison theory} 

We refer to [19], [17], and [64] for the following review of {\it Fontaine comparison theory}. This theory provides us with three rings $B^+_{\text{crys}}(W(k))$, $B_{\text{crys}}(W(k))$, and $B_{\text{dR}}(W(k))$ that are $W(k)$-algebras for which the following six properties hold:

\medskip\noindent
{\bf (i)} {\it The three rings are integral domains equipped with exhaustive and decreasing filtrations and with a Galois action. Moreover $B_{\text{dR}}(W(k))$ is a field.}

\smallskip\noindent
{\bf (ii)} {\it We have $W(k)$-monomorphisms $B^+_{\text{crys}}(W(k))\hookrightarrow B_{\text{crys}}(W(k))\hookrightarrow B_{\text{dR}}(W(k))$.}
 
\smallskip\noindent
{\bf (iii)} {\it The ring $B^+_{\text{crys}}(W(k))$ is faithfully flat over $W(k)$ and has a natural Frobenius lift that is compatible with $\sigma$ and that also extends to an endomorphism of $B_{\text{crys}}(W(k))$.}

\smallskip\noindent
{\bf (iv)} {\it There exists a functorial $B^+_{\text{crys}}(W(k))$-linear monomorphism
$$i_D^+:M\otimes_{W(k)} B^+_{\text{crys}}(W(k))\hookrightarrow H^1(D)\otimes_{\mathbb Z_p} B^+_{\text{crys}}(W(k))$$ that respects the tensor product filtrations, the Galois actions, and the tensor product Frobenius endomorphisms, with the Frobenius endomorphism of $H^1(D)$ being $1_{H^1(D)}$.}

\smallskip\noindent
{\bf (v)} {\it The functorial $B_{\text{dR}}(W(k))$-linear map $i_D:=i_D^+\otimes 1_{B_{\text{dR}}(W(k))}$ is a bijection that induces naturally a $B_{\text{dR}}(W(k))$-linear isomorphism denoted in the same way
$$i_D:\mathcal T(M)\otimes_{W(k)} B_{\text{dR}}(W(k))\rightarrow\mathcal T(H^1(D))\otimes_{\mathbb Z_p} B_{\text{dR}}(W(k)).$$}

\noindent
{\bf (vi)} {\it Each \'etale Tate-cycle $v_{\alpha}$ on $D_{B(k)}$ defines a tensor $t_{\alpha}:=i_D^{-1}(v_{\alpha})\in \mathcal T(M)\otimes_{W(k)} B_{\text{dR}}(W(k))$ which in fact belongs to $F^0(\mathcal T(M))[\frac{1}{p}]\subset\mathcal T(M)[\frac{1}{p}]$ and is fixed by $\phi$.}

\medskip
Let $\mathcal G_{B(k)}$ be the subgroup of $\pmb{\text{GL}}_{M[\frac{1}{p}]}$ that fixes $t_{\alpha}$ for all $\alpha\in\mathcal J$. As $\phi$ fixes each $t_{\alpha}$ we have $\phi(\text{Lie}(\mathcal G_{B(k)}))=\text{Lie}(\mathcal G_{B(k)})$. As we also have $\mathcal G_{\mathbb Q_p}^{\acute et}\times_{\Spec(\mathbb Q_p)} \Spec(B_{\text{dR}}(W(k)))=i_D(\mathcal G_{B(k)}\times_{\Spec(B(k))} \Spec(B_{\text{dR}}(W(k))))i_D^{-1}$, the groups $\mathcal G_{\mathbb Q_p}^{\acute et}\times_{\Spec(\mathbb Q_p)} \Spec(B(k))$ and $\mathcal G_{B(k)}$ are forms of each other. 

Let $\mathcal G$ be the schematic closure of $\mathcal G_{B(k)}$ in $\pmb{\text{GL}}_M$. It is a flat, closed subgroup scheme of $\pmb{\text{GL}}_M$. Let $\mu:\mathbb G_{m,W(k)}\rightarrow \mathcal G$ be a cocharacter that produces a direct sum decomposition $M=F^1\oplus F^0$ such that for each $i\in\{0,1\}$, every $\beta\in\mathbb G_{m,W(k)}(W(k))$ acts through $\mu$ on $F^i$ as the multiplication with $\beta^{-i}$. For instance, we can take $\mu$ to be the factorization through $\mathcal G$ of the inverse of the canonical split cocharacter $\mu_{\text{can}}:\mathbb G_{m,W(k)}\rightarrow \pmb{\text{GL}}_M$ of $(M,F^1,\phi)$ defined in [69, p. 512]; this is so as from the functorial properties in [69, p. 513] we get that $\mu_{\text{can}}$ fixes each $t_{\alpha}$. 

We identify $\text{Hom}(F^1,F^0)$ with the direct summand $\{e\in\text{End}(M)|e(F^0)=0,\;e(F^1)\subset F^0\}$ of $\text{End}(M)$. Let $U_{\text{bigg}}$ and $U$ be the smooth, unipotent, closed subgroup schemes of $\pmb{\text{GL}}_M$ and $\mathcal G$ (respectively) defined by the rule: if $\diamond$ is an arbitrary commutative $W(k)$-algebra, then $U_{\text{bigg}}(\diamond):=1_{M\otimes_{W(k)} \diamond}+\text{Hom}(F^1,F^0)\otimes_{W(k)} \diamond$ and 
$$U(\diamond):=1_{M\otimes_{W(k)} \diamond}+(\text{Lie}(\mathcal G_{B(k)})\cap \text{Hom}(F^1,F^0))\otimes_{W(k)} \diamond.$$
We have $\text{Lie}(U_{\text{bigg}})=\text{Hom}(F^1,F^0)$ and $\text{Lie}(U)=(\text{Lie}(\mathcal G_{B(k)})\cap \text{Hom}(F^1,F^0))$. 

Let $\mathcal H$ be the centralizer of $\mu$ in $\mathcal G$ and let $\mathcal P_0$ be the semidirect product of $\mathcal H$ and $U$. If $\mathcal G$ is smooth, then $\mathcal H$ is a smooth closed subgroup scheme of $\mathcal G$ (cf. [10, Lem. 2.1.5 and Prop. 2.1.8 (3)]) and thus also $\mathcal P_0$ is a a smooth closed subgroup scheme of $\mathcal G$.

\begin{lemma}\label{L12} 
{\it Let $\mu_1:\mathbb G_{m,W(k)}\rightarrow \mathcal G$ be a cocharacter such that we have a direct sum decomposition $M=F^1_1\oplus F^0_1$ with the properties that $\mathbb G_{m,W(k)}$ acts through $\mu_1$ on each $F^i_1$ via the weight $-i$ and we have $F^1_1/pF^1_1=F^1/pF^1$. Then there exists a unique element $v\in p\text{Lie}(U)$ such that for $u:=1_M+v\in\text{Ker}(U(W(k))\rightarrow U(k))$ we have $u(F^1)=F^1_1$.}
\end{lemma}

\medskip\proof
There exists a unique element $u\in \text{Ker}(U_{\text{bigg}}(W(k))\rightarrow U_{\text{bigg}}(k))$ such that we have an identity $u(F^1)=F^1_1$. We write $u=1_M+v$, where $v\in p\text{Hom}(F^1,F^0)=p\text{Lie}(U_{\text{bigg}})$. Let $\mathcal T(M)=\oplus_{i\in\mathbb Z} \tilde F^i(\mathcal T(M))$ be the direct sum decomposition such that $\mathbb G_{m,W(k)}$ acts on each $\tilde F^i(\mathcal T(M))$ through $\mu$ via the weight $-i$. The filtration $(F^i(\mathcal T(M)))_{i\in\mathbb Z}$ of $\mathcal T(M)$ defined by $F^1$ satisfies  for all $i\in\mathbb Z$ the following identity $F^i(\mathcal T(M))=\oplus_{j\ge i} \tilde F^j(\mathcal T(M))$. As $\mu$ and $\mu_1$ are two cocharacters of $\mathcal G$, they fix each $t_{\alpha}$. In particular, we have $t_{\alpha}\in \tilde F^0(\mathcal T(M))[\frac{1}{p}]$ and the tensor $u^{-1}(t_{\alpha})=(1_M-v)(t_{\alpha})$ belongs to $F^0(\mathcal T(M))[\frac{1}{p}]$. As $v\in \text{Hom}(F^1,F^0)\subset\tilde F^{-1}(\mathcal T(M))$, the component of $(1_M-v)(t_{\alpha})$ in $\tilde F^{-1}(\mathcal T(M))[\frac{1}{p}]$ is $-v(t_{\alpha})$ as well as $0$.  Thus $v$ annihilates $t_{\alpha}$ for all $\alpha\in\mathcal J$ and therefore $v\in p\text{Hom}(F^1,F^0)\cap\text{Lie}(\mathcal G_{B(k)})=p\text{Lie}(U)$.\endproof 

\begin{lemma}\label{L13}  
{We assume that the group scheme $\mathcal G$ is smooth. Then $[1_M+\frac{1}{p}\text{Lie}(U)]/U(W(k))$ is the intersection of $[1_M+\frac{1}{p}\text{Lie}(U_{\text{bigg}})]/U_{\text{bigg}}(W(k))$ and $\mathcal P_0(B(k))/\mathcal P_0(W(k))$ taken inside $\pmb{\text{GL}}_M(B(k))/\pmb{\text{GL}}_M(W(k))$. Thus, if $\mathcal G$ is a reductive group scheme over $W(k)$, then $[1_M+\frac{1}{p}\text{Lie}(U)]/U(W(k))$ is the intersection of $[1_M+\frac{1}{p}\text{Lie}(U_{\text{bigg}})]/U_{\text{bigg}}(W(k))$ and $\mathcal G(B(k))/\mathcal G(W(k))$ taken inside $\pmb{\text{GL}}_M(B(k))/\pmb{\text{GL}}_M(W(k))$.}
\end{lemma}
\begin{proof}
We check that if $\mathcal G$ is a reductive group scheme over $W(k)$, then the natural injective map 
$$
\mathcal P_0(B(k))/\mathcal P_0(W(k))\rightarrow \mathcal G(B(k))/\mathcal G(W(k))
$$ 
is a bijection. This is equivalent to the equality $\mathcal G(B(k))=\mathcal P_0(B(k))\mathcal G(W(k))$.  If $\mathcal G$ is a reductive group scheme over $W(k)$, $\mathcal P_0$ is a parabolic subgroup scheme (see [10, Lem. 2.1.5 and Prop. 2.1.8 (3)]) and thus the last equality follows from the fact that projective $W(k)$-scheme $\mathcal G/\mathcal P_0$ has the same sets of $B(k)$- and $W(k)$-valued points.

We are left to show that if $c\in \frac{1}{p}\text{Lie}(U_{\text{bigg}})$ and $g\in \mathcal P_0(B(k))$ are such that $g(M)=(1_M+c)(M)$, then the reduction $\bar X\in\text{Lie}(U_{\text{big},k})$ of $X:=pc$ modulo $p$ is in fact an element of $\text{Lie}(U_k)$. 

We consider the smooth, closed subgroup schemes $\mathcal G_1:=g\mathcal G g^{-1}$ of $\pmb{\text{GL}}_{g(M)}=\pmb{\text{GL}}_{(1_M+c)(M)}$ and $\tilde{\mathcal G}:=(1_M-c)\mathcal G_1(1_M+c)$ of $\pmb{\text{GL}}_M$. Both $U_{\text{bigg}}$ and $U$ are closed subgroup schemes of $\pmb{\text{GL}}_{(1_M+c)(M)}$ and moreover $U\leqslant\mathcal G_1$. 

As $\mathcal G$ is smooth through which $\mu$ factors, we have $U_{\text{bigg}}\cap \mathcal G=U$ (cf. [10, Lem. 2.1.5 and Prop. 2.1.8 (3)]). Similarly, as $\mathcal G_1$ is smooth through which $g\mu g^{-1}$ factors, we have $U_{\text{bigg}}\cap \mathcal G_1=U$.

All $2\times 2$ block matrices with coefficients in $W(k)$ or $k$ will be with respect to the direct sum decomposition $M=F^1\oplus F^0$ or its  reduction modulo $p$. For each $t\in W(k)$, the element 
$\binom{1+pt\;0}{0\;\;\;\;\;\;1}=\mu(W(k))((1+pt)^{-1})$ belongs to $\mathcal G_1(W(k))$. Thus $\binom{1\;\;\;0}{-c\;1}\binom{1+pt\;0}{0\;\;\;\;\;1}\binom{1\;\;0}{c\;\;1}=\binom{1+pt\;0}{-tX\;\;1}$ belongs to $\tilde{\mathcal G}(W(k))$. Therefore $\binom{1\;\;\;0}{\bar t\bar X\;1}$ belongs to $\tilde{\mathcal G}(k)$ for all $\bar t\in k$. Conjugating with $1_M+c$ we get that $\binom{1\;\;\;0}{\bar t\bar X\;1}$ belongs to $(U_{\text{bigg}}\cap \mathcal G_1)(k)=U(k)$ for all $\bar t\in k$. Thus $\bar X\in\Lie(U_k)$.\end{proof}

\begin{theorem} \label{T9} ([64, Thm. 1.2 and Ex. 4.4.1]) 
{If $p=2$, then we assume that $D$ is a direct sum of connected and \'etale Barsotti--Tate groups (e.g., this holds if $\mathcal G^{\acute et}$ is a torus). Then there exist isomorphisms 
$$
\varrho_D:(M,(t_{\alpha})_{\alpha\in\mathcal J})\rightarrow (H^1(D)\otimes_{\mathbb Z_p} W(k),(v_{\alpha})_{\alpha\in\mathcal J})
$$ 
(in the sense of Subsection 2.1). If moreover $D$ has a principal quasi-polarization, then  there exist isomorphisms $\varrho_D:(M,(t_{\alpha})_{\alpha\in\mathcal J},\psi_M)\rightarrow (H^1(D)\otimes_{\mathbb Z_p} W(k),(v_{\alpha})_{\alpha\in\mathcal J},\psi_{H^1(D)})$.}
\end{theorem}

If $\mathcal G^{\acute et}$ is a reductive group scheme over $\mathbb Z_p$ and if for $p=2$ the $2$-divisible group $D$ is connected, then Theorem \ref{T9} is also proved in [27, Cor. 1.4.3].

\begin{lemma}\label{L14} 
Let $k_1$ be an algebraically closed field that contains $k$. We assume that there exists an isomorphism $(M\otimes_{W(k)} W(k_1),(t_{\alpha})_{\alpha\in\mathcal J})\rightarrow (H^1(D)\otimes_{\mathbb Z_p} W(k_1),(v_{\alpha})_{\alpha\in\mathcal J})$. Then there exists an isomorphism $\varrho_D:(M,(t_{\alpha})_{\alpha\in\mathcal J})\rightarrow (H^1(D)\otimes_{\mathbb Z_{p}} W(k),(v_{\alpha})_{\alpha\in\mathcal J})$.\end{lemma}

\medskip\proof 
To check the existence of $\varrho_D$ we can assume that we have $t_{\alpha}\in\mathcal T(M)$ and $v_{\alpha}\in H^1(D)$ for all $\alpha\in\mathcal J$. Thus we an speak about the affine $W(k)$-scheme $\mathfrak{P}$ of finite type that parameterizes isomorphisms between $(M,(t_{\alpha})_{\alpha\in\mathcal J})$ and $(H^1(D)\otimes_{\mathbb Z_{p}} W(k),(v_{\alpha})_{\alpha\in\mathcal J})$. We know that $\mathfrak{P}$ has a $W(k_1)$-valued point. As the monomorphism $W(k)\hookrightarrow W(k_1)$ is of ramification index one, from [7, Ch. 3, Sect. 3.6, Prop. 4] we get that there exists a morphism $\mathfrak{P}^\prime\rightarrow\mathfrak{P}$ of $W(k)$-schemes such that $\mathfrak{P}^\prime$ is smooth over $W(k)$ and has a $W(k_1)$-valued point. Thus the special fibre $\mathfrak{P}^\prime_k$ is non-empty. As $\mathfrak{P}^\prime$ is smooth over $W(k)$ and has a non-empty special fibre, it has $W(k)$-valued points. Therefore $\mathfrak{P}$ also has $W(k)$-valued points and thus the isomorphism $\varrho_D$ exists.\endproof

\subsection{Group correspondences} 

Let $\mathcal F^{\acute et}_{\mathbb Q_p}$ be a reductive, closed subgroup of $\mathcal G^{\acute et}_{\mathbb Q_p}$. The restriction to $\text{Lie}(\mathcal F^{\acute et}_{\mathbb Q_p})$ of the trace bilinear form on $\text{End}(H^1(D)[\frac{1}{p}])$ is non-degenerate, cf. Lemma \ref{L11} (b). Let $\text{Lie}(\mathcal F^{\acute et}_{\mathbb Q_p})^\perp$ be the perpendicular on $\text{Lie}(\mathcal F^{\acute et}_{\mathbb Q_p})$ with respect to the trace bilinear form on $\text{End}(H^1(D)[\frac{1}{p}])$; we have a direct sum decomposition of $\mathbb Q_p$-vector spaces
$$\text{End}(H^1(D)[\frac{1}{p}])=\text{Lie}(\mathcal F^{\acute et}_{\mathbb Q_p})\oplus\text{Lie}(\mathcal F^{\acute et}_{\mathbb Q_p})^\perp.$$ 
Let $\pi^{\acute et}$ be the projector of $\text{End}(H^1(D)[\frac{1}{p}])$ on $\text{Lie}(\mathcal F^{\acute et}_{\mathbb Q_p})$ along $\text{Lie}(\mathcal F^{\acute et}_{\mathbb Q_p})^{\perp}$; it is an idempotent of $\text{End}(H^1(D)[\frac{1}{p}])$ fixed by  each subgroup of $\pmb{\text{GL}}_{H^1(D)[\frac{1}{p}]}$ that normalizes $\mathcal F^{\acute et}_{\mathbb Q_p}$. 

If $\mathcal D^{\acute et}_{\mathbb Q_p}$ normalizes $\mathcal F^{\acute et}_{\mathbb Q_p}$ (e.g., this holds if $\mathcal F^{\acute et}_{\mathbb Q_p}$ is a normal subgroup of $\mathcal G^{\acute et}_{\mathbb Q_p}$), then $\pi^{\acute et}$ is fixed by $\mathcal D^{\acute et}_{\mathbb Q_p}$ and thus also by $\text{Im}(\rho_D)$ and therefore we can speak about the projector $\pi^{\text{crys}}$ of $\text{End}(M[\frac{1}{p}])$ that corresponds to $\pi^{\acute et}$ via Fontaine comparison theory.

\begin{lemma}\label{L15} 
We assume that $\mathcal D^{\acute et}_{\mathbb Q_p}$ normalizes $\mathcal F^{\acute et}_{\mathbb Q_p}$. Then the following two properties hold:

\medskip
{\bf (a)} There exists a unique reductive subgroup $\mathcal F_{B(k)}$ of $\mathcal G_{B(k)}$ whose Lie algebra is $\text{Im}(\pi^{\text{crys}})$.

\smallskip
{\bf (b)} If the generic fibre of $\mu_{\text{can}}$ factors through $\mathcal F_{B(k)}$, then $\mathcal D^{\acute et}_{\mathbb Q_p}$ is a subgroup of $\mathcal F^{\acute et}_{\mathbb Q_p}$.
\end{lemma}

\medskip\proof
We check part (a). As $i_D^{-1}$ is a $B_{\text{dR}}(W(k))$-linear isomorphism that takes $\pi^{\acute et}$ to $\pi^{\text{crys}}$, the group $i_D^{-1}(\mathcal F^{\acute et}_{\mathbb Q_p}\times_{\Spec(\mathbb Q_p)} \Spec(B_{\text{dR}}(W(k))))i_D$ is a subgroup of
$$i_D^{-1}(\mathcal G^{\acute et}_{\mathbb Q_p}\times_{\Spec(\mathbb Q_p)} \Spec(B_{\text{dR}}(W(k)))i_D=\mathcal G_{B(k)}\times_{\Spec(B(k))} \Spec(B_{\text{dR}}(W(k)))$$
whose Lie algebra is $\text{Im}(\pi^{\text{crys}})\otimes_{B(k)} B_{\text{dR}}(W(k))$. Thus as $B_{\text{dR}}(W(k))$ is a field, from Lemma \ref{L11} (a) applied with $(\mathcal W,\mathcal L,\eta,\eta_1)=(M[\frac{1}{p}],\text{Im}(\pi^{\text{crys}}),B(k),B_{\text{dR}}(W(k)))$, we get that there exists a unique reductive subgroup $\mathcal F_{B(k)}$ of $\pmb{\text{GL}}_{M[\frac{1}{p}]}$ whose Lie algebra is $\text{Im}(\pi^{\text{crys}})$. As $\mathcal F_{B(k)}\times_{\Spec(B(k))} \Spec(B_{\text{dR}}(W(k)))$ is a subgroup of $\mathcal G_{B(k)}\times_{\Spec(B(k))} \Spec(B_{\text{dR}}(W(k)))$, the group $\mathcal F_{B(k)}$ is in fact a subgroup of $\mathcal G_{B(k)}$. Thus part (a) holds.

We check part (b). Let $l_{\text{can}}$ be the Lie algebra of the image of the generic fibre of $\mu_{\text{can}}$. As $\pi^{\text{crys}}$ is fixed by $\phi$, the Lie algebra $\text{Lie}(\mathcal F_{B(k)})=\text{Im}(\pi^{\text{crys}})$ is normalized by $\phi$. Let $\mathcal D_{B(k)}$ be the smallest connected subgroup of $\mathcal F_{B(k)}$ with the property that its Lie algebra $\text{Lie}(\mathcal D_{B(k)})$ contains $\phi^m(l_{\text{can}})$ for all $m\in\mathbb Z$. From [5, Ch. I, Sect. 7.1] we get that all conjugates of the generic fibre of $\mu_{\text{can}}$ through integral powers of $\phi$ factor through $\mathcal D_{B(k)}$ and $\mathcal D_{B(k)}$ is the smallest subgroup of $\mathcal F_{B(k)}$ that has this property. This implies that $\mathcal D_{B(k)}$ corresponds to $\mathcal D_{\mathbb Q_p}^{\acute et}$ via Fontaine comparison theory (cf. [69, Prop. 4.2.3]), i.e., we have an identity of subgroups
$$\mathcal D_{\mathbb Q_p}^{\acute et}\times_{\Spec(\mathbb Q_p)} \Spec(B_{\text{dR}}(W(k)))=i_D(\mathcal D_{B(k)}\times_{\Spec(B(k))} \Spec(B_{\text{dR}}(W(k))))i_D^{-1}$$
of $\pmb{\text{GL}}_{H^1(D)\otimes_{\mathbb Z_p} B_{\text{dR}}(W(k))}$. Thus, as $\mathcal D_{B(k)}$ is a subgroup of $\mathcal F_{B(k)}$, $\mathcal D_{\mathbb Q_p}^{\acute et}\times_{\Spec(\mathbb Q_p)} \Spec(B_{\text{dR}}(W(k)))$ is a subgroup of $\mathcal F_{\mathbb Q_p}^{\acute et}\times_{\Spec(\mathbb Q_p)} \Spec(B_{\text{dR}}(W(k)))=i_D(\mathcal F_{B(k)}\times_{\Spec(B(k))} \Spec(B_{\text{dR}}(W(k))))i_D^{-1}$. From this part (b) follows.\endproof

\subsection{Faltings deformation theory} 

\smallskip
Let $l\in\mathbb N$. Let $R=W(k)[[x_1,\ldots,x_l]]$ be the ring of formal power series in $l$ variables with coefficients in $W(k)$. Let $\Phi_R$ be the Frobenius lift of $R$ that is compatible with $\sigma$ and that takes $x_i$ to $x_i^p$ for all $i\in\{1,\ldots,l\}$. We consider the ideal $\mathfrak{I}:=(x_1,\ldots,x_l)$ of $R$. Let $\hat\Omega_{R/W(k)}=\oplus_{i=1}^l Rdx_i$ be the $\mathfrak{I}$-adic completion of the $R$-module of relative differentials $\Omega_{R/W(k)}$. Let $d\Phi_R:\hat\Omega_{R/W(k)}\rightarrow\hat\Omega_{R/W(k)}$ be the ($\mathfrak{I}$-adic completion of the) differential map of $\Phi_R$.

Let $(M_R,F^1_R,\Phi)$ be a triple such that the following four axioms hold:

\medskip\noindent
{\bf (i)} $M_R$ is a free $R$-module of rank equal to the height of $D$;

\smallskip\noindent
{\bf (ii)} $F^1_R$ is a direct summand of $M_R$ of rank equal to the rank of $F^1$;

\smallskip\noindent
{\bf (iii)} $\Phi:M_R\rightarrow M_R$ is a $\Phi_R$-linear endomorphism that induces an $R$-linear isomorphism $(M_R+\frac{1}{p}F^1_R)\otimes_R {}_{\Phi_R} R\rightarrow M_R$;

\smallskip\noindent
{\bf (iv)} the reduction of $(M_R,F^1_R,\Phi)$ modulo $\mathfrak{I}$ is canonically identified with $(M,F^1,\phi)$.

\medskip
Let $\Phi$ act in the natural tensor way on $\mathcal T(M_R)[\frac{1}{p}]$. For instance, if $e\in M_R^{\vee}:=\text{Hom}(M_R,R)$, then $\Phi(e)\in M_R^{\vee}[\frac{1}{p}]$ is the unique element such that we have $\Phi(e)(\Phi(a))=\Phi_R(e(a))\in R$ for all $a\in M_R$. 

It is known that there exists a unique connection $\nabla:M_R\rightarrow M_R\otimes_R \hat\Omega_{R/W(k)}$ such that we have an identity $\nabla\circ\Phi=(\Phi\otimes d\Phi_R)\circ\nabla$ and moreover such a connection is automatically integrable and nilpotent modulo $p$, cf. either [17, Thm. 10] or [64, Thm. 3.2 and Cor. 3.3.2]. By viewing $\mathcal T(M_R)[\frac{1}{p}]$ as a module over the Lie algebra (associated to) $\text{End}(M_R)$, we can view also $\nabla$ as a connection on the $R$-module  $\mathcal T(M_R)[\frac{1}{p}]$ and thus it makes sense to say that it annihilates some specific tensor of $\mathcal T(M_R)[\frac{1}{p}]$.

\begin{lemma}\label{L16} 
There exists a unique Barsotti--Tate group $D_R$ over $R$ which modulo the ideal $\mathfrak{I}$ is $D$ and such that its filtered $F$-crystal over $R/pR$ is $(M_R,F^1_R,\Phi,\nabla)$.
\end{lemma}

\medskip\proof
Let $J$ be an ideal of $R$ such that $R$ is complete in the $J$-adic topology (e.g., $(p)$, $\mathfrak{I}$, or $p\mathfrak{I}$). Let $\Spf(R)$ be the formal scheme which is the formal completion of $\Spec(R)$ along $\Spec(R/J)$. The categories of Barsotti--Tate groups over $\Spec(R)$ and respectively over $\text{Spf}(R)$ are canonically isomorphic, cf. [35, Ch. II, Lem. 4.16]; below we will use this fact without any extra comment. 

The existence of $D_R$ is implied by [17, Thm. 10]. The uniqueness of the fibre $D_{R/pR}$ of $D_R$ over $\Spec(R/pR)$ is implied by [3, Thm. 4.1.1]. As the ideal $p(\mathfrak{I}/\mathfrak{I}^m)$ of $R/\mathfrak{I}^m$ has a natural nilpotent divided power structure for all $m\in\mathbb N^*$, from the Grothendieck--Messing deformation theory we get that $D_R$ is the unique Barsotti--Tate group over $R$ that lifts both $D$ and $D_{R/pR}$ and whose filtered $F$-crystal is $(M_R,F^1_R,\Phi,\nabla)$.\endproof 

\begin{lemma}\label{L17} 
We assume that $D$ has a principal quasi-polarization $\lambda_D$. We also assume that there exists a perfect, alternating bilinear form $\psi_{M_R}$ on $M_R$ that lifts $\psi_M$ (i.e., which modulo $\mathfrak{I}$ is $\psi_M$), that satisfies $\psi_{M_R}(F^1_R,F^1_R)=0$ (i.e., $F^1_R$ is anisotropic with respect to $\psi_{M_R}$), and such that for all $a,b\in M_R$ we have $\psi_{M_R}(\Phi(a),\Phi(b))=p\Phi_R(\psi_{M_R}(a,b))$. Then there exists a unique principal quasi-polarization $\lambda_{D_R}$ of $D_R$ which modulo the ideal $\mathfrak{I}$ is $\lambda_D$ and whose crystalline realization is $\psi_{M_R}$.
\end{lemma}

\medskip\proof
Let $(M_R^{\text{t}},F_R^{1\text{t}},\Phi^{\text{t}},\nabla^{\text{t}})$ be the filtered $F$-crystal over $R/pR$ of the Cartier dual $D_R^{\text{t}}$ of $D_R$. The form $\psi_{M_R}$ defines naturally an isomorphism $\theta_0:(M_R^{\text{t}},F_R^{1\text{t}},\Phi^{\text{t}})\rightarrow (M_R,F^1_R,\Phi)$. As the connections $\nabla$ and $\nabla^{\text{t}}$ are uniquely determined by $(M_R,F^1_R,\Phi)$ and $(M_R^{\text{t}},F_R^{1\text{t}},\Phi^{\text{t}})$ (respectively), $\theta_0$ extends to an isomorphism $\theta:(M_R^{\text{t}},F_R^{1\text{t}},\Phi^{\text{t}},\nabla^{\text{t}})\rightarrow (M_R,F^1_R,\Phi,\nabla)$ of filtered $F$-crystals over $R/pR$. 

As the ring $R/pR$ has a finite $p$-basis $\{x_1,\ldots, x_l\}$ in the sense of [3, Def. 1.1.1],  from the fully faithfulness part of [3, Thm. 4.1.1] we get that there exists a unique principal quasi-polarization $\lambda_{D_{R/pR}}:D_{R/pR}\rightarrow D_{R/pR}^{\text{t}}$ whose crystalline realization is $\theta$; it lifts the special fibre of $\lambda_D$. As the ideal $p(\mathfrak{I}/\mathfrak{I}^m)$ of $R/\mathfrak{I}^m$ has a natural nilpotent divided power structure for all $m\in\mathbb N^*$, from the Grothendieck--Messing deformation theory we get that there exists a unique principal quasi-polarization $\lambda_{D_R}$ of $D_R$ that lifts both $\lambda_{D_{R/pR}}$ and $\lambda_D$ and whose crystalline realization is $\psi_{M_R}$.\endproof 

\subsubsection{Explicit filtered $F$-crystal with tensors} 

Let $M=F^1\oplus F^0$, $U^{\text{bigg}}$, $U$, $\mathcal H$, and $\mathcal P_0$ be as before Lemma \ref{L12}. Let $\mathcal G^\prime$ be the universal smoothening of $\mathcal G$, cf. Subsection A1. 

Until Subsection B5 we will assume that $D$ has a principal quasi-polarization $\lambda_D$, that $\mathcal G$ is a closed subgroup scheme of $\pmb{\text{GSp}}(M,\psi_M)$, and that $R=W(k)[[x_1,\ldots,x_l]]$ is the completion of the local ring of $\mathcal G^\prime$ at the identity element of $\mathcal G^\prime_k$. Thus the relative dimension of $\mathcal G$ over $W(k)$ is $l$. The closed embedding $U\hookrightarrow \mathcal G$ factors through $\mathcal G^\prime$ (cf. Subsection A1); thus $U$ is a closed subgroup scheme of $\mathcal G^\prime$. 

Let $g_{\text{univ}}\in\mathcal G^\prime(R)$ be the universal element. We define $(M_R,F^1_R):=(M,F^1)\otimes_{W(k)} R$ and $\Phi:=g_{\text{univ}}(\phi\otimes\Phi_R)$. Let
$$\mathfrak{C}_{\text{univ}}:=(M_R,F^1_R,\Phi,\nabla,(t_{\alpha})_{\alpha\in\mathcal J}).$$
The $W(k)$-algebra $R$ is complete in the $\mathfrak{I}$-topology and we have $\Phi_R(\mathfrak{I})\subset \mathfrak{I}^p$. This implies that each element of $\text{Ker}(\mathbb G_{m,W(k)}(R)\rightarrow\mathbb G_{m,W(k)}(R/\mathfrak{I}))$ is of the form $\beta\Phi_R(\beta^{-1})$ for some element $\beta\in\text{Ker}(\mathbb G_{m,W(k)}(R)\rightarrow\mathbb G_{m,W(k)}(R/\mathfrak{I}))$. As $g_{\text{univ}}$ takes $\psi_M$ to a $\text{Ker}(\mathbb G_{m,W(k)}(R)\rightarrow\mathbb G_{m,W(k)}(R/\mathfrak{I}))$-multiple of $\psi_M$, we get that there exists a $\text{Ker}(\mathbb G_{m,W(k)}(R)\rightarrow\mathbb G_{m,W(k)}(R/\mathfrak{I}))$-multiple $\psi_{M_R}$ of the perfect, alternating bilinear form $\psi_M$ on $M_R$ such that we have an identity
$$\psi_{M_R}(\Phi(a),\Phi(b))=p\Phi_R(\psi_{M_R}(a,b))$$ 
for all $a, b\in M_R$. As $1$ is the only element of $\text{Ker}(\mathbb G_{m,W(k)}(R)\rightarrow\mathbb G_{m,W(k)}(R/\mathfrak{I}))$ fixed by $\Phi_R$, this $\text{Ker}(\mathbb G_{m,W(k)}(R)\rightarrow\mathbb G_{m,W(k)}(R/\mathfrak{I}))$-multiple $\psi_{M_R}$ of $\psi_M$ is uniquely determined. 

We have the following three properties:

\medskip\noindent
{\bf (i)} {\it The connection on $\mathcal T(M_R)=\mathcal T(M)\otimes_{W(k)} R$ induced naturally by $\nabla$ (and denoted in the same way) annihilates the tensor $t_{\alpha}\in\mathcal T(M)\otimes_{W(k)} R[\frac{1}{p}]$ for all $\alpha\in\mathcal J$.}

\smallskip\noindent
{\bf (ii)} {\it The connection $\nabla$ is of the form $\delta+\gamma$, where $\delta$ is the flat connection on $M_R=M\otimes_{W(k)} R$ that annihilates $M\otimes 1$ and where $\gamma\in(\text{Lie}(\mathcal G_{B(k)})\cap\text{End}(M))\otimes_{W(k)} \hat\Omega_{R/W(k)}$.} 

\smallskip\noindent
{\bf (iii)} {\it The Kodaira--Spencer map of the connection $\nabla$ has an image $\Theta$ which is the direct summand $\text{Lie}(U)\otimes_{W(k)} R$ of $\text{Lie}(U_{\text{bigg}})\otimes_{W(k)} R\rightarrow\text{Hom}(F^1,M/F^1)\otimes_{W(k)} R$.}

\medskip
As $\phi$ fixes $t_{\alpha}$ and $\nabla\circ \Phi=(\Phi\otimes d\Phi_R)\circ\nabla$, we have $\nabla(t_{\alpha})=(\Phi\otimes d\Phi_R)(\nabla(t_{\alpha}))$.  As $d\Phi_R(x_i)=px_i^{p-1}dx_i$, by induction on $n\in\mathbb N^*$ we get that $\nabla(t_{\alpha})\in M\otimes_{W(k)} \mathfrak{I}^n (\oplus_{i=1}^l Rdx_i)[\frac{1}{p}]$. This implies that (i) holds. Property (ii) follows from the property (i) and the fact that $\text{Lie}(\mathcal G_{B(k)})\cap\text{End}(M)$ is the Lie subalgebra of $\text{End}(M)$ which annihilates $t_{\alpha}$ for all $\alpha\in\mathcal J$. 

To check (iii), we first remark that the property (ii) implies that $\Theta$ is contained in the image of  $(\text{Lie}(\mathcal G_{B(k)})\cap\text{End}(M))\otimes_{W(k)} R$ in $\text{Lie}(U_{\text{bigg}})\otimes_{W(k)} R\rightarrow\text{Hom}(F^1,M/F^1)\otimes_{W(k)} R$ and thus it is contained in $\text{Lie}(U)\otimes_{W(k)} R$. It is easy to see that $\gamma$ modulo $(p,\mathfrak{I}^{p-1})$ is $g_{\text{univ}}^{-1}dg_{\text{univ}}$ modulo $(p,\mathfrak{I}^{p-1})$ (for instance, this follows from [64, Equations (11) and (12)]). Thus, as $U$ is a closed subgroup scheme of $\mathcal G^\prime$ and as $g_{\text{univ}}\in\mathcal G^\prime(R)$ is the universal element, we get that $\Theta$ surjects onto $\text{Lie}(U)\otimes_{W(k)} R/(p,\mathfrak{I})$. From this and the inclusion $\Theta\subset \text{Lie}(U)\otimes_{W(k)} R$ we get that the property (iii) holds. 

Let $m\in\mathbb N$, $R_1:=W(k)[[x_1,\ldots,x_m]]$, and $Z:=\Spec(R_1)$. Let $\Phi_{R_1}$ be the Frobenius lift of $R_1$ that is compatible with $\sigma$ and that takes $x_i$ to $x_i^p$ for all $i\in\{1,\ldots,m\}$. We consider the ideal $\mathfrak{I}_1:=(x_1,\ldots,x_m)$ of $R_1$. 

Let $(M_1,F_1^1,\Phi_1,\nabla_1)$ be a filtered $F$-crystal over $R_1/pR_1$. Thus:

\medskip\noindent
{\bf (iv)} {\it $\Phi_1$ induces an $R_1$-linear isomorphism $(M_1+\frac{1}{p}F_1^1)\otimes_{R_1} {}_{\Phi_{R_1}} R_1\rightarrow M_1$.}

\medskip 
Let $\mathfrak{C}_1:=(M_1,F_1^1,\Phi_1,\nabla_1,(t_{1,\alpha})_{\alpha\in\mathcal J})$, where $(t_{1,\alpha})_{\alpha\in\mathcal J}$ is a family of tensors $(t_{1,\alpha})_{\alpha\in\mathcal J}$ of $\mathcal T(M_1)[\frac{1}{p}]$ such that the following two axioms hold (here $\mathcal T(M_1)$ is as in Subsection 2.1):

\medskip\noindent
{\bf (v)} {\it Each tensor $t_{1,\alpha}$ is fixed by $\Phi_1$, is annihilated by $\nabla_1$, and belongs to $F^0(\mathcal T(M_1))[\frac{1}{p}]$ (here $(F^i(\mathcal T(M_1)))_{i\in\mathbb Z}$ is the filtration of $\mathcal T(M_1)$ defined by $F^1_1$, cf. Subsection 2.1).}

\smallskip\noindent
{\bf (vi)} {\it Its reduction modulo the ideal $\mathfrak{I}_1$ is $(M,F^1,\phi,(t_{\alpha})_{\alpha\in\mathcal J})$.}

\medskip
The $R_1$-module $M_1$ is free of rank equal to the rank of $M$, cf. property (vi). We consider the closed embedding $z_Z:\Spec(W(k))\hookrightarrow Z$ defined by the ideal $\mathfrak{I}_1$ of $R_1$.

\begin{theorem} \label{T10}
The following two properties hold:

\medskip
{\bf (a)} There exists a morphism $i_Z:Z\rightarrow\Spec(R)$ of $W(k)$-schemes such that $g_{\text{univ}}\circ i_Z\circ z_Z$ is the identity section of $\mathcal G^\prime$ and $\mathfrak{C}_1$ is isomorphic to $i_Z^*(\mathfrak{C}_{\text{univ}})$ under an isomorphism which modulo the ideal $\mathfrak{I}_1$ becomes the identity automorphism $1_M$ of $M$.

\smallskip
{\bf (b)} If there exists a perfect, alternating bilinear form $\psi_{M_1}$ on $M_1$ which modulo $\mathfrak{I}_1$ is $\psi_M$ and which is a principal quasi-polarization of the filtered $F$-crystal $(M_1,F^1,\Phi_1,\nabla_1)$ over $R_1/pR_1$, then $(\mathfrak{C}_1,\psi_1)$ is isomorphic to $i_Z^*(\mathfrak{C}_{\text{univ}},\psi_{M_R})$ under an isomorphism which modulo the ideal $\mathfrak{I}_1$ becomes the identity automorphism $1_M$ of $M$.
\end{theorem}

\medskip\proof
If $\mathcal G$ is smooth, then part (a) is a particular case of [17, Thm. 10 and Rm. iii) after it]. To prove part (a) in the general case, we follow the proof of [64, Thm. 5.3]. Let $D_{R_1}$ be the unique Barsotti--Tate group over $R_1$ which modulo the ideal $\mathfrak{I}_1$ is $D$ and whose filtered $F$-crystal over $R_1/pR_1$ is $(M_1,F_1^1,\Phi_1,\nabla_1)$, cf. Lemma \ref{L17}. 

By induction on $s\in\mathbb N^*$ we show that there exists a morphism $i_{Z,s}:\Spec(R_1/\mathfrak{I}_1^s)\rightarrow\Spec(R)$ of $W(k)$-schemes which at the level of rings maps $\mathfrak{I}$ to $\mathfrak{I}_1/\mathfrak{I}_1^s$ and such that $i_{Z,s}^*(D_R)$ is isomorphic to $D_{R_1}$ modulo $\mathfrak{I}_1^s$ under a unique isomorphism $\mathcal I_s$ that has the following two properties:

\medskip\noindent
{\bf (i)} {\it it lifts the identity automorphism of $D$;}

\smallskip\noindent
{\bf (ii)} {\it its crystalline Dieudonn\'e realization defines an isomorphism $\mathcal E_s$ between $\mathfrak{C}_1$ modulo $\mathfrak{I}_1^s$ and $i_{Z,s}^*(\mathfrak{C}_{\text{univ}})$ which modulo $\mathfrak{I}_1/\mathfrak{I}_1^s$ is the identity automorphism $1_M$ of $M$.} 

\medskip
As $\Phi_{R_1}(\mathfrak{I}_1)\subset \mathfrak{I}_1^p$ and the ideal $\mathfrak{I}_1/\mathfrak{I}_1^s$ is complete, such an isomorphism $\mathcal E_s$ is unique. We take $i_{Z,1}$ to be defined by the $W(k)$-epimorphism $R\twoheadrightarrow R/\mathfrak{I}=W(k)=R_1/\mathfrak{I}_1$ and we take $\mathcal I_1$ and $\mathcal E_1$ to be defined by the identity automorphism of $D$ and by $1_M$ (respectively). Thus the existence and the uniqueness of $i_{Z,1}$ and $\mathcal I_1$ are obvious. 

For $s\ge 2$ the passage from $s-1$ to $s$ goes as follows. We endow the ideal $\mathfrak{J}_s:=\mathfrak{I}_1^{s-1}/\mathfrak{I}_1^s$ of $R_1/\mathfrak{I}_1^s$ with the trivial divided power structure; thus $\mathfrak{J}_s^{[2]}=0$. The uniqueness of $\mathcal I_s$ is implied by the uniqueness of $\mathcal I_{s-1}$ and $\mathcal E_{s}$, cf. Grothendieck--Messing deformation theory. To end the induction, we check that we can choose $i_{Z,s}$ such that $\mathcal I_s$ and $\mathcal E_s$ exist. 

Let $\tilde i_{Z,s}:\Spec(R_1/\mathfrak{I}_1^s)\rightarrow\Spec(R)$ be an arbitrary morphism of $W(k)$-schemes through which $i_{Z,s-1}$ factors naturally. We write 
$$\tilde i_{Z,s}^*(\mathfrak{C}_{\text{univ}})=(M\otimes_{W(k)} R_1/\mathfrak{I}_1^s,F^1\otimes_{W(k)} R_1/\mathfrak{I}_1^s,{}_s\Phi,{}_s\nabla,(t_{\alpha})_{\alpha\in\mathcal J}).$$ 
Due to the existence of $\mathcal I_{s-1}$, there exists (cf. Grothendieck--Messing deformation theory) a direct summand ${}_sF^1$ of $M\otimes_{W(k)} R_1/\mathfrak{I}_1^s$ that lifts $F^1\otimes_{W(k)} R_1/\mathfrak{I}_1^{s-1}$ and such that the quintuple $(M_1,F_1,\Phi_1,\nabla_1)$ modulo $\mathfrak{I}_1^s$ is isomorphic to the quintuple 
$(M\otimes_{W(k)} R_1/\mathfrak{I}_1^s,{}_sF^1,{}_s\Phi,{}_s\nabla)$ under an isomorphism $\tilde{\mathcal E}_s$ that lifts the one defined by $\mathcal E_{s-1}$. Let $t_{1,\alpha,s}\in\mathcal T(M\otimes_{W(k)} R_1/\mathfrak{I}_1^s)$ be the image under $\tilde{\mathcal E}_s$ of $t_{1,\alpha}$.  As $t_{1,\alpha}$ is fixed by $\Phi_1$, $t_{1,\alpha,s}$ is fixed by ${}_s\Phi$. As $\tilde{\mathcal E}_s$ lifts $\mathcal E_{s-1}$, the reductions modulo $\mathfrak{J}_s$ of $t_{\alpha}$ and $t_{1,\alpha,s}$ coincide. As ${}_s\Phi(\mathcal T(M)\otimes_{W(k)} \mathfrak{J}_s)=0$, inside $\mathcal T(M)\otimes_{W(k)} R_1/\mathfrak{I}_1^s$ we have
$$t_{1,\alpha,s}-t_{\alpha}={}_s\Phi(t_{1,\alpha,s}-t_{\alpha})\in {}_s\Phi(\mathcal T(M)\otimes_{W(k)} \mathfrak{J}_s)=0.$$
Thus we have $t_{1,\alpha,s}=t_{\alpha}\in\mathcal T(M)\otimes_{W(k)} R_1/\mathfrak{I}_1^s$ for all $\alpha\in\mathcal J$. 

Let $v_s\in \text{Lie}(U_{\text{bigg}})\otimes_{W(k)} \mathfrak{J}_s$ be the unique element such that we have 
$$(1_{M\otimes_{W(k)} R_1/\mathfrak{I}_1^s}+v_s)(F^1\otimes_{W(k)} R_1/\mathfrak{I}_1^s)={}_sF^1.$$ 
As each $t_{1,\alpha,s}=t_{\alpha}$ belongs to the $F^0$-filtrations defined by either ${}_sF^1$ or $F^1\otimes_{W(k)} R_1/\mathfrak{I}_1^s$ and as the $W(k)$-module $\mathfrak{I}_s$ is torsionless, as in [64, proof of Thm. 5.3, bottom of p. 241 and top of p. 242] we argue that $v_s\in\text{Lie}(U)\otimes_{W(k)} \mathfrak{J}_s$. Based on this and the property (iii), as in [64, proof of Thm. 5.3, p. 242] we argue that we can replace $\tilde i_{Z,s}$ by another morphism $i_{Z,s}:\Spec(R_1/\mathfrak{I}_1^s)\rightarrow\Spec(R)$ through which $i_{Z,s-1}$ factors and for which ${}_sF^1$ gets replaced by (i.e., becomes) $F^1\otimes_{W(k)} R_1/\mathfrak{I}_1^s$. From Grothendieck--Messing deformation theory we get that $i_{Z,s}^*(D_R)$ is isomorphic to $D_{R_1}$ modulo $\mathfrak{I}_1^s$ under an isomorphism $\mathcal I_s$ which lifts $\mathcal I_{s-1}$ and which defines an isomorphism $\mathcal E_s$ between $\mathfrak{C}_1$ modulo $\mathfrak{I}_1^s$ and $i_{Z,s}^*(\mathfrak{C}_{\text{univ}})$. As $\mathcal I_s$ lifts $\mathcal I_{s-1}$, the uniqueness of $\mathcal E_{s-1}$ implies that $\mathcal E_s$ lifts $\mathcal E_{s-1}$. This ends the induction. 

We take $i_Z:Z\rightarrow\Spec(R)$ such that it lifts $i_{Z,s}$ for all $s\in\mathbb N^*$. From the very definition of $i_{Z,1}$ we get that $g_{\text{univ}}\circ i_Z\circ z_Z$ is the identity section of $\mathcal G^\prime$. Moreover, $i_Z^*(\mathfrak{C}_{\text{univ}})$ is isomorphic to $\mathfrak{C}_1$ under an isomorphism that lifts $\mathcal E_s$ for all $s\in\mathbb N^*$. Thus part (a) holds. 

Part (b) follows from part (a) and the fact that $\psi_{M_1}$ is the unique principal quasi-polarization of $(M_1,F^1,\Phi,\nabla_1)$ which modulo $\mathfrak{I}_1$ is $\psi_M$.\endproof
 
\subsubsection{Variant of Subsubsection B4.1 and Theorem \ref{T10}} 

Let $d\in\mathbb N$ be the rank of $\text{Lie}(U)=\text{Lie}(\mathcal G_{B(k)})\cap\text{Hom}(F^1,F^0)$. Let $S:=W(k)[[x_1,\ldots,x_d]]$ and $\mathfrak{I}_0:=(x_1,\ldots,x_d)$ be its ideal. We consider an arbitrary closed embedding $\Spec(S)\hookrightarrow\Spec(R)$ such that the following two properties hold:

\medskip\noindent
{\bf (i)} at the level of $W(k)$-algebras, the ideal $\mathfrak{I}$ of $R$ maps to the ideal $\mathfrak{I}_0$ of $S$;

\smallskip\noindent
{\bf (ii)} the pullback $\mathfrak{D}_{\text{univ}}$ of $\mathfrak{C}_{\text{univ}}$ via the closed embedding $\Spec(S)\hookrightarrow\Spec(R)$, has a Kodaira--Spencer map which is injective and whose image equals to the direct summand $\text{Lie}(U)\otimes_{W(k)} S$ of $\text{Lie}(U_{\text{bigg}})\otimes_{W(k)} S\simeq\text{Hom}(F^1,M/F^1)\otimes_{W(k)} S$. 

\medskip
The proof of Theorem \ref{T10} applies to give us that there exists a morphism $j_Z:Z\rightarrow\Spec(S)$ of $W(k)$-schemes such that $\mathfrak{C}_1$ is isomorphic to $j_Z^*(\mathfrak{D}_{\text{univ}})$  under an isomorphism which modulo $\mathfrak{I}_1$ becomes the identity automorphism $1_M$ of $M$. As the Kodaira--Spencer map of $\mathfrak{D}_{\text{univ}}$ is injective, the morphism $j_Z$ is unique. In simpler words, we can choose $i_Z:Z\rightarrow\Spec(R)$ to factor through the closed embedding $\Spec(S)\hookrightarrow\Spec(R)$ and the resulting factorization is our unique morphism $j_Z:Z\rightarrow\Spec(S)$. 

In this paragraph we assume that $\mathcal G$ is smooth over $W(k)$. This assumption implies that the normalizer $\mathcal P=\mathcal P_1$ of $F^1$ in $\mathcal G$ is smooth over $W(k)$ and the product morphism $U\times_{\Spec(W(k))} \mathcal P\rightarrow \mathcal G$ is an open embedding, cf. [10, Lem. 2.1.5 and Prop. 2.1.8 (3)]. Thus we can view $g_{\text{univ}}\in \mathcal G(R)$ as an $R$-valued point of $U\times_{\Spec(W(k))} \mathcal P\rightarrow \mathcal G$ as well as of the quotient $[U\times_{\Spec(W(k))} \mathcal P]/\mathcal P=U$. So the closed embedding $\Spec(S)\hookrightarrow\Spec(R)$ can be any closed embedding with the property that 
the morphism $\Spec(S)\rightarrow [U\times_{\Spec(W(k))} \mathcal P]/\mathcal P=U$ induced by $g_{\text{univ}}$ is formally \'etale.

\subsubsection{On the $p=2$ case}

The following result complements Theorem \ref{T9} for $p=2$. 

\begin{theorem} \label{T11}
We assume that $p=2$ and that one of the following two conditions holds:

\medskip\noindent
{\bf (i)} the group scheme $\mathcal G$ is reductive;

\smallskip\noindent
{\bf (ii)} the $2$-divisible group $D_k$ is ordinary. 

\medskip
{\bf (a)} Then there exists a $2$-divisible group $D^{\prime}$ over $W(k)$ which lifts $D_k$, whose filtered $F$-crystal over $k$ is as well the triple $(M,F^1,\phi)$, and for which there exists an isomorphism $\varrho_{D^{\prime}}:(M,(t_{\alpha})_{\alpha\in\mathcal J})\rightarrow (H^1(D^{\prime})\otimes_{\mathbb Z_2} W(k),(v_{\alpha})_{\alpha\in\mathcal J})$. Here $v_{\alpha}\in\mathcal T(H^1(D^{\prime}))[{1\over 2}]=\mathcal T(H^1(D))[{1\over 2}]$ is the tensor that corresponds to $t_{\alpha}$ via Fontaine comparison theory for either $D^{\prime}$ or $D$ (cf. the canonical identification $H^1(D^{\prime})[{1\over 2}]=H^1(D)[{1\over 2}]$ induced by the $B_{\text{dR}}(W(k))$-linear isomorphism $i_{D^\prime}\circ i_D^{-1}$).

\smallskip
{\bf (b)} Let $\lambda_{D_k}$ be the principal quasi-polarization of $D_k$ which is the pullback of the principal quasi-polarization $\lambda_D$ of $D$. Then we can assume that $D^{\prime}$ and $\varrho_{D^{\prime}}$ are such that there exists a principal quasi-polarization $\lambda_{D^{\prime}}$ of $D^{\prime}$ which lifts $\lambda_{D_k}$ and whose \'etale realization is a perfect, alternating bilinear form $\psi_{H^1(D^{\prime})}$ on $H^1(D^{\prime})$ such that $\varrho_{D^{\prime}}$ is in fact an isomorphism $\varrho_{D^{\prime}}:(M,(t_{\alpha})_{\alpha\in\mathcal J},\psi_M)\rightarrow (H^1(D^{\prime})\otimes_{\mathbb Z_2} W(k),(v_{\alpha})_{\alpha\in\mathcal J},\psi_{H^1(D^{\prime})})$.

\smallskip
{\bf (c)} If (ii) holds, then we moreover assume that $\mathcal G$ is smooth. Then the number of $D^{\prime}$'s (resp. of $(D^{\prime},\lambda_{D^{\prime}})$'s) for which part (a) (resp. (b)) holds is $2^a$, where $a$ is the multiplicity of the Newton polygon slope $-1$ for $(\text{Lie}(\mathcal G)[{1\over 2}],\phi)$. Moreover, if we can take $D^{\prime}=D$, then each other such $D^{\prime}$ is the pullback of the $2$-divisible group $D_R$ of Lemma \ref{L16} via a uniquely determined morphism $\Spec(W(k))\rightarrow \Spec(R)$ that factors through the closed embedding $\Spec(S)\hookrightarrow \Spec(R)$ introduced in Subsubsection B4.2.

\smallskip
{\bf (d)} We assume that (ii) holds and that $\phi(F^1)=2F^1$. Then referring to part (a), as $D^{\prime}$ we can take the canonical lift of $D_k$.
\end{theorem}

\medskip\proof We prove part (a). We consider the direct sum decomposition 
$$(M,\phi)=(M_0,\phi)\oplus (M_{>0},\phi)$$ 
such that $\phi(M_0)=M_0$ and $\phi:M_{>0}\rightarrow M_{>0}$ is topologically nilpotent. 

In this paragraph we check that there exists a cocharacter $\tilde \mu:\mathbb G_{m,W(k)}\rightarrow\mathcal G$ which normalizes the descending Newton polygon slope filtration of $(M,\phi)$ (in particular, it normalizes $M_{>0}$) and which produces naturally a direct sum decomposition $M=\tilde F^1\oplus \tilde F^0$ such that $\tilde F^1/2\tilde F^1=F^1/2F^1$ (for each $i\in\{0,1\}$, every $\beta\in\mathbb G_{m,W(k)}(W(k))$ acts through $\tilde\mu$ on $\tilde F^i$ as the multiplication by $\beta^{-i}$); this implies that we have $\tilde F^1\subset M_{>0}$. If $\mathcal G$ is a reductive group scheme over $W(k)$, then the existence of $\tilde\mu$ is a particular case of [61, Thm. 1.3.1 or Cor. 1.3.2 (a)]. If $D_k$ is ordinary, then we have $\phi(M_{>0})=2M_{>0}$ and we can take $\tilde F^1=M_{>0}$ and $\tilde F^0=M_0$; the resulting cocharacter $\tilde\mu:\mathbb G_{m,W(k)}\rightarrow\pmb{\text{GL}}_M$ fixes each $t_{\alpha}$ (as $\tilde\mu$ is the inverse of the Newton cocharacter of $(M,\phi)$ and as we have $\phi(t_{\alpha})=t_{\alpha}$ for all $\alpha\in\mathcal J$), and therefore it factors through $\mathcal G$ as desired. 

Let $\tilde D=\tilde D_0\oplus \tilde D_{>0}$ be the unique $2$-divisible group over $W(k)$ such that the filtered $F$-crystals of $\tilde D_0$ and $\tilde D_{>0}$ are $(M_0,0,\phi)$ and $(M_{>0},\tilde F^1,\phi)$ (respectively), cf. [64, Prop. 2.2.6] for the uniqueness of $\tilde D_{>0}$. If $D_k$ is ordinary, then $\tilde D$ is the canonical lift of $D_k$. From Theorem \ref{T9} we get the existence of an isomorphism $\varrho_{\tilde D}:(M,(t_{\alpha})_{\alpha\in\mathcal J})\rightarrow (H^1(\tilde D)\otimes_{\mathbb Z_2} W(k),(\tilde v_{\alpha})_{\alpha\in\mathcal J})$, where $\tilde v_{\alpha}\in\mathcal T(H^1(\tilde D))[{1\over 2}]$ corresponds to $t_{\alpha}$ via Fontaine comparison theory for $\tilde D$. Thus, if $\tilde F^1=F^1$, then we can take $D^{\prime}=\tilde D$.

In the general case (thus $\tilde F^1$ could now be different from $F^1$), we will use the deformation theory of Subsection B6 for $\tilde D$ in order to prove that $D^{\prime}$ exists. If $\mathcal G$ is a reductive group scheme, then we have $\mathcal G^{\prime}=\mathcal G$. Let $R$, $\mathfrak{I}$, $M_R$, $\Phi$, $\nabla$ be as in Subsubsection B4.1. Let $\tilde F^1_R:=\tilde F^1\otimes_{W(k)} R$. There exists a unique $2$-divisible group $\tilde D_R$ over $R$ which modulo the ideal $\mathfrak{I}$ is $\tilde D$ and whose filtered $F$-crystal over $R/2R$ is $(M_R,\tilde F^1_R,\Phi,\nabla)$, cf. Lemma \ref{L16} applied to $(\tilde D,\tilde F^1_R)$ instead of $(D,F^1_R)$. Let $\tilde{\mathfrak{C}}_{\text{univ}}:=(M_R,\tilde F^1_R,\Phi,\nabla,(t_{\alpha})_{\alpha\in\mathcal J})$ be the last filtered $F$-crystal endowed with the family $(t_{\alpha})_{\alpha\in\mathcal J}$ of crystalline tensors. Let $z:\Spec(W(k))\rightarrow\Spec(R)$ be the closed embedding defined by the ideal $\mathfrak{I}$ of $R$. We have $z^*(\tilde D_R)=\tilde D$. We emphasize that the pullbacks of $\tilde D_R$ and $D_R$ to $\Spec(R/2R)$ coincide, cf. [3, Thm. 4.1.1]. Thus a closed embedding $\Spec(S)\hookrightarrow \Spec(R)$ chosen as in Subsubsection B4.2 working with $D_R$ works as well for $\tilde D_R$.  

Let $K$ be the field of fractions of $R$. From [64, Subsubsect. 3.4.2 and Lem. 3.4.3] we get that for each $\alpha\in\mathcal J$ there exists an \'etale Tate-cycle $\tilde{\mathcal V}_{\alpha}\in\mathcal T(H^1(\tilde D_K))[{1\over 2}]$ on $\tilde D_K$ which corresponds to $t_{\alpha}$ via Fontaine comparison theory for $\tilde D_R$. If $z_1:\Spec(W(k))\rightarrow\Spec(R)$ is a closed embedding, then the filtered $F$-crystal of $D_1:=z_1^*(\tilde D_R)$ is of the form $(M,F^1_1,\phi)$ for a suitable direct summand $F^1_1$ of $M$ which lifts $F^1/2F^1$ and moreover to each $t_{\alpha}$ corresponds an \'etale Tate-cycle $v_{1,\alpha}\in \mathcal T(H^1(D_1))[{1\over 2}]$ on $D_{1,B(k)}$ in such a way that we have a canonical identification $(H^1(\tilde D_K),(\tilde{\mathcal V}_{\alpha})_{\alpha\in\mathcal J})=(H^1(D_1),(v_{1,\alpha})_{\alpha\in\mathcal J})$ (see proof of [64, Lem. 3.4.3]). 

Thus we have a canonical identification $(H^1(\tilde D),(\tilde v_{\alpha})_{\alpha\in\mathcal J})=(H^1(D_1),(v_{1,\alpha})_{\alpha\in\mathcal J})$. Therefore the existence of $\varrho_{\tilde D}$ implies the existence of an isomorphism 
$$\varrho_{D_1}:(M,(t_{\alpha})_{\alpha\in\mathcal J})\rightarrow\break (H^1(D_1)\otimes_{\mathbb Z_2} W(k),(v_{1,\alpha})_{\alpha\in\mathcal J}).$$ 
Thus to end the proof of part (a) it suffices to show that we can choose $z_1$ such that we have $F^1_1=F^1$ (and then we can take $D^{\prime}=D_1$). Let $v\in 2\text{Lie}(U)$ be the unique element such that for $u:=1_M+v\in \text{Ker}(U(W(k))\rightarrow U(k))$ we have $u(\tilde F^1)=F^1$, cf. Lemma \ref{L12}. By denoting $z_{1,0}:\Spec(k)\hookrightarrow\Spec(R)$ the closed point of $\Spec(R)$, by induction on $s\in\mathbb N^*$ we check that there exists a morphism $z_{1,s}:\Spec(W_s(k))\rightarrow\Spec(R)$ which lifts $z_{1,s-1}$ and such that the Hodge filtration of $z_{1,s}^*(\tilde D_R)$ is the direct summand $F^1/2^sF^1$ of $M/2^sM$. We can take $z_{1,1}:=z_{1,0}$. For $s\ge 2$, assuming that $z_{1,s-1}$ exists, the existence of the lift $z_{1,s}$ of $z_{1,s-1}$ is implied by the property (iii) of Subsubsection B4.1 and the relation $v\in 2\text{Lie}(U)$ (the arguments for these are the same as the ones of the proof of [63, Prop. 6.4.6 (b)] and rely on the fact that our field $k$ is algebraically closed).

To prove part (b), we consider a direct sum decomposition 
$$(M_{>0},\phi)=(M_{(0,1)},\phi)\oplus (M_1,\phi)$$ 
such that $\phi(M_1)=2M_1$ and all Newton polygon slopes of $(M_{(0,1)},\phi)$ belong to $(0,1)\cap\mathbb Q$. As $\tilde\mu$ normalizes the descending Newton polygon slope filtration of $(M,\phi)$, we have $M_1\subset \tilde F^1$. Thus 
$$(M_{>0},\tilde F^1,\phi)=(M_{(0,1)},M_{(0,1)}\cap \tilde F^1,\phi)\oplus (M_1,M_1,\phi)$$ 
and therefore we have a uniquely determined direct sum decomposition $\tilde D_{>0}=\tilde D_{(0,1)}\oplus \tilde D_1$: the filtered $F$-crystals over $k$ of $\tilde D_{(0,1)}$ and $\tilde D_1$ are $(M_{(0,1)},M_{(0,1)}\cap \tilde F^1,\phi)$ and $(M_1,M_1,\phi)$ (respectively). 

As $\tilde\mu$ factors through $\mathcal G$ and as $\mathcal G$ normalizes the $W(k)$-span of $\psi_M$, $\tilde F^1$ is a maximal isotropic direct summand of $M$ with respect to $\psi_M$. Due to this and the uniqueness properties of $\tilde D=\tilde D_0\oplus \tilde D_{(0,1)}\oplus \tilde D_1$, there exists a unique principal quasi-polarization $\lambda_{\tilde D}$ of $\tilde D$ which lifts $\lambda_{D_k}$. The \'etale realization of $\lambda_{\tilde D}$ is a perfect, alternating bilinear form $\psi_{H^1(\tilde D)}$ on $H^1(\tilde D)$. We choose $\varrho_{\tilde D}$ such that we have an isomorphism  $\varrho_{\tilde D}:(M,(t_{\alpha})_{\alpha\in\mathcal J},\psi_M)\rightarrow (H^1(\tilde D)\otimes_{\mathbb Z_2} W(k),(\tilde v_{\alpha})_{\alpha\in\mathcal J},\psi_{H^1(\tilde D)})$, cf. Theorem \ref{T9}. Let $\lambda_{\tilde D_R}$ be the unique principal quasi-polarization of $\tilde D_R$ whose reduction modulo the ideal $\mathfrak{I}$ is $\lambda_{\tilde D}$ and whose crystalline realization is the perfect, alternating bilinear form $\psi_{M_R}$ on $M_R$, cf. Lemma \ref{L17} applied to $(\tilde D,\lambda_{\tilde D},\tilde F^1_R)$ instead of $(D,\lambda_{D},F^1_R)$.

The remaining part of the proof of part (b) is the same as of part (a). Briefly, it goes as follows. If $\tilde F^1=F^1$, then we take $(D^{\prime},\lambda_{D^{\prime}})=(\tilde D,\lambda_{\tilde D})$. If $\tilde F^1\neq F^1$, then we have to consider the filtered principally quasi-polarized $F$-crystal $(M,F^1_1,\phi,\psi_M)$ of $(D_1,\lambda_{D_1}):=z_1^*(\tilde D_R,\lambda_{\tilde D_R})$ and the \'etale realizations $\psi_{H^1(\tilde D_K)}$ and $\psi_{H^1(D_1)}$ of $(\lambda_{\tilde D_R})_K$ and $\lambda_{D_1}$ (respectively); as above one gets a canonical identification $(H^1(\tilde D),(\tilde v_{\alpha})_{\alpha\in\mathcal J},\psi_{H^1(\tilde D)})=(H^1(D_1),(v_{1,\alpha})_{\alpha\in\mathcal J},\psi_{H^1(D_1)})$. If $z_1:\Spec(W(k))\rightarrow\Spec(R)$ is such that $F^1_1=F^1$, then by taking $(D^{\prime},\lambda_{D^{\prime}})=(D_1,\lambda_{D_1})$ we get that part (b) holds.

To prove part (c), based on the proof of part (b) it suffices to consider only the non-principally quasi-polarized case. To ease notation we can assume that $D$ is one of the $D^{\prime}$'s, cf. part (a). Thus there exists an isomorphism $\varrho_{D}:(M,(t_{\alpha})_{\alpha\in\mathcal J})\rightarrow (H^1(D)\otimes_{\mathbb Z_2} W(k),(v_{\alpha})_{\alpha\in\mathcal J})$. We will consider two cases, the first one being only a particular case of the second (general) one.

\medskip
{\bf Case 1.} We assume that $F^1=\tilde F^1$ and $D=\tilde D$ and thus also $D_R=\tilde D_R$. To the direct sum decomposition $D=\tilde D=\tilde D_0\oplus \tilde D_{(0,1)}\oplus \tilde D_1$, corresponds a direct sum decomposition $H^1(D)=H^1(D)_0\oplus H^1(D)_{(0,1)}\oplus H^1(D)_1$. If $D^{\prime}$ is a $2$-divisible group for which part (a) holds, then we have short exact sequences $0\rightarrow \tilde D_1\rightarrow D^{\prime}\rightarrow \tilde D_{(0,1)}\oplus \tilde D_0\rightarrow 0$ and $0\rightarrow \tilde D_1\oplus\tilde D_{(0,1)}\rightarrow D^{\prime}\rightarrow\tilde D_0\rightarrow 0$ and $H^1(D^{\prime})$ is a $\mathbb Z_2$-submodule of ${1\over 2}H^1(D)$ that contains $2H^1(D)$ (as one can easily check based on [64, Prop. 2.2.6] and the proof of [64, Lem. 2.2.5]). We get the existence of an element $c\in {1\over 2}\text{Hom}(H^1(D)_1,H^1(D)_0)$ such that $H^1(D^{\prime})=(1_M+c)(H^1(D))$; it is uniquely determined modulo $\text{Hom}(H^1(D)_1,H^1(D)_0)$.
But as there exists an isomorphism $\varrho_{D^{\prime}}:(M,(t_{\alpha})_{\alpha\in\mathcal J})\rightarrow (H^1(D^{\prime})\otimes_{\mathbb Z_2} W(k),(v_{\alpha})_{\alpha\in\mathcal J})$, there exists $g\in \mathcal G^{\acute{et}}(B(k))$ such that $H^1(D^{\prime})\otimes_{\mathbb Z_2} W(k)=g(H^1(D)\otimes_{\mathbb Z_2} W(k))$. 

We claim that we can assume that we have 
$$c\in {1\over 2}[\text{Hom}(H^1(D)_1,H^1(D)_0)\cap \text{Lie}(\mathcal G_{\mathbb Q_2}^{\acute{et}})],$$ 
i.e., the image of $c$ in the quotient group ${1\over 2}\text{Hom}(H^1(D)_1,H^1(D)_0)/\text{Hom}(H^1(D)_1,H^1(D)_0)$ belongs to the following subgroup 
$${1\over 2}[\text{Hom}(H^1(D)_1,H^1(D)_0)\cap \text{Lie}(\mathcal G_{\mathbb Q_2}^{\acute{et}})]/[\text{Hom}(H^1(D)_1,H^1(D)_0)\cap \text{Lie}(\mathcal G_{\mathbb Q_2}^{\acute{et}})].$$ This is only a variant of the Lemma \ref{L13} over $\mathbb Z_2$ instead of $W(k)$ which gets reduced to the Lemma \ref{L13} as follows. We can assume that $\varrho_D$ maps $M_0$, $M_{(0,1)}$, and $M_1$ onto $H^1(D)_0\otimes_{\mathbb Z_2} W(k)$, $H^1(D)_{(0,1)}\otimes_{\mathbb Z_2} W(k)$, and $H^1(D)_1\otimes_{\mathbb Z_2} W(k)$ (respectively), cf. Theorem \ref{T9}. Thus we can assume that the element 
$$\varrho_D g\varrho_D^{-1}\pmb{\text{GL}}_M(W(k))=\varrho_D (1_M+c)\varrho_D^{-1}\pmb{\text{GL}}_M(W(k))\in\pmb{\text{GL}}_M(B(k))/\pmb{\text{GL}}_M(W(k))$$ 
belongs to the intersection of 
$$[1_M+{1\over 2}\text{Lie}(U_{\text{bigg}})]/U_{\text{bigg}}(W(k))\cap \mathcal P_0(B(k))/\mathcal P_0(W(k))$$ taken inside $\pmb{\text{GL}}_M(B(k))/ \pmb{\text{GL}}_M(W(k))$ and thus it is an element of $[1_M+{1\over 2}\text{Lie}(U)]/U(W(k))$ (cf. Lemma \ref{L13}). We note that if (ii) holds, then based on Theorem \ref{T9} we can assume that $\varrho_Dg\varrho_D^{-1}$ fixes $W(k)$-bases of $M_0$ and $M/M_0$ formed by elements fixed by $\phi$ and $p^{-1}\phi$ (respectively) and thus in fact in the above intersection we can replace $\mathcal P_0(B(k))/\mathcal P_0(W(k))$ by $U(B(k))/U(W(k))$. 

This implies that there exists an element $c_{\text{crys}}\in {1\over 2}[\text{Hom}(M_1,M_0)\cap\text{Lie}(\mathcal G)]$ such that $\varrho_D c\varrho_D^{-1}-c_{\text{crys}}\in\text{Hom}(M^1,M^0)$. Thus 
$$c-\varrho_D^{-1}c_{\text{crys}}\varrho_D\in \text{Hom}(H^1(D)_1,H^1(D)_0)\otimes_{\mathbb Z_2} W(k)$$ and moreover 
$$\varrho_D^{-1}c_{\text{crys}}\varrho_D\in {1\over 2}[\text{Hom}(H^1(D)_1,H^1(D)_0)\cap \text{Lie}(G_{\mathbb Q_2}^{\acute et})]\otimes_{\mathbb Z_2} W(k).$$ Therefore by replacing $c$ with $\varrho_D^{-1}c_{\text{crys}}\varrho_D$, we get that the claim follows. 

The group ${1\over 2}[\text{Hom}(H^1(D)_1,H^1(D)_0)\cap \text{Lie}(\mathcal G_{\mathbb Q_2}^{\acute{et}})]/[\text{Hom}(H^1(D)_1,H^1(D)_0)\cap \text{Lie}(\mathcal G_{\mathbb Q_2}^{\acute{et}})]$ has order $2^a$. We conclude that the number of $\mathbb Z_2$-lattices $H^1(D^{\prime})$ of $H^1(D)[{1\over 2}]$ such that the above properties hold (equivalently, the number of $D^{\prime}$'s as in part (a)), is precisely $2^a$. The fact that all of them are pullbacks of $D_R$ via $W(k)$-valued points of $\Spec(R)$ follows from the fact that there exists a closed embedding $\Spec(R_1)\hookrightarrow\Spec(R)$ defined by an ideal of $R$ contained in the ideal $(x_1,\ldots,x_l)$ and with $R_1=W(k)[[x_1,\ldots,x_a]]$, such that the restriction $\tilde D_{R_1}$ of $\tilde D_R=D_R$ to $\Spec(R_1)$ is a direct sum $\tilde D_{R_1}=\tilde D_{(0,1),R_1}\oplus \tilde D_{0,1,R_1}$, where $\tilde D_{0,1,R_1}$ sits in a short exact sequence $0\rightarrow \tilde D_{1,R_1}\rightarrow\tilde D_{0,1,R_1}\rightarrow \tilde D_{0,R_1}\rightarrow 0$ which is a versal deformation of $\tilde D_1\oplus \tilde D_0$ and which endows $\Spf(R_1)$ with the structure of a formal subtorus of dimension $a$ of the formal torus over $\Spf(W(k))$ of deformations of the ordinary $2$-divisible group $\tilde D_{1,k}\oplus \tilde D_{0,k}$ over $k$ (here versal is used in the sense that the Kodaira--Spencer map is injective and has an image which is a direct summand of its codomain). More precisely, if $U_1$ is the smooth, connected, closed subgroup scheme of $U$ whose Lie algebra is $\text{Hom}(M_1,M_0)\cap\text{Lie}(U)=\text{Hom}(M_1,M_0)\cap\text{Lie}(\mathcal G)$ (and thus has rank $a$), then the filtered $F$-crystal of $\tilde D_{R_1}$ endowed with tensors is 
$\tilde{\mathfrak{D}}_1:=(M\otimes_{W(k)} R_1,\tilde F^1\otimes_{W(k)} R_1,\Phi_1,\nabla,(t_{\alpha})_{\alpha\in\mathcal J})$, where $\Phi_1=u_1(\phi\otimes\Phi_{R_1})$ with $\Phi_{R_1}$ as in Subsubsection B4.1 for $m:=a$ and with $u_1\in U_1(R_1)$ a universal element which identifies $R_1$ with the completion of the local ring of $U_1$ at the identity element of $U_1(k)$. This is so  as $\tilde{\mathfrak{D}}_1$ is the pullback of $\tilde{\mathfrak{C}}_{\text{univ}}$ via a $W(k)$-morphism $\Spec(R_1)\rightarrow\Spec(R)$ which is a closed embedding and which at the level of rings maps the ideal $(x_1,\ldots,x_l)$ of $R$ to the ideal $(x_1,\ldots,x_a)$ of $R_1$, cf. Theorem \ref{T10} (a) and the fact that $\tilde{\mathfrak{D}}_1$ is versal. Each $\Spf(W(k))$-valued point of the formal torus $\Spf(R_1)$ which is of order $1$ or $2$ corresponds uniquely to a $D^{\prime}$ as in part (a) and therefore indeed we have precisely $2^a$ such $D^{\prime}$'s as in part (a) and all of them are pullbacks of $D_R$ via $W(k)$-valued points of $\Spec(R)$. Thus part (c) holds if $D=\tilde D$. From the uniqueness part of Subsubsection B4.2 we get that we can assume that $\Spec(R_1)$ is as well a closed subscheme of the closed subscheme $\Spec(S)$ of $\Spec(R)$ chosen in Subsubsection B4.2; therefore all $D^{\prime}$'s as in part (a) are pullbacks of $D_R$ via uniquely determined $W(k)$-valued points of $\Spec(R)$ that factor through $\Spec(S)$. 

\medskip
{\bf Case 2.} We check that part (c) holds in the general case (i.e., without assuming that $F^1=\tilde F^1$ and $D=\tilde D$). Let $\tilde D_S$ be the pullback of $\tilde D_R$ constructed above via the closed embedding $\Spec(S)\rightarrow \Spec(R)$ of the Subsubsection B4.2. Let $\tilde{\mathfrak{D}}_{\text{univ}}$ be the pullback to $S/2S$ of $\tilde{\mathfrak{C}}_{\text{univ}}$. As in the proof of part (a) we argue that there exists a morphism $\tilde z^{\prime}:\Spec(W(k))\rightarrow\Spec(R)$ such that the Hodge filtration of $M$ defined by $\tilde D^{\prime}:=(\tilde z^{\prime})^*(\tilde D_R)$ is $\tilde F^1$ and there exists an isomorphism $\varrho_{\tilde D^{\prime}}:(M,(t_{\alpha})_{\alpha\in\mathcal J})\rightarrow (H^1(\tilde D^{\prime})\otimes_{\mathbb Z_2} W(k),(\tilde v_{\alpha})_{\alpha\in\mathcal J})$. Let $\tilde{\mathfrak{I}}^{\prime}$ be the ideal of $R$ that defines $\tilde z^{\prime}$. Let $y_1,\ldots,y_l$ be regular parameters of $R$ such that we have an identity $\tilde{\mathfrak{I}}^{\prime}=(y_1,\ldots,y_l)$ between ideals of $R$. Let $\tilde\Phi_{R,1}$ be the Frobenius lift of $\Spec(R)$ which is compatible with $\sigma$ and takes each $y_i$ to $y_i^2$. Based on Case 1, we can assume that the morphism $\tilde z^{\prime}:\Spec(W(k))\rightarrow\Spec(R)$ factors through the closed embedding $\Spec(S)\hookrightarrow \Spec(R)$. Let $\tilde z^{\prime}_S:\Spec(W(k)(\rightarrow \Spec(S)$ be the resulting factorization. 

From Theorem \ref{T10} (a) we get that $D_R$ and $\mathfrak{C}_{\text{univ}}$ are the pullbacks of $\tilde D_S$ and $\tilde{\mathfrak{D}}_{\text{univ}}$ (respectively) via a morphism $h:\Spec(R)\rightarrow \Spec(S)$ that satisfies the identity $h\circ \tilde z^{\prime}=\tilde  z^{\prime}_S$ (for this part we have to consider new Frobenius lifts of $R$ and $S$; like for $R$ we would have to replace $\Phi_R$ by $\tilde\Phi_{R,1}$). Due to the uniqueness part of Subsubsection B4.2 and the identity $h\circ \tilde z^{\prime}=\tilde z^{\prime}_S$, the closed embedding $\Spec(S)\hookrightarrow \Spec(R)$ is a section of $h:\Spec(R)\rightarrow \Spec(S)$.  

Due to the existence of $h$, to prove part (c) in the general case it suffices to show that there exist exactly $2^a$ morphisms $z_1:\Spec(W(k))\rightarrow \Spec(S)$ such that the Hodge filtration of $M$ defined by $z_1^*(\tilde D_S)$ is $F^1$. Fixing such a morphism $z_{1,0}$ (it exists, cf. proof of part (a)), any other such morphism $z_1$, induces a unique isomorphism $h_1:\Spec(S)\rightarrow\Spec(S)$ with the properties that $\tilde D_S=h_1^*(\tilde D_S)$ and we have $h_1\circ z_1=z_{1,0}$. But the number of isomorphisms $h_2:\Spec(S)\rightarrow\Spec(S)$ with the property that $\tilde D_S=h_2^*(\tilde D_S)$ is uniquely determined by the property that under it the ideal $\mathfrak{I}_0$ of $S$ that defines $\tilde D$ is mapped to one of the $2^a$ ideals of $S$ under which one gets a $2$-divisible group over $W(k)$ whose Hodge filtration is $\tilde F^1$ (cf. Case 1 applied to $\tilde D$). Thus we have $2^a$ such $h_2$'s and $z_1$'s and therefore part (c) holds in the general case.

Part (d) follows from Theorem \ref{T9}.\endproof 

\subsection{On abelian schemes} 

\smallskip
We assume that $D$ is the Barsotti--Tate group of an abelian scheme $A$ over $W(k)$. It is known that we have two canonical and functorial identifications: 

\medskip\noindent
{\bf (i)} $H^1_{\text{dR}}(A/W(k))=M$ of $W(k)$-modules (see [2, Ch. V, Subsect. 2.3] and [4, Prop. 2.5.8]);

\smallskip\noindent
{\bf (ii)} $H^1(D)=H^1_{\acute et}(A_{\overline{B(k)}},\mathbb Z_p)$ of $\text{Gal}(B(k))$-modules. 

\medskip\noindent
The crystalline conjecture (see [19]) provides a $B_{\text{crys}}(W(k))$-linear isomorphism
$$i_A:H^1_{\text{dR}}(A/W(k))\otimes_{W(k)} B_{\text{crys}}(W(k))\rightarrow H^1_{\acute et}(A_{\overline{B(k)}},\mathbb Z_p)\otimes_{\mathbb Z_p} B_{\text{crys}}(W(k))$$
that is compatible with the tensor product filtrations, with the $\text{Gal}(B(k))$-actions, and with the Frobenius endomorphisms. See [54, Subsubsect. 5.2.15] for a proof of the following property (strictly speaking, the paragraphs before loc. cit. work with a prime $p\ge 3$ but the arguments of loc. cit. work for all primes):

\medskip\noindent
{\bf (iii)} under the identifications of (i) and (ii), we have $i_A=i_D^+\otimes 1_{B_{\text{crys}}(W(k))}$.

\subsection{On Hodge cocharacters} 

\smallskip
In this subsection we assume that we have a monomorphism $W(k)\hookrightarrow \mathbb C$ and that $D$ is the Barsotti--Tate group of an abelian scheme $A$ over $W(k)$. 

We recall that we have canonical identifications
$$M\otimes_{W(k)} \mathbb C=H^1_{\text{dR}}(A/W(k))\otimes_{W(k)} \mathbb C=H^1_{\text{dR}}(A_{\mathbb C}/\mathbb C)=F^{1,0}\oplus F^{0,1},\leqno (6)$$
where the last identity is the usual Hodge decomposition. Under (6) we can identify
$$F^1\otimes_{W(k)} \mathbb C=F^{1,0}.$$
Let $A_{\mathbb C}^{\text{an}}$  be the complex manifold associated to $A_{\mathbb C}$. Let $W:=H_1(A_{\mathbb C}^{\text{an}},\mathbb Q)$ be the first Betti homology group of $A_{\mathbb C}^{\text{an}}$ with rational coefficients. Let $W^{\vee}:=\text{Hom}(W,\mathbb Q)$. We identify naturally $W^{\vee}\otimes_{\mathbb Q} \mathbb C$ with the first Betti cohomology group $H^1(A_{\mathbb C}^{\text{an}},\mathbb C)$ and thus also with $H^1_{\text{dR}}(A_{\mathbb C}/\mathbb C)=M\otimes_{W(k)} \mathbb C$. Let $\mu_A:\mathbb G_{m,\mathbb C}\rightarrow \pmb{\text{GL}}_{W^{\vee}\otimes_{\mathbb Q} \mathbb C}$ be the Hodge cocharacter that fixes $F^{0,1}$ and that acts on $F^{1,0}$ via the weight $-1$.  

\begin{lemma}\label{L18} 
Let the cocharacter $\mu:\mathbb G_{m,W(k)}\rightarrow\mathcal G$ be as in Subsection B2. We assume that for every $\alpha\in\mathcal J$ the tensor $t_{\alpha}\in \mathcal T(M)[\frac{1}{p}]=\mathcal T(H^1_{\text{dR}}(A/W(k)))[\frac{1}{p}]$ is the de Rham component of a Hodge cycle on $A_{B(k)}$. We also assume that $\mathcal G_{B(k)}$ is a reductive group. Then the cocharacter $\mu_A:\mathbb G_{m,\mathbb C}\rightarrow \pmb{\text{GL}}_{M\otimes_{W(k)} \mathbb C}$ factors through $\mathcal G_{\mathbb C}$ and this factorization $\mu_A:\mathbb G_{m,\mathbb C}\rightarrow \mathcal G_{\mathbb C}$ is $\mathcal G(\mathbb C)$-conjugate to $\mu_{\mathbb C}$. Thus, if $\mathcal G_{B(k)}$ is a torus, then we have $\mu_A=\mu_{\mathbb C}$.
\end{lemma}
\begin{proof}
Let $v^B_{\alpha}\in \mathcal T(W^{\vee})$ be the Betti realization of $t_{\alpha}$; it is fixed by $\mu_A$. The identity $W^{\vee}\otimes_{\mathbb Q} \mathbb C=M\otimes_{W(k)} \mathbb C$ produces an identity $\mathcal T(W^{\vee}\otimes_{\mathbb Q} \mathbb C)=\mathcal T(M\otimes_{W(k)} \mathbb C)$ under which the tensors $t_{\alpha}$ and $v_{\alpha}^B$ are as well identified. Thus the cocharacter $\mu_A:\mathbb G_{m,\mathbb C}\rightarrow \pmb{\text{GL}}_{W^{\vee}\otimes_{\mathbb Q} \mathbb C}$ fixes $t_{\alpha}$ for all $\alpha\in\mathcal J$ and therefore it factors through $\mathcal G_{\mathbb C}$. Let $\mathcal P_{\mathbb C}$ be the parabolic subgroup of $\mathcal G_{\mathbb C}$ that normalizes  $F^1\otimes_{W(k)} \mathbb C=F^{1,0}$. Both the cocharacters $\mu_A:\mathbb G_{m,\mathbb C}\rightarrow \pmb{\text{GL}}_{M\otimes_{W(k)} \mathbb C}$ and $\mu_{\mathbb C}$ factor through $\mathcal P_{\mathbb C}$ and thus a $\mathcal P_{\mathbb C}(\mathbb C)$-conjugate $\mu^{\prime}_{\mathbb C}$  of $\mu_{\mathbb C}$ commutes with $\mu_A$. As the commuting cocharacters $\mu^{\prime}_{\mathbb C}$ and $\mu_A$ of $\mathcal P_{\mathbb C}$ act on $F^1\otimes_{W(k)} \mathbb C=F^{1,0}$ and on $M\otimes_{W(k)} \mathbb C/(F^1\otimes_{W(k)} \mathbb C)=H^1_{\text{dR}}(A_{\mathbb C}/\mathbb C)/F^{1,0}$ in the same way, we have $\mu^{\prime}_{\mathbb C}=\mu_A$. Thus the cocharacters $\mu_{\mathbb C}$ and $\mu_A$ are $\mathcal P_{\mathbb C}(\mathbb C)$-conjugate and therefore they are also $\mathcal G(\mathbb C)$-conjugate.\end{proof}

\end{appendix}

\end{document}